\documentclass[a4paper,11pt]{article}
\usepackage{amsfonts,amssymb,amsbsy,amsmath,epsf,eepic,graphics,graphicx}
\usepackage{color}
\usepackage{ mathrsfs }
\usepackage{bbm}
\usepackage{enumerate}
\usepackage[english]{babel}
\usepackage[hidelinks]{hyperref}
 \makeatletter
 \@addtoreset{enunciato}{section}
 \makeatother

 \newcounter{enunciato}[section]

 \newtheorem{ittheorem}{Theorem}
 \newtheorem{ittheoremm}{Result}
 \newtheorem{itlemma}{Lemma}
 \newtheorem{itproposition}{Proposition}
 \newtheorem{itdefinition}{Definition}
 \newtheorem{itconjecture}{Conjecture}
 \newtheorem{itassumption}{Assumption}
 \newtheorem{itcorollary}{Corollary}

\newenvironment{remark}{\addtocounter{enunciato}{1} \noindent \textbf{Remark \thesection.\arabic{enunciato}}\,\,}
{\smallskip }

\newenvironment{theorem}{\addtocounter{enunciato}{1}
\begin{ittheorem}}{\end{ittheorem}}

\newenvironment{theoremm}{\addtocounter{enunciato}{1}
\begin{ittheoremm}}{\end{ittheoremm}}

 \newenvironment{lemma}{\addtocounter{enunciato}{1}
 \begin{itlemma}}{\end{itlemma}}

 \newenvironment{proposition}{\addtocounter{enunciato}{1}
 \begin{itproposition}}{\end{itproposition}}

 \newenvironment{definition}{\addtocounter{enunciato}{1}
 \begin{itdefinition}}{\end{itdefinition}}

 \newenvironment{conjecture}{\addtocounter{enunciato}{1}
 \begin{itconjecture}}{\end{itconjecture}}

\newenvironment{assumption}{\addtocounter{enunciato}{1}
\begin{itassumption}}{\end{itassumption}}

\newenvironment{corollary}{\addtocounter{enunciato}{1}
\begin{itcorollary}}{\end{itcorollary}}

 \newcommand{\be}[1]{\begin{equation}\label{#1}}
 \newcommand{\ee}{\end{equation}}
 
 \newcommand{\ba}[1]{\begin{align}\label{#1}}
 \newcommand{\ea}{\end{align}}

 \newcommand{\bea}[1]{\begin{eqnarray}\label{#1}}
 \newcommand{\eea}{\end{eqnarray}}

 \newcommand{\bl}[1]{\begin{lemma}\label{#1}}
 \newcommand{\el}{\end{lemma}}

 \newcommand{\bt}[1]{\begin{theorem}\label{#1}}
 \newcommand{\et}{\end{theorem}}

 \newcommand{\bd}[1]{\begin{definition}\label{#1}}
 \newcommand{\ed}{\end{definition}}

 \newcommand{\bp}[1]{\begin{proposition}\label{#1}}
 \newcommand{\ep}{\end{proposition}}

 \newcommand{\bc}[1]{\begin{corollary}\label{#1}}
 \newcommand{\ec}{\end{corollary}}

 \newcommand{\bcj}[1]{\begin{conjecture}\label{#1}}
 \newcommand{\ecj}{\end{conjecture}}

\newcommand{\bass}[1]{\begin{assumption}\label{#1}}
\newcommand{\eass}{\end{assumption}}

 \newcommand{\bpr}{\begin{proof}}
 \newcommand{\epr}{\end{proof}}

 \newcommand{\bprl}[1]{\begin{proofof}{\it\ref{#1}}.\,\,}
 \newcommand{\eprl}{\end{proofof}}
 
  \newcommand{\bprp}[1]{\begin{proofofp}{\it\ref{#1}}.\,\,}
 \newcommand{\eprp}{\end{proofofp}}

\newcommand{\br}[1]{\begin{remark}\label{#1}}
\newcommand{\er}{\end{remark}}

 \newcommand{\bi}{\begin{itemize}}
 \newcommand{\ei}{\end{itemize}}

 \newcommand{\ben}{\begin{enumerate}}
 \newcommand{\een}{\end{enumerate}}

\newenvironment{proof}{\noindent {\em Proof}.\,\,}
{\hspace*{\fill}$\square$\medskip}
\newenvironment{proofof}{\noindent {\em Proof of Lemma\,\,}}
{\hspace*{\fill}$\square$\medskip}
\newenvironment{proofofp}{\noindent {\em Proof of Proposition\,\,}}
{\hspace*{\fill}$\halmos$\medskip}
\newcommand{\halmos}{\rule{1ex}{1.4ex}}

\newcommand\restr[2]{{
		\left.\kern-\nulldelimiterspace 
		#1 
		\vphantom{\big|} 
		\right|_{#2} 
}}

\newcommand\smallrestr[2]{{
		\left.\kern-\nulldelimiterspace 
		#1 
		\vphantom{|} 
		\right|_{#2} 
}}

\def \R {{\mathbb R}}
\def \N {{\mathbb N}}

\def \Ra {\longrightarrow}
\def \ra {\rightarrow}
\def \to {\mapsto}
\def \hra {\rightharpoonup}

\def \barr{\begin{array}}
\def \earr {\end{array}}

\def\L {\mathrm{L}}
\def\Id {\mathrm{Id}}
\def\Ind {\mathbbm{1}}

\DeclareMathOperator*{\argmin}{\arg\!\min}

\newcommand{\Eb}[1]{
	\mathbb{E}\left[#1\right]
}
\newcommand{\Pb}[1]{
	\mathbb{P}\left[#1\right]
}
\def \P {{\mathbb P}}

\def \cV {{\mathcal V}}
\def \cL {{\mathcal L}}
\def \cH {{\mathcal H}}
\def \cT {{\mathcal T}}
\def \cS {{\mathcal S}}
\def \cC {{\mathcal C}}
\def \cD {{\mathcal D}}
\def \cP {{\mathcal P}}
\def \cJ {{\mathscr J}}

\def \cM {{\mathcal M}}
\def \cW {{\mathcal W}}
\def \cA {{\mathcal A}}
\def \cN {{\mathcal N}}
\def \cK {{\mathcal K}}
\def \cF {{\mathcal F}}

\def\b {\beta}

\def\s {\sigma}
\def\t {\theta}

\def\e {\varepsilon}

\newcommand{\Pw}{\mathcal{P}_2(\R)}
\newcommand{\PwN}{\mathcal{P}_2(\R^N)}
\newcommand{\Pwn}{\mathcal{P}_2(\R^n)}
\newcommand{\Pwa}{\mathcal{P}_2^a(\R)}
\newcommand{\PL}{\mathcal{P}_2^{\mathrm{L}}(\Td \times \R)}
\newcommand{\PLone}{\mathcal{P}_2^{\mathrm{L}}(\Tone \times \R)}
\newcommand{\PLo}{\mathcal{P}_2^{\mathrm{L}}(\Tone \times \R)}
\newcommand{\PLa}{\mathcal{P}_2^{\mathrm{L},a}(\Td \times \R)}
\newcommand{\ML}{\cM_1^\mathrm{L}(\TdR)}
\newcommand{\MLo}{\cM_1^\mathrm{L}(\ToR)}

\def \Td {{\mathbb T}^d}
\def \Tone {{\mathbb T}}
\def \To {{\mathbb T}}
\def \TdR {{\mathbb T}^d \times \R}
\def \ToR {{\mathbb T} \times \R}
\newcommand{\pb}{\textbf{p}}
\newcommand{\AC}[1]{\cA\cC(#1)}

\newcommand{\ACLT}{\cA\cC((0,T);\PL)}

\newcommand{\ACwT}{\cA\cC((0,T);\Pw)}

\newcommand{\TanL}[1]{\mathrm{Tan}_{#1}\PL}
\newcommand{\Tanw}[1]{\mathrm{Tan}_{#1}\Pw}
\newcommand\Opt[1]{\mathrm{Opt}^\mathrm{L}(#1)}
\newcommand\Coupl[1]{\mathrm{Cpl}^\mathrm{L}(#1)}
\newcommand\OptW[1]{\mathrm{Opt}(#1)}
\newcommand\CouplW[1]{\mathrm{Cpl}(#1)}
\newcommand{\LebTd}{\mathrm{Leb}_{\Td}}
\newcommand{\Leb}[1]{\mathrm{Leb}_{#1}}
\newcommand{\WL}{\mathrm{W}^\mathrm{L}}

\newcommand{\muc}{(\mu_t)_t}
\def\T {\mathrm{T}}
\def \bW {{\mathbb W}}

\definecolor{darkgreen}{rgb}{0,.6,0}
\definecolor{darkagenta}{rgb}{.5,0,.5}
\definecolor{darkred}{rgb}{1,0,0}
\definecolor{darkblue}{rgb}{0,0,.4}
\definecolor{black}{rgb}{0,0,0}
\definecolor{gray}{rgb}{.4,.4,.4}
\definecolor{white}{rgb}{0.99,0.99,0.99}
\definecolor{geel1}{rgb}{1,1,0.3}

\numberwithin{equation}{section}
\begin{document}

\bibliographystyle{abbrv}


\title{%
Gradient flow approach to local mean-field spin systems}

\author{
\renewcommand{\thefootnote}{\arabic{footnote}}
K.\ Bashiri\footnotemark[1] , A. Bovier\footnotemark[2]
}

\footnotetext[1]{
Institut f\"ur Angewandte Mathematik,
Rheinische Friedrich-Wilhelms-Universit\"at, 
Endenicher Allee 60, 53115 Bonn, Germany.
Email:
bashiri@iam.uni-bonn.de
}
\footnotetext[2]{
Institut f\"ur Angewandte Mathematik,
Rheinische Friedrich-Wilhelms-Universit\"at, 
Endenicher Allee 60, 53115 Bonn, Germany.
Email:
bovier@uni-bonn.de
}

\date{\today}
\maketitle


\begin{abstract}
\let\thefootnote\svthefootnote{

It is well-known that many diffusion equations can be recast as  Wasserstein gradient flows.
Moreover, in recent years, by modifying the Wasserstein distance appropriately, this technique has been transferred to further evolution equations and systems; see e.g.\ \cite{Maas}, \cite{FathiSSEP}, \cite{ErbarBoltz}.  In this paper we establish such a \emph{gradient flow representation} for evolution equations that depend on a non-evolving parameter. These equations are connected to a local mean-field interacting spin system. We then use this gradient flow representation to prove a \emph{large deviation principle} for the empirical process associated to this system. 
This is done by using the criterion established by Fathi in \cite{FathiLDP}.
Finally, the corresponding \emph{hydrodynamic limit} is  shown by using the  approach initiated in \cite{SandierSer} and \cite{Serfaty} by Sandier and Serfaty.


}
\let\thefootnote\relax\footnotetext{
The research in this paper is partially supported by the German Research Foundation in the Collaborative Research Centre 1060 ``The
Mathematics of Emergent Effects'' and the Bonn International Graduate School in Mathematics (BIGS) in the Hausdorff
Center for Mathematics (HCM).}
\let\thefootnote\svthefootnote
\end{abstract}


\section{Introduction}
\label{S1}

Many classes of diffusion equations can be represented as \emph{gradient flows} in the  space of probability measures equipped with the $ \L^2 $-\emph{Wasserstein distance}. 
This fact was first discovered in the seminal works \cite{otto} and \cite{jko}.
The gradient flow representation entails a lot of useful properties such as contractivity, stability (with respect to gamma-convergence), regularisation estimates and a variational  characterisation (as a minimum of an ``energy-dissipation functional''). 
See \cite{ambgigsav} for a comprehensive treatment of this concept. 
Moreover, in recent years, this gradient flow formalism was translated to other  systems such as \emph{discrete Markov chains} (\cite{Maas} or \cite{Mielke}) and the \emph{Boltzmann equation} (\cite{ErbarBoltz}).

Gradient flow representations  can also be used to study the asymptotic behaviour  of sequences of interacting particle systems.
For example,  \emph{hydrodynamic limit} results were proven in \cite{ErbarSchlichting}, \cite{FathiSSEP} or  \cite{ErbarBoltz}.
These results were obtained via an approach that was  introduced first  in the papers  \cite{SandierSer} and \cite{Serfaty}.
This approach
relies on stability properties of certain functionals that appear in the variational characterisation of the  respective gradient flows.
Furthermore, in the case of sequences of reversible diffusion processes, Fathi shows in \cite{FathiLDP} that  gamma-convergence of those functionals is sufficient to prove a \emph{large deviation principle} for the sequence. 

We establish the following results in the present paper.
\begin{itemize}
	\item 
	In Chapter \ref{S2} we \emph{modify the Wasserstein distance} and establish a \emph{gradient flow formalism} with respect to the resulting metric. 
	Then we show that gradient flows in this modified Wasserstein space correspond to partial differential equations, which depend on a non-evolving parameter.
	In particular, we investigate a special example, which will represent the limiting system in the forthcoming chapters. 
	\item 	
	In Chapter \ref{LDPChap} we use the criterion in \cite{FathiLDP} and the results of Chapter \ref{S2} to prove a \emph{large deviation principle} for a local mean-field interacting spin system, which will be introduced in Section \ref{S1.1.1} and more rigorously in Section  \ref{3.0.1}. 
	\item	
	 In Chapter \ref{HDLChap} we adapt the approach of \cite{SandierSer} and \cite{Serfaty} to prove a \emph{hydrodynamic limit} result for the system in Section \ref{S1.1.1}. 
	Although this result already follows from the large deviation principle from Chapter \ref{LDPChap}, we reprove the statement in order to obtain the hydrodynamic limit result for a slightly larger class of initial values and with respect to the stronger topology  of the Wasserstein distance.
\end{itemize}
  
  The results from Chapter \ref{S2} are new, whereas some of the results from Chapter \ref{LDPChap} and \ref{HDLChap} have already been proven in \cite{HDLAntonDima} and \cite{PatLDP} via different approaches. 
  For instance, the large deviation principle was proven via the approach of the paper \cite{dawgae} and the hydrodynamic limit was proven via the  relative entropy method  (see  \cite{HDLAntonDima}).
The main purpose of the Chapters \ref{LDPChap} and \ref{HDLChap} is to see that gradient flow methods can be used to provide elegant proofs for hydrodynamic limit results and large deviation principles. 
  However, the representation of the rate function in Chapter \ref{LDPChap}  differs from the one in \cite{PatLDP}. We also show here that the rate function admits a unique minimum point, which is not shown in \cite{PatLDP}. Moreover, in Chapter \ref{HDLChap} we also establish the convergence in the stronger topology of the Wasserstein distance.

\section{The model and main results}
  \label{S1.1}

  This chapter is organized as follows.
  In Section \ref{S1.1.1} we introduce the microscopic spin system. 
  In Section \ref{S1.1.2} we define the macroscopic limiting system and show how to modify the Wasserstein distance to obtain a gradient flow representation for this system. In Section \ref{S1.1.3} we give a first formulation of the main results of this paper and sketch the ideas of the proofs.
  In Section \ref{S1.2} we list  notations and definitions we use 
  throughout this paper. 
  To avoid too much terminology in this introductory treatment, we  only provide  rough formulations of the setting and the main results here.
  
  \subsection{A local mean-field interacting spin system}
    \label{S1.1.1} 
      Let $ T \in (0, \infty) $ and $ N\in\N  $. 
      We denote by $ \Tone  $  the one-dimensional unit torus. 
      Let $ \Psi :\R\ra\R $ and 
      $ J:\Tone \ra \R $ be two functions that satisfy Assumption \ref{Fass}  below.
      Moreover, let 
      $ B^N = (B^{i,N})_{i=0,\dots,N-1}  $ be an $ N $-dimensional Brownian motion and $ \mu_0^N \in \cM_1(\R^N) $, i.e. $  \mu_0^N $ is a probability measure on $ \R^N $.
      In this paper we consider a system of $ N  $  
  coupled stochastic differential equations given by
   
  \begin{align}\label{SpinSystem}
  	&d\t_t^{i,N} = - \Psi' \left(\t_t^{i,N} \right) \, dt + \frac{1}{N} \sum_{j=0	}^{N-1} J\left(\frac{i-j}{N}\right) \t_t^{j,N}\, dt  + \sqrt{2}\,  dB_t^{i,N},  \ \ \ t \in (0,T],\ 0\leq i \leq N-1, \nonumber\\
  	&(\t_0^{0,N}, \dots \t_0^{N-1,N}) \sim 	\mu_0^N.
  \end{align} 
  For each $ i = 0,\dots,N-1 $ and $ t \in [0,T] $,  we call $ \t_t^{i,N} $ the  \emph{spin value} at time $ t $ of a particle, which is \emph{located} at $ i/N \in \To $.
  For a detailed historical review on such models we refer to \cite[Section 1.1]{PatLDP}.

Define a \emph{microscopic Hamiltonian} $ H^N:\R^N \ra \R $ by
\begin{align}\label{Hamiltonian}
	H^N(\varTheta) &= 
	\sum_{i=0}^{N-1} \left( \Psi \left(\t^i \right)  - \frac{1}{2N} \sum_{j=0	}^{N-1} J\left(\frac{i-j}{N}\right) \t^{i} \t^{j} \right).
\end{align} 
Let  $ \Theta^N_t := (\t_t^{i,N})_{i=0,\dots,N-1} $  denote the  vector of all  $ N $ spins. Then we observe that 
\begin{align}\label{SpinSystemN}
	\begin{split}	
		d\Theta_t^{N} &= - \nabla H^N (\Theta_t^N) \,  dt  + \sqrt{2} \, dB_t^{N},\\
		\Theta^N_0 &\sim \mu_0^N.
	\end{split}
\end{align} 
Let  $ \mu_t^N $ denote the law  of $ \Theta_t^N  $ for each $ t\in[0,T] $. 
It is well-known that $ (\mu_t^N)_{t\in[0,T]} $ can be represented as a \emph{Wasserstein gradient flow} (see, e.g.\ \cite{jko} or \cite{ambgigsav}).
Moreover,  for each $ t $, $ \mu_t^N $ has a density $ \rho_t^N $ with respect to the Lebesgue measure on $ \R^N $ and   $ (\rho_t^N)_{t\in[0,T]} $ is a weak solution to the \emph{Fokker-Planck equation }
\begin{align}\label{FPE}
	\partial_t \rho^N_t = \Delta \rho^N_t + \mathrm{div} \left(	\nabla H^N \rho^N_t 		\right).
\end{align}  
In this paper we  focus on  curves of laws rather than on the path-wise solutions of  systems of stochastic differential equations.
Hence, instead of the systems \eqref{SpinSystem} and \eqref{SpinSystemN} we study $ (\mu_t^N)_{t\in[0,T]} $ and $ (\rho_t^N)_{t\in[0,T]} $. 
However, it is also possible to specify  the roles of \eqref{SpinSystem} and \eqref{SpinSystemN}  in the results of this paper; see \cite{HDLAntonDima}. 

In order  to analyse the curves $ (\mu_t^N)_{t\in[0,T]} $ as $ N \ra \infty $, we push all measures into the same space via the map $ K^N $ that sends a vector to the corresponding \emph{empirical pair measure}, i.e.\
\begin{align}
	\begin{split}
		K^N : \R^N &\ra \cM_1(\Tone \times \R)	\\
		\varTheta = (\t^k)_{k=0}^{N-1} &\to \frac{1}{N} \sum_{k=0}^{N-1} \delta_{\left(\frac{k}{N}, \theta^k\right)}.
	\end{split}
\end{align}
The  goal is to state a hydrodynamic limit result and a large deviation principle for the sequence $ \{((K^N)_\# \mu^N_t)_{t\in[0,T]} \}_N $, where 
$ (K^N)_\# \mu^N_t  $ denotes the  image measure of $ \mu_t^N $ under $ K^N $.

  \subsection{The limiting object}
  \label{S1.1.2}

We  want to  explain intuitively, what the limiting system should be.
Note that \eqref{SpinSystem} is of the form
\begin{align} 
d\t_t^{i,N} =
b\left(\frac{i}{N} \, ,\,  \t_t^{i,N} \, ; \,  K^N(\varTheta^N_t)\right) dt  + \sqrt{2}\, dB_t^{i,N},
\end{align}
where $ b : \Tone \times \R \times  \cM_1(\Tone \times \R) \ra \R $ is given by 
\begin{align}
b(x,\t;\nu ) = -\Psi'(\t ) + \int_{\Tone \times \R} J(x-x') \t' d\nu(x',\t').
\end{align}
This suggests that the limiting system should be 
\begin{align} 
d \hat{\t}_t^{x}  =
b\left( x \, , \,  \hat{\t}_t^{x} \,  ; \,  \mu_t \right) dt  + \sqrt{2}\, dB_t^{x}, \ \ \ \ \ x \in \To,
\end{align}
where  $\mu_t \in \cM_1(\Tone \times \R)$ is of the form  $ \mu_t = \mu_t^x \, dx $  and such that $ \mu_t^x $ is the law of $ \hat{\t}_t^{x} $ for all $ t $ and  $ x $. 
However, this in turn suggests that $ \mu_t  $ should have a density $ \rho_t $ with respect to the Lebesgue measure on $ \ToR $ for all $ t \in (0,T]$ and $ (\rho_t)_{t\in[0,T]} $ should be a weak solution of a partial differential equation of the form
\be{pdeEq}  
 \partial_t \rho_t(x,\t) = \partial_{\t\t}^2 \rho_t(x,\t) + \partial_\t \left(	\rho_t(x,\t) \left( \Psi'(\t) - 	\int J(x-\bar{x})\bar{\t}   \rho_t(\bar{x},\bar{\t})   \, d\bar{\t} d	\bar{x}	\right)					\right).
 \ee 
It is not possible to find a representation of this partial differential equation in the usual  Wasserstein setting, since there are no partial derivatives with respect to $ x $.
Hence, we have to modify the Wasserstein distance in such a way that the new metric takes into account that there is no evolution in this parameter.
It turns out that the correct distance is given by
\be{MetricDefIntro}
\WL(\mu,\nu)^2 :=  \int_{\To} W_2(\mu^x,\nu^x)^2 \, dx,
\ee 
where $ \mu= \mu^x \, dx \in \cM_1(\ToR) $ and $\nu = \nu^x \, dx \in \cM_1(\ToR) $ are suitable and	$ W_2 $ denotes the Wasserstein distance on the space of square-integrable probability measures on $ \R $; see Chapter \ref{S2} for the details.  
Now we have to rebuild the whole gradient flow theory as in the Wasserstein space in order to show that we can represent \eqref{pdeEq} in this new  framework. This is the content of Chapter \ref{S2}.

  \subsection{Results}
\label{S1.1.3}

  In this section, we  state our main results and sketch the ideas of the corresponding proofs.  
  The first  result is the gradient flow formulation of \eqref{pdeEq}.
  
  \begin{theoremm}[Gradient flow formulation, cf.\ Theorems \ref{FProp}, \ref{Lyapunov} and \ref{pde}]\mbox{}	\\  
  		Define  $ \cF: \cM_1(\Tone \times \R) \ra (-\infty, \infty] $ by
  	\begin{align}\label{EqThmGFWL}
  	\cF(\mu) :=
  	\cH( \mu | e^{-\Psi(\theta)} dx d\theta ) - \frac{1}{2} \int_{(\Tone \times \R)^2} J(x-x') \theta \theta' d\mu(x,\theta) d \mu(x',\theta'),
  	\end{align}
  	where $ \cH $ is the relative entropy functional (see \eqref{RelEntropyFunctional} below).
  	Let $ \mu_0 \in D(\cF) $, i.e.\ $ \cF(\mu_0) < \infty $. Then
 there exists a unique \emph{$ \WL$-gradient flow} $ (\mu_t)_{t\in[0,T]} $ for $ \cF $ with initial value $ \mu_0 $.
  		Moreover, for all $ t \in [0,T]$, $ \mu_t  $ has a density $ \rho_t $ with respect to the Lebesgue measure on $ \ToR $ and 
  			 $ (\rho_t)_{t \in [0,T]} $  is a weak solution to \eqref{pdeEq}.
Finally,
  	  $ (\mu_t)_{t \in [0,T]} $ is the unique $ \WL $-continuous curve  such that $ \lim_{t \downarrow 0} \WL(\mu_t,\mu_0)=0 $ and $ \cJ[(\mu_t)_{t \in [0,T]}]=0 $, where $  \cJ: C([0,T] \, ; \, \cM_1(\ToR) ) \ra [0,\infty]  $ and for smooth curves $ (\nu_t)_{t \in [0,T]} $, $ \cJ[(\nu_t)_{t \in [0,T]}] $ is given by 
\begin{align}
\cJ[(\nu_t)_{t \in [0,T]}]:=\cF(\nu_T) - 		\cF(\nu_0) + \frac 12 \int_0^T\big(|\partial \cF|^2(\nu_t)  + |\nu'|^2(t)\big)\,dt,
\end{align}  
where the objects $ |\partial \cF| $ and $ |\nu'| $ will be introduced in \eqref{MetricSlopeDefEq} and \eqref{MetricDerEq}, respectively. 
  	\end{theoremm}
  
  To prove this result, we have to develop the same theory for $ \WL $ as in the  book \cite{ambgigsav} for the Wasserstein space. 
  To this end, we first show that $ \WL $ is lower semi-continuous  with respect to weak convergence (Lemma \ref{LSCLem}) and that 
   $ (\PLo , \WL)$ is a Polish space (Paragraph \ref{Polish}), where $ \PLo $ is defined in \eqref{SpaceDef} below.
   Then we analyse curves in $ (\PLo , \WL)$ and characterize $ \WL $-absolutely continuous curve via distributional solutions of certain partial differential equations (Proposition \ref{AbsContProp}). This characterisation will later be the key fact to build the bridge to \eqref{pdeEq}.
In Section \ref{GradFlowSec}, we introduce a subdifferential calculus in $ (\PLo , \WL)$ and define the notion of gradient flows with it.
This allows us to apply the abstract theory of Part I of the book \cite{ambgigsav} to show existence, uniqueness and further  properties of $ \WL $-gradient flows in Theorem \ref{ExistenceThm}.
In Section \ref{McKean}, we finally consider the special case of the functional $ \cF $ and apply the previous results for this case and arrive at Result I.

  \begin{theoremm}[Large deviation principle, cf.\ Theorem \ref{LDP}] \mbox{}	\\  	
For all $ N \in \N $, let $ (\mu^N_t)_{t\in[0,T]} $ be 
 defined as in Section \ref{S1.1.1}.
 Let $ (\mu^N_0)_N $ satisfy Assumption \ref{AssI}.
 Then $ (\{(K^N)_\# \mu^N_t\}_{t\in[0,T]} )_N $ satisfies a large deviation principle  in $ C([0,T] ; \cM_1(\Tone \times \R) ) $ with rate function 
\begin{align}
	  (\nu_t)_t \to I[(\nu_t)_t] := \frac{1}{2} \cJ[(\nu_t)_t] + \cH (\nu_0 | \mu_0) 
\end{align}
for some $ \mu_0 \in D(\cF) $ (see Theorem \ref{LDP} for details).
  \end{theoremm}

The proof is based on the paper \cite{FathiLDP}  in the following way.
For each $ N $, an analogous statement as Result I holds true for $ (\mu^N_t)_{t\in[0,T]} $ with respect to some functional $ \cJ^N $; see e.g.\ \cite[11.2.1]{ambgigsav}.
Then, the results in \cite{FathiLDP} (combined with some additional arguments that we will provide in the proof of Theorem  \ref{LDP} on page \pageref{LDP}) show   that in order to prove the large deviation principle for $ (\{(K^N)_\# \mu^N_t\}_{t\in[0,T]} )_N $ it 
is equivalent to show that the following two claims hold:
\begin{itemize}
	\item If $ (\nu_t)_{t\in[0,T]}  \in C([0,T] ; \cM_1(\ToR)) $ and   $  (\nu^N_t)_{t\in[0,T]} \in C([0,T]; \cM_1(\R^N)) $ for all $ N \in \N $ are such that  $ (K^N)_\# \nu_t^N  \hra  \delta_{\nu_t} $ for all $ t \in[0,T] $, then 
	\be{LBEqIntro}
	\liminf_{N\ra \infty }  \frac{1}{N} \left(  \frac{1}{2}\cJ^N[ (\nu_t^N)_{t\in[0,T]} ]  + \cH(\nu_0^N	\, | \, \mu_0^N 	)	\right)	\geq  I[ (\nu_t)_{t\in[0,T]} ].
	\ee
	\item For all $ (\nu_t)_{t\in[0,T]}  \in C([0,T] ; \cM_1(\ToR)) $ there exists   $  (\nu^N_t)_{t\in[0,T]} \in C([0,T]; \cM_1(\R^N)) $ for all $ N \in \N $  such that  $ (K^N)_\# \nu_t^N  \hra  \delta_{\nu_t} $ for all $ t \in[0,T] $, and 
	\be{UBEqIntro}
	\limsup_{N\ra \infty }  \frac{1}{N} \left(  \frac{1}{2}\cJ^N[ (\nu_t^N)_{t\in[0,T]} ]  + \cH(\nu_0^N	\, | \, \mu_0^N 	)	\right)	\leq  I[ (\nu_t)_{t\in[0,T]} ].
	\ee
\end{itemize}
These two claims will be shown in Section \ref{LowerBoundSec} and \ref{UpperBoundSec}, respectively.
Therefore, the large deviation principle is related to a (variant of) gamma-convergence result of the functionals $ (\nu_t^N)_{t\in[0,T]} \to \frac{1}{2}\cJ^N[ (\nu_t^N)_{t} ]  + \cH(\nu_0^N	\, | \, \mu_0^N 	) $. 
We  explain this in more detail in Chapter \ref{LDPChap}.

  	\begin{theoremm}[Hydrodynamic limit; cf.\ Theorem \ref{HDL}]\mbox{}	\\  
  Let  $ (\mu_t)_{t\in[0,T]} $ be the $ \WL$-gradient flow for $ \cF $ with initial value $ \mu_0 \in D(\cF) $.
 For all $ N \in \N $, let $ (\mu^N_t)_{t\in[0,T]} $ be 
 defined as in Section \ref{S1.1.1}.
 Suppose that the sequence of initial conditions $ (\mu^N_0)_N $ is  such that
 $ ((K^N)_\# \mu^N_0 )_N $ converges to $ \delta_{\mu_0} $ weakly in $\cM_{1} ( \cM_{1}(\Tone \times \R))$ and 
 \begin{align}
 \lim_{N\ra \infty } \frac{1}{N} \cH(\mu_0^N\, |\, e^{-H^N} \Leb{\R^N} ) = \cF(\mu_0) .
 \end{align} 
 Then
 	$ ((K^N)_\# \mu^N_t )_N $ converges to $ \delta_{\mu_t} $ weakly in $\cM_{1} ( \cM_{1}(\Tone \times \R))$ for all $  t \in [0,T] $ and 
 	 \begin{align}
\lim_{N\ra \infty } \frac{1}{N} \cH(\mu_t^N\, |\, e^{-H^N} \Leb{\R^N} ) = \cF(\mu_t) \ \ \  \text{ for all $  t \in [0,T] $}.
 	\end{align} 	
 	Moreover, under some additional assumption on  $ \Psi $, the convergence holds even in a stronger topology, which is induced by the Wasserstein topology on $ \cM_{1}(\Tone \times \R) $.
  \end{theoremm}

The assumption on the initial configurations  here is weaker than in Assumption \ref{AssI}.
The proof uses the same strategy as in \cite{SandierSer} and \cite{Serfaty}.
Again,  the characterisation of $ (\mu^N_t)_t $ and $ (\mu_t)_{t \in [0,T]} $ as the unique minimizers of $ \cJ^N $ and $ \cJ $, respectively, plays an important role. There are three main ingredients: \emph{compactness},  \emph{superposition} and \emph{lower semi-continuity}. 
The compactness of $ (\{(K^N)_\# \mu^N_t\}_{t\in[0,T]} )_N $ follows from the Arzel\`a-Ascoli theorem.
For the superposition principle, which states that all limit points of $ (\{(K^N)_\# \mu^N_t\}_{t\in[0,T]} )_N $ can be represented via a probability measure $ \varUpsilon $ on $ C([0,T] ; \cM_1(\Tone \times \R) ) $, we use  \cite[Theorem 5]{Lisini}. Finally, the lower semi-continuity states that
\be{LSCHDLIntro}
\int \,  \Ind_{\mu_0}(\eta_0) \cdot  \cJ[\, (\eta_t)_t \,]\ d\varUpsilon((\eta_t)_t  ) \leq \liminf_{n \ra \infty  } \frac{1}{n} \cJ^n [\, (\mu_t^n)_t \, ] = 0 . 
\ee
This fact is an extension  of \eqref{LBEqIntro}. Since $ \cJ[\cdot] $ is a non-negative functional with unique minimizer $ (\mu_t)_{t \in [0,T]} $, \eqref{LSCHDLIntro} yields the claim. For more details, see Chapter \ref{HDLChap}.

  \begin{remark}
Most of the statements that we  prove in this paper can be extended easily. For instance, it 
	is possible to add a random environment, which is drawn according to $ 
	\varsigma \in \cM_1(\R) $ or to replace $ \To $ by a compact Riemannian manifold $ M $  or to allow the spins to  take values in $ \R^d $, for some $ d >1 $.
	The corresponding metric should then be of the form
	\be{IntroMetricDefExt}
	\mathrm{W}^{M,\varsigma}(\mu,\nu)^2 :=  \int_{M}\int_{\R} W_2(\mu^{m,\omega},\nu^{m,\omega})^2 \, d\varsigma(\omega) \, d\mathrm{vol}(m),
	\ee 
	Moreover, 
	it is possible to generalize  \eqref{pdeEq} in various ways without much additional work. 
	For instance, we could add a term of the form $  \partial_{\t\t}^2L_F(\rho_t(x,\t)) $ for some function $ L_F:[0,\infty) \ra [0,\infty)  $ as in  \cite[Example 9.3.6 and Subsection 10.4.3]{ambgigsav}, or we could include a diffusion coefficient as in \cite{FathiLDP}.
	It is also straightforward to see that the single-site potentials $ \Psi $ could also be dependent on the space parameter $ x $, and the quadratic interaction (given by the factor $ - \theta \theta' $ in \eqref{EqThmGFWL})
	 could be replaced by a more general class of interactions. 
	However, we try to keep the notation as simple as possible and did not try to optimize our results. 
  \end{remark}

\subsection{Notation and some definitions}  
\label{S1.2}

  In the following let $ n \in \N $ and $ (Y,\mathrm{d}), (\bar{Y},\mathrm{e}) , (Y_1,\mathrm{d}_1), \dots , (\bar{Y}_n,\mathrm{d}_n) $ be Polish spaces.

\vspace{-0.3cm}
\paragraph*{Measure theoretic notations.}
 \begin{itemize}
  	\setlength\itemsep{-0.3em} 
  	\item  $ \cM_1(Y)  $ denotes the space of Borel probability measures on $ Y $. 
  	We equip $ \cM_1(Y)  $ with the topology of weak convergence, where we say that 
  	$ (\mu_n)_{n\in \N} \subset \cM_1(Y)  $ \emph{converges weakly in $ \cM_1(Y)  $} to $ \mu \in \cM_1(Y)  $ (and write $ \mu_n \hra \mu $) if $ \int_Y f d\mu_n \ra \int_Y f d\mu $ for all $ f\in C_b(Y) $, i.e.\ for all continuous and bounded functions $ f : Y \ra \R $.
  	To emphasize the particular metric on $ Y $, we sometimes write that $ (\mu_n)_{n} $ \emph{converges weakly in $ \cM_1((Y,\mathrm{d}))  $} to $ \mu $.
  	\item For $ \mu \in \cM_1(Y)  $ and a Borel map $ f : Y\ra \bar{Y} $,  $ f_\#\mu  $ is the image measure of $ \mu $ by $ f $. 
  \item If $ Y \subset \R^d $ for some $ d \in \N $, we denote by $ \Leb{Y} $ the Lebesgue measure restricted to $ Y $.	
  \item We  denote elements in $ \R $ by $ \t $ or $ \bar{\t} $ and write  $ d\t $ instead of $ \Leb{\R} $. 
  In the same manner, for  $ N\in \N $, we  denote elements in $ \R^N $ by $ \varTheta = (\t^k)_{k=0}^{N-1} $ and write  $ d\varTheta $ instead of $ \Leb{\R^N} $. 
  \item Let $ \Td $ denote the $ d $-dimensional unit torus. 
  We usually denote elements in $ \Td $ by $ x $ or $ \bar{x} $
  and  write $ dx $ instead of $ \LebTd $.
  \item Define 
  \be{SpaceLeftMargDef}
  \cM_1^\mathrm{L}(\Td \times Y)  := \big\{	\mu \in \cM_1(\Td \times Y)	\, \big|\, \pb^1_\# \mu = \Leb{\Td}  \big\}.
  \ee
  By the disintegration theorem (see e.g.\ \cite[5.3.1]{ambgigsav}), for each $ \mu \in 	\cM_1^\mathrm{L}(\Td \times Y)   $, there exists a family $ (\mu^x)_{x\in \Td } $ of probability measures on $ Y $ such that $ x\to\mu^x $ is Borel-measurable and $ \mu = \mu^x \, dx  $, i.e.\ 
  \be{Disintegr}
  \int_{\Td \times Y} f(x,y) \, d \mu(x,y) =  \int_{\Td }\int_{Y} f(x,y) \,d \mu^x(y) dx   	
  \ee
  for all measurable and bounded $f:\Td \times Y \ra \R $. 
  \item Let $ \mu  $ and $ \nu $ be two measures on $ Y $. 
  Define the \emph{relative entropy} between $ \mu  $ and $ \nu $ by
  \begin{align}\label{RelEntropyFunctional}
  \cH(\mu \, | \, \nu ) := 
  \begin{cases}
  \int_{\TdR}  \log\left( \frac{d\mu}{d\nu}  \right) d\mu &: \mu \ll \nu, \\
  \infty   &: \text{else}.
  \end{cases} 
  \end{align}
  By  abuse of notation we use the same letter $ \cH $ for  all Polish spaces. 
\end{itemize}
\vspace{-0.6cm}
\paragraph*{Wasserstein spaces.}
  \begin{itemize}   	\setlength\itemsep{-0.3em}
  	\item 
  	By abuse of notation, for all Polish spaces $ (Y,\mathrm{d}) $, $ W_2 $ denotes the \emph{$ \L^2$-Wasserstein distance induced by $ \mathrm{d} $} on $ \cM_1(Y)  $, i.e.
  	\be{Wass}
  	W_2(\mu,\nu)^2 := \inf_{\gamma \in \CouplW{\mu,\nu} }  \int_{Y^2} \mathrm{d}(y, y')^2 \,  d\gamma(y,y') ,
  	\ee
  	where $ \mu, \nu \in \cM_1(Y) $ and $ \CouplW{\mu,\nu} $ denotes the space of all probability measures on $ Y^2 $ that have $ \mu  $ and $ \nu $ as marginals. We denote by $ \OptW{\mu,\nu} \subset \CouplW{\mu,\nu}  $ the set of all measures that realize the infimum  in \eqref{Wass} (cf.\ \cite[4.1]{vil}).
  	\item Set $ \cP_2(Y) := \{\mu \in \cM_1(Y) \, | \, \exists y_0 \in Y \, : \,  \int_Y d(y, y_0)^2 d\mu(y) < \infty 		\} $. Then $ (\cP_2(Y), W_2)$ is a Polish space (cf.\ \cite[6.18]{vil}).
  	If $ Y \subset \R^d $, then we denote by $ \cP_2^a(Y)$ the subset of $ \cP_2(Y)$ that consists of those measures that are absolutely continuous with respect to  $ \Leb{Y} $.
  	\item $ \widetilde{W} $ denotes the $ \L^2$-Wasserstein distance on $ \cM_1(Y)  $ induced by the distance $ \tilde{d}= d/(d+1) $. 
  	Then $ \widetilde{W} $ metrizes the weak topology on $ \cM_1(Y)  $ (cf.\ \cite[6.13]{vil}).
\end{itemize}
\vspace{-0.6cm}
\paragraph*{Some maps.}
\begin{itemize}   	\setlength\itemsep{-0.3em}
  	\item  For $ i \leq n $, let $ \pb^i : Y_1 \times \dots \times Y_n  \ra Y_i $ denote the projection on the $ i $-th component, i.e.\ $ \pb^i( y_1,\dots,y_n) = y_i$. 
  	Whenever it is necessary, we write $ \pb^i_{Y_1 \times \dots \times Y_n} $ instead of $ \pb^i $ to be able to distinguish different projection maps.
  	\item For $ t >0$, we denote by $ e_t $ the evaluation map at $ t $, i.e.\ $ e_t(f)= f(t) $ for all  $ f: (0,\infty) \ra Y $.
  	\item $ \Id_Y : Y\ra Y $ denotes the identity map on Y.
  \end{itemize}
\vspace{-0.6cm}
\paragraph*{Abbreviations.}
\begin{itemize}   	\setlength\itemsep{-0.3em}
	\item A function is $ \mathrm{d} $-l.s.c.\ if it is lower semi-continuous with respect to $ \mathrm{d} $.
\item For $ \varphi\in C^{1,0,1}((0,T) \times \TdR )   $ we often write $ \partial_t $ and $ \partial_\t $ to denote the partial derivative with respect to the parameter in $ (0,T) $ and $ \R $, respectively.
\item For $ a\in[-\infty,\infty] $, let $ a^+:=\max\{0,a\} $ and  $ a^-:=\max\{0,-a\} $.
\item We sometimes write $ (y_t)_t:= (y_t)_{t\in [0,T]}$ for curves $ (y_t)_{t\in [0,T]} \subset Y $.
  \end{itemize}


  
  \section{Gradient flow representation}
  \label{S2}

\subsection{Preliminaries}
  \label{S2.1}

In this section we will introduce a modification of the Wasserstein space and list some of its metric properties.
This space will provide  the framework to derive a gradient flow representation for the system in Section \ref{S1.1.2}.

The underlying space for this  representation is given by 
\be{SpaceDef}
\PL  := \left\{	\mu \in \ML	\, \Bigg|\,  \int_{\TdR} |\t|^2 d\mu(x,\t) < \infty			\right\}.
\ee
We equip $ \PL $ with the  distance
\be{MetricDef}
 \WL(\mu,\nu)^2 :=  \int_{\Td} W_2(\mu^x,\nu^x)^2 \, dx, \ \text{  } \mu,\nu \in \PL.
\ee
Here we have used that the map $ x\to W_2(\mu^x,\nu^x) $ is measurable.
This is true, since, by the measurable selection lemma (\cite[5.22]{vil}), for all $ \mu,\nu \in \PL $ there exists a family $ (\pi^x)_{x\in \Td } $ of probability measures on $ \R^2 $ such that $ x\to\pi^x $ is Borel-measurable and $ \pi^x \in \OptW{\mu^x,\nu^x}  	\text{  for almost every\ } x\in \Td  $. 
Defining $ \pi \in \cM_1^{\L}(\TdR\times\R) $ by 
$ \pi = \pi^x \, dx $, we observe that the set 
%
\begin{align}\label{OptDef}
\nonumber
\Opt{\mu,\nu} :=  \Big\{\pi \in \cM_1^{\L}(\TdR\times\R) 	\, \Big|\, 
& \pi = \pi^x \, dx, \text{ where } x\to\pi^x  \text{ is Borel-measurable and  } \\  &\pi^x \in \OptW{\mu^x,\nu^x}  	\text{  for almost every } x\in \Td 	\Big\}
\end{align}
is non-empty.
Note that 
\be{MetricOptRel}
 \WL(\mu,\nu)^2 =  \int_{\TdR\times \R} |\t - \t'|^2 d\pi(x,\t,\t') \ \text{ for all  } \pi \in \Opt{\mu,\nu}.
\ee

Moreover, $ \WL $ can be connected more directly to an optimal transportation problem, since
\cite[12.4.6]{ambgigsav} shows that for all $ \mu,\nu \in \PL $
\be{MetricCouplRel}
 \WL(\mu,\nu)^2 = \inf_{\gamma\in \Coupl{\mu,\nu}} \int_{\TdR\times \R} |\t - \t'|^2 d\gamma(x,\t,\t') ,
\ee
where 
\be{CouplDef}
\Coupl{\mu,\nu} := \Big\{\gamma \in \cM_1^{\L}(\TdR\times\R) 	\, \Big|\, 
\pb^{1,2}_\#\gamma=  \mu, \, \pb^{1,3}_\#\gamma=  \nu	\Big\}.
\ee
Using \eqref{MetricCouplRel}, it is easy to extend the definition of $ \WL $ to the whole space $ \ML $. 
Further,  \cite[5.3.2]{ambgigsav} yields that 
\begin{align}\label{CouplReprEq}
\nonumber \Coupl{\mu,\nu} = \Big\{\gamma \in \cM_1^{\L}(\TdR\times\R) 	\, \Big|\, 
&\gamma = \gamma^x \, dx, \text{ where } x\to\gamma^x  \text{ is Borel-measurable and  } \\  &\gamma^x \in \CouplW{\mu^x,\nu^x}  	\text{  for almost every\ } x\in \Td 	\Big\}. 
\end{align}
This implies that $ \Opt{\mu,\nu} \subset \Coupl{\mu,\nu} $. 
Therefore, it is easy to see that $ \Opt{\mu,\nu} $ is the set of minimizers in \eqref{MetricCouplRel}.
From now on, we call the elements of $ \Opt{\mu,\nu}  $ \emph{L-optimal plans  between $ \mu $ and $ \nu $}, and the elements of $ \Coupl{\mu,\nu} $ \emph{L-couplings  of $ \mu $ and $ \nu $}.

\paragraph{Comparison between $ \WL $ and $ W_2 $.}
Let   $ W_2 $ denote the  Wasserstein distance  on $ \cP_2(\TdR) $.
Then we have 
\be{CompWassEq}
 \WL(\mu,\nu) \geq W_2(\mu,\nu) \ \text{ for all }  \mu,\nu \in \PL .
\ee
Indeed, this can be shown by estimating the Wasserstein distance by the L$ ^2 $-norm with respect to  $ (\pb^1,\pb^2,\pb^1,\pb^3)_\#\pi \in \CouplW{\mu,\nu}$, where $ \pi \in \Opt{\mu,\nu} $.
However, there is no equality in general as it can be seen from the following example. 
Let $ A:= \{	x\in \Td \, | \, x_1 \leq \frac{1}{2}		\} $ and define $  \mu,\nu \in \PL  $ by 
\begin{align}\label{CounterEx}
\begin{split}
&\mu(dx,d\t) :=\Ind_A(x) \delta_0(d\t) dx + \Ind_{A^c}(x) \delta_1(d\t) dx, 
\\ 
&\nu(dx,d\t) :=\Ind_A(x) \delta_1(d\t) dx + \Ind_{A^c}(x) \delta_0(d\t) dx.
\end{split}
\end{align}
Then it is easy to see that $  \WL(\mu,\nu)=1 $ and $ W_2(\mu,\nu)\leq \frac{1}{4} $.

\paragraph{The absolutely continuous case.}
Let us consider the special case, when the measures are absolutely continuous with respect to $ \Leb{\TdR} $. 
Set 
\be{SpaceAbsContDef}
\PLa  = \big\{	\mu \in \PL	\, \big|\,  \mu  \ll \Leb{\TdR}  \big\}.
\ee
It is clear that, if $ \mu \in \PLa $, then $ \mu^x \in  \Pwa $ for almost every $ x \in \Td  $.
Consequently, if $ \nu \in \PL $, then $ \OptW{\mu^x, \nu^x}  = \{  (\mathrm{Id_\R},\T_{\mu^x}^{\nu^x})_\# \mu^x	\}$ for some $ \T_{\mu^x}^{\nu^x} \in \mathrm{L}^2(\mu^x) $ for almost every $ x \in \Td  $ (cf.\ \cite[10.42]{vil}). Hence, $ \Opt{\mu, \nu}  = \{  (\mathrm{Id_\R},\T_{\mu^x}^{\nu^x})_\# \mu^x	\,  dx\}$.
\bl{OptimalMapLem}
Let $ \mu \in \PLa $ and $ \nu \in \PL $.
Then there exists a unique map $ \T_{\mu}^{\nu} \in \mathrm{L}^2(\mu) $ such that 
\bi
\item $ \T_{\mu}^{\nu}(x,\t) = \T_{\mu^x}^{\nu^x}(\t)  $ for almost every $ x \in \Td  $,
\item $  \WL(\mu,\nu) = \| \pb^2 - \T_{\mu}^{\nu} 	\|_{\mathrm{L}^2(\mu)} $.
\ei
In the following we call  $ \T_{\mu}^{\nu} $ the \emph{L-optimal map} between $ \mu  $ and $ \nu  $.
\el
\bpr
Let $ \pi \in \Opt{\mu,\nu} $. 
Define a linear map $ L : \mathrm{L}^2(\mu) \ra \R $ by
\be{OptimalMapLemEq}
L(g) :=  \int_{\TdR\times \R} g(x, \t) (\t - \t') d\pi(x,\t,\t').
\ee
Due to the monotone-class theorem and the fact that $ x\to\pi^x $ is Borel-measurable, the integrand is measurable. 
Next we apply the Cauchy-Schwartz inequality to obtain
\be{OptimalMapLemEq2}
|L(g)| \leq   \| g \|_{\mathrm{L}^2(\pi)}  \WL(\mu,\nu) = \| g \|_{\mathrm{L}^2(\mu)}  \WL(\mu,\nu).
\ee
Hence, the Riesz representation theorem yields the existence of a unique element  $ f \in \mathrm{L}^2(\mu) $ such that $ L(g) = \int fg \, d\mu$ for all $ g \in \mathrm{L}^2(\mu) $.
Thus
\begin{align}\label{OptimalMapLemEq3}
\begin{split}
\int_{\TdR} fg \, d\mu = L(g) &= \int_{\TdR\times \R} g(x, \t) (\t - \t') d\pi^xdx\\
&=  \int_{\TdR\times \R} g(x, \t) (\t - \T_{\mu^x}^{\nu^x}(\t) ) d\mu.
\end{split}
\end{align}
Hence, $ f(x,\t) = \t - \T_{\mu^x}^{\nu^x}(\t)$ $ \mu $-a.e. 
Defining $ \T_{\mu}^{\nu} := \pb^2 - f $ yields the desired results.
\epr

\paragraph{Stability of L-couplings and L-optimal plans.} 
First we want to show that a sequence of L-couplings converges weakly if the corresponding sequences of marginals converge. 
For $ \cK, \cL \subset \cM_1^\mathrm{L}(\TdR)  $, define
\begin{align}\label{CouplSetReprEq}
\Coupl{\cK,\cL} := \big\{\gamma \in \cM_1(\TdR\times\R) 	\, \big|\, 
\exists \mu\in\cK, \nu \in \cL \, : \, \gamma \in \Coupl{\mu,\nu} 	\big\}. 
\end{align}

\bl{StabilityCouplLem}
\begin{enumerate}[(i)]
	\item If $ \cK $ and $ \cL $ are both tight subsets  of  $ \ML $, then 
	$ \Coupl{\cK,\cL} $ is  a tight subset of $ \cM_1(\TdR\times\R) $.
	\item If $ \cK $ and $ \cL $ are both compact with respect to the weak topology   in  $ \ML $, then 
	$ \Coupl{\cK,\cL} $ is compact with respect to the weak topology   in $ \cM_1(\TdR\times\R) $. 
\end{enumerate}
\el
\bpr
We skip this proof as it is a straightforward modification of the analogous result in the setting of the Kantorovich problem; see e.g.\ \cite[4.4]{vil}. 
\epr

We prove the analogous result for L-optimal plans only in the following special case. 
\bl{StabilityOptLem}
Let $ (\mu_n)_{n\in \N} \subset \PL$ and $ \mu \in \PL$ be such that 
for all subsequences $ (\mu_k)_k$, there exists a subsequence $ (\mu_{k_l})_{l}$ and a $ \LebTd$-nullset $ \cN_k $ such that 
\be{StabilityOptLemEq1}
\mu_{k_l}^x \hra \mu^x \ \  \text{ for all } x \in \Td \setminus \cN_k.
\ee
Let  $ \pi_n \in \Opt{\mu_n,\mu} $ for all $ n $.
Then
\begin{align}\label{StabilityOptLemEq}
\begin{split}
\pi_n  \hra  (\Id_\R, \Id_\R)_\# \mu^x dx. 
\end{split}
\end{align}
\el
\bpr
Let $ (\mu_k)_k $ be a subsequence. 
From the assumptions and from the stability of optimal plans in $ (\Pw,W_2) $ (\cite[5.21]{vil}) and since $ \OptW{\mu^x,\mu^x} = \{ (\Id_\R, \Id_\R)_\# \mu^x\} $, we have that
\be{StabilityOptLemEq2}
\pi_{k_l}^x \hra  (\Id_\R, \Id_\R)_\# \mu^x \ \  \text{ for all } x \in \Td \setminus \cN_k.
\ee
Let $ f \in C_b(\TdR\times \R) $. 
Then the dominated convergence theorem yields
\begin{align}\label{StabilityOptLemEq3}
\begin{split}
\lim_{l \ra \infty}  \int_{\TdR\times \R} f \, d\pi_{k_l} = 
 \int_{\Td}  \lim_{l \ra \infty} \int_{\R\times \R} f d\pi_{k_l}^x dx 
 =   \int_{\Td}  \int_{\R\times \R} f  \, d(\Id_\R, \Id_\R)_\# \mu^x dx 
\end{split}
\end{align}
Hence, $ \pi_{k_l}  \hra  (\Id_\R, \Id_\R)_\# \mu^x dx $.
And since the weak topology in $ \cM_1(\TdR\times\R) $ is metrizable, we infer the weak convergence of the whole sequence $ (\pi_n)_n $ towards $ (\Id_\R, \Id_\R)_\# \mu^x dx $.
\epr

\paragraph{Weak lower semi-continuity of $ \WL $.}  

\bl{LSCLem}
Let $ (\mu_n)_n,(\nu_n)_n \subset \cM_1^\mathrm{L}(\TdR)  $ and $ \mu,\nu \in \cM_1^\mathrm{L}(\TdR)  $ be such that 
$ \mu_n \hra \mu $ and $ \nu_n \hra \nu $. 
Then:
\be{LSCLemEq}
\liminf_{n \ra \infty}  \WL(\mu_n,\nu_n) \geq  \WL(\mu,\nu).
\ee
\el
\bpr
Consider a subsequence such that $ \lim_{k \ra \infty}  \WL(\mu_k,\nu_k)= \liminf_{n \ra \infty}  \WL(\mu_n,\nu_n) $.
Let $ \pi_k \in \Opt{\mu_k,\nu_k} $ for all  $ k $.
 Lemma \ref{StabilityCouplLem} yields the existence of a subsequence $ (\pi_{k_l})_l $ such that $ \pi_{k_l} \hra \pi $ for some $ \pi \in \Coupl{\mu,\nu} $.
Then
\begin{align}\label{LSCLemEq1}
	\begin{split}
 \liminf_{n \ra \infty}  \WL(\mu_n,\nu_n)^2 &=	\lim_{l \ra \infty}  \WL(\mu_{k_l},\nu_{k_l})^2
 = 	\lim_{l \ra \infty} \int_{\TdR\times\R} |\t-\t'|^2 \, d\pi_{k_l} \\
&\geq \int_{\TdR\times\R} |\t-\t'|^2 \, d\pi \, \geq  \,  \WL(\mu,\nu),
 	\end{split}
\end{align}
where the first inequality is due to a standard lower semi-continuity result for integrals (see e.g.\  \cite[5.1.7]{ambgigsav}) and the second inequality is due to \eqref{MetricCouplRel}. 
\epr

\paragraph{Characterization of convergence in $ \boldsymbol{(\PL,\WL)} $.}
Convergence with respect to the Wasserstein distance can be characterized by weak convergence plus convergence of the moments. 
A similar fact is true for convergence in $ (\PL,\WL) $. 

\bp{ConvCharacLem}
Let $ (\mu_n)_n\subset \PL  $ and $ \mu \in \PL  $.
Then $ \lim_{n \ra \infty}  \WL(\mu_n,\mu) = 0 $ if and only if 
\begin{enumerate}[(i)]
	\item $ \lim_{n \ra \infty} \int_{\TdR} |\t|^2 d\mu_n = \int_{\TdR} |\t|^2 d\mu $, and \label{ConvCharacLemPropMom}
	\item \label{ConvCharacLemPropWeak}	For all subsequences $ (\mu_k)_k$, there exists a subsequence $ (\mu_{k_l})_{l}$ and a $ \LebTd$-nullset $ \cN_k $ such that 
	\be{ConvCharacLemEq0}
\mu_{k_l}^x \hra \mu^x \ \  \text{ for all } x \in \Td \setminus \cN_k.
\ee
\end{enumerate}
\ep
\bpr
Assume that $ \lim_{n \ra \infty}  \WL(\mu_n,\mu) = 0 $.
(\ref{ConvCharacLemPropMom}) is a simple consequence of the triangle inequality for $ \WL $, which we will prove below in Lemma \ref{MetricSpaceLem}. Indeed,
\begin{align}\label{ConvCharacLemEq}
 \left|\left( \int_{\TdR} |\t|^2 d\mu_n \right)^{\frac{1}{2}}- \left(\int_{\TdR} |\t|^2 d\mu\right)^{\frac{1}{2}} \right|  
&= 
 \left| \WL(\mu_n, \delta_0 \otimes \LebTd) -  \WL(\mu, \delta_0 \otimes \LebTd) \right| \nonumber \\
&\leq  \WL(\mu_n, \mu) \ \Ra 0. 
\end{align}
To show (\ref{ConvCharacLemPropWeak}), let $ (\mu_k)_k$ be a subsequence.
Note that the  function $ x \to W_2(\mu_k^x,\mu^x) $ converges to $ 0 $ in $ \L^2(\Td) $.
Hence, there exists a further subsequence $ (\mu_{k_l})_l $ and a $ \LebTd$-nullset $ \cN_k $ such that 
\be{ConvCharacLemEq7}
\lim_{l \ra \infty} W_2(\mu_{k_l}^x,\mu^x) = 0 \ \  \text{ for all } x \in \Td \setminus \cN_k.
\ee
This yields \eqref{ConvCharacLemEq0}, since Wasserstein convergence implies weak convergence.

Conversely, assume (\ref{ConvCharacLemPropMom}) and (\ref{ConvCharacLemPropWeak}). 
Let $ \pi_n \in \Opt{\mu_n,\mu} $ for all $ n $.
 Lemma  \ref{StabilityOptLem} shows  that
(\ref{ConvCharacLemPropWeak}) implies
\begin{align}\label{ConvCharacLemEq2}
\pi_n  \hra  (\Id_\R, \Id_\R)_\# \mu^x dx. 
\end{align}
It is a simple consequence of (\ref{ConvCharacLemPropWeak}), the dominated convergence theorem and the metrizability of weak convergence that     
\begin{align}\label{ConvCharacLemEq3}
\mu_n  \hra   \mu. 
\end{align}
Proceeding exactly as in the Wasserstein case (see e.g.\ the last part of the proof of \cite[6.9]{vil}), 
we can show that (\ref{ConvCharacLemPropMom}), \eqref{ConvCharacLemEq2} and \eqref{ConvCharacLemEq3}
imply  $ \lim_{n \ra \infty}  \WL(\mu_n,\mu) = 0 $. 
Again, we skip the details as there will be no new insights. 
\epr

\paragraph{$ \boldsymbol{(\PL,\WL)} $ is a Polish space.} \label{Polish}
\bl{MetricSpaceLem}
$ (\PL,\WL) $ is a metric space.
\el
\bpr
$ \WL $ is well-defined on $ \PL $, since for all $ \mu,\nu  \in \PL $
\begin{align}
 \WL(\mu,\nu)^2 &
\leq  \int_{\Td} (W_2(\mu^x,\delta_0) + W_2(\delta_0,\nu^x))^2 dx
\leq 4 \int_{\Td} (W_2(\mu^x,\delta_0)^2 + W_2(\delta_0,\nu^x)^2 ) dx  \nonumber \\
&= 
 4 \int_{\TdR} |\t|^2 d\mu + 4  \int_{\TdR} |\t|^2 d\nu < \infty.
\end{align}
$ \WL $ is symmetric, since the Wasserstein distance on $ \R $ is symmetric.  
Let $ \mu,\nu  \in \PL. $ If $ \mu = \nu $, then $ \mu^x=\nu^x $ for a.e.\ $ x\in \Td $ by the uniqueness claim in the disintegration theorem, and therefore $  \WL(\mu,\nu) = 0 $. 
And if $  \WL(\mu,\nu) = 0 $, then necessarily $ W_2(\mu^x,\nu^x)= 0  $ for a.e.\ $ x\in \Td $. This implies that $ \mu^x=\nu^x $ for a.e.\ $ x\in \Td $, and hence $ \mu = \nu $. It remains to show the triangle inequality. Let $ \mu,\nu,\sigma  \in \PL. $ Then
\begin{align}
 \WL(\mu,\nu) 
 &= \left(\int_{\Td} W_2(\mu^x,\nu^x)^2 dx \right)^{\frac{1}{2}}
\leq  \left(\int_{\Td} (W_2(\s^x,\mu^x) + W_2(\s^x,\nu^x))^2 dx		 \right)^{\frac{1}{2}}   \\\
&\leq  \left(\int_{\Td} W_2(\s^x,\mu^x) ^2 dx\right)^{\frac{1}{2}}+  \left(\int_{\Td} W_2(\s^x,\nu^x)^2 dx\right)^{\frac{1}{2}}
=  \WL(\s,\mu) +  \WL(\s,\nu),\nonumber
\end{align}
where we have used the triangle inequality for the Wasserstein distance and Minkowski's inequality.
\epr

\bl{CompleteLem}
$ (\PL,\WL) $ is complete.
\el
\bpr
Let $ (\mu_n)_n $ be a Cauchy sequence in $ (\PL,\WL) $.
Let $ \e >0 $.
There exists $ N_\e >0 $ such that $  \WL(\mu_n, \mu_m) < \e  $ for all $ n,m \geq N_\e $.
Then if $ n \geq  N_e $
\begin{align}\label{CompleteLemEq}
\begin{split}
\Bigg( \int_{\TdR} |\t|^2 d\mu_n \Bigg)^{\frac{1}{2}}
\leq  \WL(\mu_n,\mu_{N_\e}) +  \WL(\mu_{N_\e}, \delta_0 \otimes \LebTd) \leq \e + \max_{i\leq  N_\e} \Bigg( \int_{\TdR} |\t|^2 d\mu_i \Bigg)^{\frac{1}{2}}.
\end{split}
\end{align}
Therefore, 
\begin{align}\label{CompleteLemEq1}
\begin{split}
\sup_{n\in \N} \Bigg( \int_{\TdR} |\t|^2 d\mu_n \Bigg)^{\frac{1}{2}}
\leq \e + \max_{i\leq  N_\e} \Bigg( \int_{\TdR} |\t|^2 d\mu_i \Bigg)^{\frac{1}{2}} < \infty,
\end{split}
\end{align}
and we infer the existence of a weakly converging subsequence $ (\mu_k)_k $ with  limit point $ \hat{\mu} \in \ML $.
The weak lower semi-continuity of $ \nu \to \int |\t|^2 d\nu $ and \eqref{CompleteLemEq1} imply that even 
$ \hat{\mu} \in \PL $.\
Finally, the weak lower semi-continuity of $ \WL $ yields 
\be{CompleteLemEq2}
\lim_{n \ra \infty}  \WL(\mu_n, \hat{\mu}) \leq \lim_{n \ra \infty} \liminf_{k \ra \infty}  \WL(\mu_n, \mu_k) =0,
\ee
since $ (\mu_n)_n $ is Cauchy.
Thus $ (\mu_n)_n $ is a converging sequence in $ (\PL,\WL) $.
\epr

\bl{SeparableLem}
$ (\PL,\WL) $ is separable.
\el
\bpr
To simplify the notation, we only give the proof for the case $ d = 1 $.
Let $ D \subset \Pw$ be countable and dense with respect to $ W_2 $.
Let for all $ n\in \N $ and $ k \leq 2^n-1 $, $ A_{k,n} = [k 2^{-n},(k+1) 2^{-n}) $.
Define
\be{SeparableLemEq}
\cD := \bigcup_{n\in \N} \ \bigcup_{	\{\nu_{k}^{n}	\}_{k=0, \dots,2^{n}-1} \subset D	} \left\{ \sum_{k=0	}^{2^{n}-1}		\Ind_{A_{k,n}} (x) \, \nu_k^n\, dx	\right\}
\ee
Then $ \cD $ is countable and $ \cD \subset \PLone  $.
In the following we show that $ \cD $ is dense  in $ (\PLone,\WL) $.

Define for all $ n $, the operator $ \mathrm{S}_n : \PLone \ra \PLone $ by 
\bea{SeparableLemEq1}
 \mathrm{S}_n(\mu) := \sum_{k=0	}^{2^{n }-1}	\, 	\Ind_{A_{k,n}} (x) \, S_{k,n}(\mu)\,  dx, 	\ \ \mu \in \PLone,
\eea
where for all $ k \leq  2^{n }-1$,  $ S_{k,n} : \PLone \ra \Pw $ is the operator that sends $ 	\mu \in \PLone $ to the averaged measure $ S_{k,n}(\mu)  = 2^{n }	\int_{A_{k,n}}   d\mu^x dx $ defined by
\bea{SeparableLemEq2}
\int_{\R} f \, dS_{k,n}(\mu) = 2^{n }	\int_{A_{k,n}} \int_{\R} f \,  d\mu^x dx , \ \text{ for all measurable, bounded }f:\R \Ra \R. 	
\eea
Let $ \mu \in \PLone $.
The proof of this lemma consists of showing the following two facts.

\noindent
\begin{itemize}
\item[(i)] 
For all 
 $ \e >0 $ and $ n\in\N $ there exists $ \nu^n \in \cD $ such that  $  \WL( \mathrm{S}_n(\mu),\nu^n) < \e $.
 \item[(ii)]
$ \lim_{n\uparrow\infty}  \WL( \mathrm{S}_n(\mu),\mu) =0 $.
\end{itemize}
Indeed,  statements (i) and (ii) imply that  for any $\mu$, there exists a sequence $(\nu^n)_n\subset \cD$ such 
that $  \WL( \nu^n,\mu) \ra 0 $, that is, $ \cD $ is dense in $ \PLone $.
 
We now show statement (i). 
Since $ D $ is dense in $ \Pw $, there exists $ \nu_{k,n} \in D$ such that $ W_2(\nu_{k,n} , S_{k,n}(\mu)) < \e $ for all $ k \leq 2^n-1$.
Set $ \nu^n = \sum_{k=0	}^{2^{n } -1}		\Ind_{A_{k,n}} (x) \nu_{k,n} \, dx $. 
We immediately observe that $  \WL( \mathrm{S}_n(\mu),\nu^n) < \e $.

Next we prove (ii). 
In view of Proposition \ref{ConvCharacLem}, it will be enough to show that 
\begin{enumerate}[(A)]
	\item $ \int |\t|^2 d \mathrm{S}_n(\mu) =   \int |\t|^2 d\mu $ for all $ n $, and 
	\item $  \mathrm{S}_n(\mu)^x \hra \mu^x \  \text{ for almost every } x \in \Tone. $
\end{enumerate}
(A) is a simple consequence of   \eqref{SeparableLemEq2}. 
It remains to show (B), which will be done in six steps.
The main problem is to avoid the non-separability of the space $ C_b(\R) $. 
We do this in a standard way, which was done e.g.\ in the proof of \cite[11.4.1]{dudley}. 
%
%
This  means, we will push the measures down from $ \Tone \times \R $  to a bounded set. 
Consider $ h(\t) = \arctan(\t) $ and abbreviate $ O:= (\pi/2,\pi/2) $.
Set  $ \sigma = (\pb^1, h)_\# \mu$.
Consequently, $ \sigma $ is supported in $ \Tone \times O. $ 
Let $ \mathrm{BL}(O) $ be the set of real-valued bounded Lipschitz functions on $ O $.

\noindent
\textbf{Step 1.} \big[$ \forall f \in  \mathrm{BL}(O)  \, \exists \, \text{nullset }\cN^f : \int f \, d\mathrm{S}_n(\sigma)^x \ra \int f \, d\sigma^x  \ \ \ \forall x \in \Tone \setminus \cN^f $.\big]

\noindent
Let $  \Tone \setminus \cN^f $ be the set of Lebesgue-points of   $ x \mapsto  \int_\R f d\sigma^x \in \L^1(\Tone)$. 
For each  $ x \in \Tone $, let $ k_x(n) = \lfloor x 2^{n} \rfloor $.
Hence,  $ x \in A_{k_x(n),n} $ for each $ n $.
Denote by $ B(x,  2^{-n}) $  the ball of radius $  2^{-n} $ around $ x \in \Tone $.
Then we observe that for each $ x \in \Tone \setminus \cN^f $
\begin{align}
\begin{split}
\left|\int_O f d\mathrm{S}_n(\sigma)^x - \int_O f d\sigma^x\right| 
&=  \left|\int_O f dS_{k_x(n),n}(\sigma) - \int_O f d\sigma^x\right|\\
&\leq  2^{n} \int_{A_{k_x(n),n}} \left|\int_O f d\sigma^y - \int_O f d\sigma^x\right| dy \\
&\leq  \frac{ 2}{\Leb{\Tone}(	B(x,  2^{-n})	)} \int_{B(x,  2^{-n})} \left|\int_O f d\sigma^y - \int_O f d\sigma^x\right| dy \\
&\Ra 0 \quad \text{ as } \ n \ra \infty, 
\end{split}
\end{align}
since $ x $ is a Lebesgue point. 

\noindent
\textbf{Step 2.} 
\big[\emph{Let $ \iota : O \ra \bar{O} $ be the canonical inclusion, then }

\hspace{1cm} $  \forall \bar{f} \in  \mathrm{BL}(\bar{O})  \, \exists \, \text{nullset }\cN^{\bar{f}}: \int \bar{f} \, d\iota_\#\mathrm{S}_n(\sigma)^x \ra \int f \, d\iota_\#\sigma^x  \ \ \  \forall x \in \Tone \setminus \cN^{\bar{f}}. \big]$

\noindent
$ \bar{f} $ has the representation
\be{SeparableLemEq3}
\bar{f}(\t) = \inf_{\vartheta \in \bar{O}} \bar{f}(\vartheta ) + \mathrm{Lip} ( \bar{f}) \, | \t - \vartheta | = \inf_{\vartheta \in O} \bar{f}(\vartheta ) + \mathrm{Lip} ( \bar{f}) \, | \t - \vartheta | ,
\ee
where $ \mathrm{Lip} ( \bar{f})  $ is the Lipschitz-constant of $ \bar{f} $.
Define $ f \in  \mathrm{BL}(O)  $ by $ f(\t ):=  \inf_{\vartheta \in O} \bar{f}(\vartheta ) + \mathrm{Lip} ( \bar{f}) \, | \t - \vartheta |  $.
Then $ \bar{f} = f $ on $ O $.
Set $ \cN^{\bar{f}} := \cN^f $, where $ \cN^f $ is the nullset from Step 1. 
Then, since $ \iota_\#\mathrm{S}_n(\sigma) $ and $ \iota_\#\sigma $ are supported on $ O $, we obtain that for all $x \in \Tone \setminus \cN^{\bar{f}} $
\begin{align}
\begin{split}
\lim_{n\ra \infty} \int_{\bar{O}}  \bar{f} \, d\iota_\#\mathrm{S}_n(\sigma)^x &= \lim_{n\ra \infty} \int_{O}  f \, d\iota_\#\mathrm{S}_n(\sigma)^x
= \lim_{n\ra \infty} \int_{O}  f \, d\mathrm{S}_n(\sigma)^x\\
&=\int_{O}  f \, d\sigma^x = \int_{\bar{O}}  \bar{f} \, d\iota_\#\sigma^x.
\end{split}
\end{align}

\noindent
\textbf{Step 3.} 
\big[$  \exists \, \text{nullset }\cN: \int \bar{f} \, d\iota_\#\mathrm{S}_n(\sigma)^x \ra \int f \, d\iota_\#\sigma^x  \ \ \  \forall x \in \Tone \setminus \cN \ \ \  \forall \bar{f} \in  \mathrm{BL}(\bar{O})  $.\big]

\noindent
$ \mathrm{BL}(\bar{O}) $ is separable, i.e.\ there exists a countable set $ E \subset  \mathrm{BL}(\bar{O}) $, which is  dense with respect to $ \| \cdot \|_\infty $.
Set $ \cN := \cup_{\bar{f} \in E} \cN_k^{\bar{f}} $. 
Since $ E $ is dense in $ \mathrm{BL}(\bar{O})  $, this concludes the claim. 

\noindent
\textbf{Step 4.} 
\big[$  \exists \, \text{nullset }\cN: \int f \, d\mathrm{S}_n(\sigma)^x \ra \int f \, d\sigma^x  \ \ \  \forall x \in \Tone \setminus \cN \ \ \  \forall f \in  \mathrm{BL}(O)  $.\big]

\noindent
Using \cite[6.1.1]{dudley} we know that there exists 	$ \bar{f} \in  \mathrm{BL}(\bar{O})  $ such that $ \bar{f} = f $ on $ O $.
Now the claim follows immediately from Step 3. 

\noindent
\textbf{Step 5.} 
\big[$  \exists \, \text{nullset }\cN: \int f \, d\mathrm{S}_n(\sigma)^x \ra \int f \, d\sigma^x  \ \ \  \forall x \in \Tone \setminus \cN \ \ \  \forall f \in  C_b(O)  $.\big]

\noindent
The claim follows from Step 4 and  \cite[11.3.3]{dudley}.

\noindent
\textbf{Step 6.} 
\big[$  \exists \, \text{nullset }\cN: \int f \, d\mathrm{S}_n(\mu)^x \ra \int f \, d\mu^x  \ \ \  \forall x \in \Tone \setminus \cN \ \ \  \forall f \in  C_b(\R) $.\big]

\noindent
Note that $ \mathrm{S}_n(\sigma)^x = (h^{-1})_\#\mathrm{S}_n(\mu)^x  $ and $ \sigma^x = (h^{-1})_\#\mu^x  $ for all $ x \in \Tone \setminus \cN  $.
Hence, the claim follows from the continuous mapping theorem (see e.g.\ \cite[5.2.1]{ambgigsav}).
This concludes the proof. 
\epr

\subsection{\texorpdfstring{Curves in $ (\PL,\WL) $}{Curves in (P\_2\textasciicircum L ,W\textasciicircum L)}}  \label{S2.2}

In this section we analyse 
 \emph{geodesics} and \emph{absolutely continuous curves} in $ (\PL,\WL) $. 
 For the latter we will show that these curves are characterized by weak solutions of some type of \emph{continuity equation} and we introduce a notion of \emph{tangent velocity} at these curves. 
 This fact will be the main key later to represent weak solutions of   \eqref{pdeEq} as gradient flows in  $ (\PL,\WL) $ (see Paragraph \ref{BridgeToPde}). 

\paragraph{Geodesics.}
Let $ T \in (0, \infty) $. A curve $ (\mu_t)_{t\in [0,T ]} $ in a metric space $ (X,\mathrm{d}) $ is called \emph{geodesic} (between $ \mu_0 $ and $ \mu_T $) if $ \mathrm{d}(\mu_s, \mu_t) = (t-s)/T $ for all $ 0\leq s\leq t\leq T $.
In the following we show that $ (\PL,\WL) $ is a geodesic space, i.e.\ between each pair of measures there exists a geodesic. 

\bp{GeodLem}
Let $ \mu_0,\mu_T \in \PL $. 
Let $ \pi \in \Opt{\mu_0,\mu_T} $. 
Define the curve $ (\mu_t)_{t\in [0,T]} $ by 
\be{GeodLemEq}
 \mu_t  := (\pb_{\TdR \times \R}^1 \, , \,  (1-t)\pb_{\TdR \times \R}^2 + t \pb_{\TdR \times \R}^3 )_\# \pi, \ \ t\in [0,T].
\ee
Then $ (\mu_t)_{t} $ is   a geodesic.
Moreover, if in addition $ \mu_0 \ll \Leb{\TdR} $, then $  \mu_t  = (\pb^1_{\TdR}, (1-t)  \pb^2_{\TdR} + t \,  \mathrm{T}_{\mu_0}^{\mu_T})_\# \mu_0 $ and we also have that $ \mu_t \ll \Leb{\TdR}  $ for all $ t \in (0,T)$. 
\ep
\bpr
Note that for each $ t $, the disintegration of $ \mu_t $ with respect to $ \LebTd $ is given by $ \mu_t^x = ((1-t)\pb^1_{\R \times \R} + t\pb^2_{\R \times \R} ) _\# \pi^x $ for almost every $ x \in \Td $.
So that we know that $ (\mu_t^x)_{t}  $ is a geodesic in $ (\Pw,W_2)  $ for almost every $ x $ (see e.g.\ \cite[7.2.2]{ambgigsav}).
We infer
\be{GeodLemEq1}
\WL(\mu_s, \mu_t)^2 = 
			\frac{(t-s)^2}{T^2}	 \int_{\Td} W_2(\mu_0^x,\mu_T^x)^2 dx = 	\frac{(t-s)^2}{T^2}	\WL(\mu_0, \mu_T)^2.
\ee
The second claim follows from the  observation that $ \pi = (\pb^1_{\TdR},\pb^2_{\TdR},\T_{\mu_0}^{\mu_T}) _\# \mu_0$.
The third claim follows from the analogue statement in the Wasserstein space (see e.g.\ \cite[2.4]{ambsav}). 
\epr

\paragraph{Absolutely continuous curves.}
Let $ T \in (0, \infty)$ (or $ T=\infty $)  and let $ \mathrm{I} \subset  ( 0,T ) $ be a bounded (or unbounded) interval. A curve $  (\mu_t)_{t\in \mathrm{I}} $ in a metric space $ (X,\mathrm{d}) $ is called \emph{absolutely continuous} and we write $ \muc \in \AC{\mathrm{I};X} $ if there exists $ m \in \L^2(\mathrm{I}) $ (or $ m \in \L^2_\mathrm{loc}(\mathrm{I}) $ if $ \mathrm{I} $ is unbounded) such that 
\be{AbsContEq}
\mathrm{d}(\mu_s, \mu_t) \leq \int_s^t m(r ) dr  \ \ \ \ \forall \,  s,t\in \mathrm{I},  s\leq t.
\ee 
If $ (X,\mathrm{d}) $ is a Polish space, \cite[1.1.2]{ambgigsav} yields the existence of the \emph{metric derivative} $ |\mu'| \in \L^2(\mathrm{I}) $ (or $ |\mu'| \in \L^2_\mathrm{loc}(\mathrm{I}) $ if $ \mathrm{I} $ is unbounded) defined by 
\be{MetricDerEq}
|\mu'|(t) = \lim_{s\ra t} \frac{\mathrm{d}(\mu_s, \mu_t)}{|s-t|} \ \ \ \ \text{ for almost every } t\in \mathrm{I}.
\ee 
In the following we analyse absolutely continuous curves in $ (\PL,\WL) $ and 
show that some analogous results as in Wasserstein spaces (cf.\ \cite[Chapter 8]{ambgigsav}) hold true.

\bp{AbsContProp} 
\begin{itemize}
\item [\textbf{(A)}]
Let $ T \in (0, \infty) $ (or $ T=\infty $) and $ \muc\in \ACLT $. 
Then there exists $ v: (0,T) \times \Td \times \R \ra \R $ jointly measurable such that 
\begin{enumerate}[(i)]
	\item 
	$ \partial_t \mu_t + \partial_\t (\mu_t\,  v) = 0  $ in $ (0,T) \times \Td \times \R $ in the sense of distributions, i.e.
	\be{AbsContPropEq}
	\int_{ (0,T) \times   \Td \times \R} \Big(\partial_t \varphi_t(x,\t) + \partial_\t \varphi_t(x,\t) \, v_t^x(\t) \Big) d\mu_t(x,\t) dt =0    \ \  \forall \,  \varphi \in C_c^\infty((0,T) \times \Td \times \R),
	\ee
	\item
	$ \| v_t\|_{\L^2(\mu_t)}  \leq |\mu'|(t)			 $ for almost every $ t $,
	\item
	$ v_t \in \overline{\{	 \partial_\t \varphi\, | \, \varphi \in C^\infty_c(\TdR)			\}}^{\L^2(\mu_t)} $ for almost every $ t $,
	\item
	$ v_t^x \in \overline{\{	  \varphi' \, | \, \varphi \in C^\infty_c(\R)			\}}^{\L^2(\mu_t^x)} $ for almost every $ t $ and $ x $,
\end{enumerate} 
\item [\textbf{(B)}]
 Conversely, let $  (\mu_t)_{t\in (0,T)} \subset  \PL $  and let $ v   \in 	\L^2((0,T)\times \TdR \, ; \, \mu_t dt) 	 $ (or $ t \mapsto \|  v_t\|_{\L^2(\mu_t)}    \in 	\L_\mathrm{loc}^2((0,T)) 	 $ if $ T = \infty $).
 Suppose that \eqref{AbsContPropEq} holds. 
Then
$ \muc \in \ACLT $ and $ \| v_t\|_{\L^2(\mu_t)}  \geq |\mu'|(t)			 $ for almost every $ t $.
\end{itemize}
%
\ep
\bpr
Without restriction we can assume that $ T < \infty $, since otherwise we can exhaust $ (0,T) $ with bounded intervals. 

We now show \emph{\textbf{(A)}}.
We proceed analogously to the proof of \cite[8.3.1]{ambgigsav}.
Let $ \cT = \{	 \partial_\t \varphi\, | \, \varphi \in C^\infty_c((0,T) \times \TdR)			\} $. Define a linear map $ L:  \cT \ra \R $ by 
 	\be{AbsContPropEq1}
 L(\partial_\t \varphi ) := \int_{ (0,T) \times   \Td \times \R} \partial_t \varphi \, d\mu_t \, dt.
 \ee
Performing the very same steps as in the proof of  \cite[8.3.1]{ambgigsav}, we obtain 
 	\be{AbsContPropEq2}
|L(\partial_\t \varphi )| \leq    \| \,  |\mu'|\,  \|_{\L^2((0,T))} \  \| \partial_\t \varphi \|_{\L^2((0,T) \times \TdR \, ;\, \mu_t dt)},
\ee
which resembles equation (8.3.10)  in \cite{ambgigsav}. 
Note that we have tacitly used Lemma \ref{StabilityOptLem}.
Let $\overline{	\cT		} $ denote the closure of $ \cT $  with respect to  $ \| \cdot \|_{\L^2((0,T) \times \TdR \, ;\, \mu_t dt)} $.
Then, using the Riesz representation theorem,  \eqref{AbsContPropEq2} implies that there exists a unique $ v  \in \overline{	\cT		} $ such that 
\be{AbsContPropEq3}
L( w ) = \int_{ (0,T) \times   \Td \times \R} \,  v  \, w \,  d\mu_t \,  dt \ \ \ \  \forall \, w \in \overline{	\cT	}.
\ee
In particular, since we can take $ w = \partial_\t \varphi $ for $ \varphi \in C^\infty_c((0,T) \times \TdR)	 $, \eqref{AbsContPropEq3} yields \emph{(i)}.
Again, using  the same arguments as in \cite[8.3.1]{ambgigsav}, we obtain that for all intervals $ J \subset (0,T) $
 	\be{AbsContPropEq4}
\int_J  \| v_t\|_{\L^2(\mu_t)} ^2 dt \leq  \int_J |\mu'|^2(t)		dt,
\ee
which is equation (8.3.13)  in \cite{ambgigsav}. 
As $ J  $ was arbitrary, this implies \emph{(ii)}.
To show \emph{(iii)}, take $ (\varphi_n)_n \subset C^\infty_c((0,T) \times \TdR)	$ such that $ \partial_\t \varphi_n \ra v  $ in $ \L^2((0,T) \times \TdR \, ;\, \mu_t dt) $. 
Hence, the function $ t \mapsto \| \partial_\t \varphi_n(t, \cdot ) - v_t  \|_{\L^2(  \TdR \, ;\, \mu_t)}  $ converges to $ 0 $ in  $ \L^2((0,T) \, ;\, dt) $. 
This yields that, up to subsequences, $ t \mapsto \| \partial_\t \varphi_n(t, \cdot ) - v_t  \|_{\L^2(  \TdR \, ;\, \mu_t)}  $ converges to $ 0 $ point-wise almost everywhere.
Since $ \varphi_n(t, \cdot ) \in C^\infty_c( \TdR)	 $ for all $ t $, we conclude the proof of \emph{(iii)}.
In the same way, one proves the claim \emph{(iv)}.

Next we prove \textbf{\emph{(B)}}.
Let $ D \subset C^\infty_c((0,T) \times \R) $ be countable and dense with respect to $ \| \cdot \|_\infty $. 
Let $ \varphi \in D $.
Then \eqref{AbsContPropEq}  implies that 
	\be{AbsContPropEq5}
\int_{ \Td} \zeta(x)	\int_{ (0,T) \times \R} (\partial_t \varphi + \partial_\t \varphi \, v^x ) d\mu_t^x \, dt \,dx =0    \ \ \ \ \forall \,  \zeta  \in C_c^\infty( \Td).
\ee
Hence, there exists a $ \LebTd $-nullset $ \cN^\varphi $ such that 
	\be{AbsContPropEq51}
	\int_{ (0,T) \times \R} (\partial_t \varphi + \partial_\t \varphi \, v^x ) d\mu_t^x \, dt  =0    \ \ \ \ \forall \,  x  \in \Td \setminus \cN^\varphi.
\ee
Set $ \cN'= \cup_{\varphi  \in D } \cN^\varphi $. 
Moreover, the assumption that $ t \mapsto \|  v_t\|_{\L^2(\mu_t)}  \in 	\L^2((0,T)) 	 $  assures that there exists a further nullset $ \cN'' $ such that 
\be{AbsContPropEq6}
\int_{ (0,T) \times \R}  |v^x|^2 d\mu_t^x \, dt  < \infty     \ \ \ \forall \,  x  \in \Td \setminus \cN''.
\ee
Using that $ D $ is dense, the dominated convergence theorem yields  that 
	\be{AbsContPropEq52}
\int_{ (0,T) \times \R} (\partial_t \varphi + \partial_\t \varphi \, v^x ) d\mu_t^x \, dt  =0    \ \ \ \forall \,  \varphi  \in C_c^\infty( (0,T) \times \R) \ \forall \,  x  \in \Td \setminus ( \cN'\cup\cN'').
\ee
Therefore, for each $ x  \in \Td \setminus ( \cN'\cup\cN'') $, the pair $ \big((\mu_t^x)_t , (v_t^x)_t\big) $ fulfils the assumptions of the converse implication of \cite[2.5]{erbar}.
In particular, we obtain 
	\be{AbsContPropEq7}
W_2(\mu_s^x,\mu_t^x)^2 \leq (t-s) \int_{ s}^t  \|  v_r ^x \|_{\L^2(\mu_r^x)}^2 \, dr     \ \ \forall \,  0< s\leq t< T \ \  \forall  \,  x  \in \Td \setminus ( \cN'\cup\cN'') .
\ee
This inequality was shown at the end of the proof of  \cite[2.5]{erbar}.
\eqref{AbsContPropEq7} easily implies  that for all $ 0< s\leq t< T $
	\be{AbsContPropEq8}
\WL(\mu_s,\mu_t)^2 \leq (t-s) \int_{ s}^t  \|  v_r  \|_{\L^2(\mu_r)}^2 \, dr 
 \leq \left( \int_{ s}^t  \max \{1,\|  v_r  \|_{\L^2(\mu_r)}^2\} \, dr \right)^2.
\ee
We infer that $ (\mu_t)_{t\in (0,T )} $ is an absolutely continuous curve in $ (\PL,\WL) $.
Finally, the first inequality in \eqref{AbsContPropEq8} shows that $ \| v_t\|_{\L^2(\mu_t)}  \geq |\mu'|(t)			 $ almost everywhere.
\epr

The previous result introduced a few important objects that have to emphasized.
\bd{TangDef}
Let $ \mu \in \PL $,  $ T \in (0, \infty) $ (or $ T=\infty $) and $ \muc \in \ACLT $. 
Define
\begin{enumerate}[(i)]
	\item $ \TanL{\mu} := \overline{\{	 \partial_\t \varphi\, | \, \varphi \in C^\infty_c(\TdR)			\}}^{\L^2(\mu)}  $,	
	the \emph{tangent space} at $ \mu $,
	\item $ \Tanw{\mu^x} := \overline{\{	  \varphi' \, | \, \varphi \in C^\infty_c(\R)			\}}^{\L^2(\mu^x)}  $  for   $ x \in \Td $,
	\item $ v: (0,T) \times \Td \times \R \ra \R $ is called \emph{tangent velocity}  for $ (\mu_t)_{t}$ if 
	\begin{itemize}
		\item $ v   \in 	\L^2((0,T)\times \TdR \, ; \, \mu_t dt) 	 $ (or $ t \mapsto \|  v_t\|_{\L^2(\mu_t)}    \in 	\L_\mathrm{loc}^2((0,T)) 	 $ if $ T = \infty $),
		\item 
		$ \partial_t \mu_t + \partial_\t (\mu_t\,  v) = 0  $ in $ (0,T) \times \Td \times \R $ in the sense of distributions,
		\item
		$ v_t \in \TanL{\mu_t} $ for almost every $ t $.
	\end{itemize}
\end{enumerate}
\ed

The following lemma is an easy consequence of the above definition and can be proven exactly as in  \cite[Chapter 8.4]{ambgigsav}.

\bl{TangPropLem}
	\begin{enumerate}[(i)]
		\item \label{TangPropLemItem1} $\TanL{\mu} = \{w\in \L^2(\mu_t) \, | \,  \partial_\t (w \mu) =0  \}^\perp$ for $ \mu \in \PL $,
		where $ \partial_\t  $ is meant in the sense of distributions.
		\item \label{TangPropLemItem0} $v \in \TanL{\mu}$ if and only if $\| v \|_{L^2(\mu)} = \inf \{\| v +w\|_{L^2(\mu)} \, | \, w\in L^2(\mu), \partial_\t (w \mu) =0 \}$.
		\item \label{TangPropLemItem} Let $ \mu \in \PL $, $v \in \TanL{\mu}$ and $w\in \L^2(\mu)$ be such that $ \partial_\t (w \mu) =0$. Then 
		$\| v \|_{L^2(\mu)} =\| v +w\|_{L^2(\mu)}$ if and only if $\| w\|_{L^2(\mu)}=0$.
	\end{enumerate}
\el

We can summarize the previous results in the following statement. 
\bc{TangCor}
Let $ T \in (0, \infty] $. $ (\mu_t)_{t\in (0,T)} $ is absolutely continuous  in $ (\PL,\WL) $ if and only if there exists a tangent velocity $ v $ for $ (\mu_t)_{t} $.
Moreover, $ \| v_t\|_{\L^2(\mu_t)}  = |\mu'|(t)			 $ for almost every $ t $ and $ v $ is uniquely determined $ \Leb{(0,T)} $-a.e.
\ec
\bpr
Obviously,  Proposition \ref{AbsContProp} shows each claim except of the uniqueness result. 
Let $ w $ 
be an other tangent velocity for $ (\mu_t)_{t} $. 
Note that $ \partial_\t ((w_t-v_t) \mu_t) =0 $ for almost every $ t $.
Therefore, Lemma \ref{TangPropLem} (\ref{TangPropLemItem0}) implies that $ \| w_t \|_{L^2(\mu_t)} \leq \| w_t + (v_t -w_t)\|_{L^2(\mu_t)}  =  \| v_t \|_{L^2(\mu_t)}$  for almost every $ t $.
Analogously, applying Lemma \ref{TangPropLem} (\ref{TangPropLemItem0}) for $ v $ shows that $ \| v_t \|_{L^2(\mu_t)} =  \| w_t \|_{L^2(\mu_t)}= \| v_t + (w_t -v_t)\|_{L^2(\mu_t)} $ for almost every $ t $.
Using  Lemma \ref{TangPropLem} (\ref{TangPropLemItem}), this yields that $\| w_t -v_t\|_{L^2(\mu_t)} = 0$ for almost every $ t $. 
\epr

\paragraph{L-optimal maps vs. $ \boldsymbol{\TanL{\mu}} $.}
In the following we show that, if $ \mu\in\PLa $ and $ \nu \in \PL $, then $ \T_{\mu}^\nu - \pb^2 \in \TanL{\mu} $.
This will be a consequence of the following observation. 
\bl{TangentProjLem}
Let $ \mu\in\PL $ and $ w\in L^2(\mu) $. 
Then $ w \in \TanL{\mu}^\perp $ if and only if $ w(x,\cdot ) \in \Tanw{\mu^x}^\perp $ for almost every $ x \in \Td $.
\el
\bpr
The proof relies on Lemma \ref{TangPropLem} (\ref{TangPropLemItem1}).
Note that the same statements as in  Lemma \ref{TangPropLem} also hold for $ \Tanw{\mu^x} $ (\cite[Chapter 8.4]{ambgigsav}).
Therefore, the ``if''-part is trivial. 
To show the ``only if''-part, 
we apply the same arguments as in the proof of Proposition \ref{AbsContProp} \emph{(B)} to obtain a $ \LebTd- $nullset $ \cN $ such that 
for all $ x  \in \Td \setminus \cN $
	\be{TangentProjLemEq5}
\int_{  \R} \varphi'\  w(x,\cdot)\,  d\mu^x   =0    \ \ \ \forall \,  \varphi  \in C_c^\infty( \R).
\ee
We conclude that $ w(x,\cdot ) \in  \{w\in \L^2(\mu_t) \, | \,  \partial_\t (w \mu) =0  \}=\Tanw{\mu^x}^\perp $ for almost every $ x \in \Td $.
\epr
\bc{LOptTangLem}
Let $ \mu\in\PLa $ and $ \nu \in \PL $.
Then $ \T_{\mu}^\nu - \pb^2 \in \TanL{\mu} $.
\ec
\bpr
It is enough to show that for all $  w \in \TanL{\mu}^\perp  $
\be{LOptTangLemEq}
 \int_{\TdR} (\T_{\mu}^\nu - \pb^2) \,  w\  d \mu = 0  .
 \ee
\cite[8.5.2]{ambgigsav} states that $   \T_{\mu}^\nu(x, \cdot ) - \pb^2 = \T_{\mu^x}^{\nu^x} - \Id_{\R} \in \Tanw{\mu^x}  $ for almost every $ x \in \Td $. 
Therefore, Lemma \ref{TangentProjLem} implies that $ \int_\R (\T_{\mu}^\nu - \pb^2)(x,\cdot) w(x,\cdot) d \mu^x = 0  $ for almost every $ x $, which immediately implies \eqref{LOptTangLemEq}.
\epr

\paragraph{$ \boldsymbol{\ACLT} $ vs. $ \boldsymbol{\ACwT} $.}
Here we show that a curve $ (\mu_t)_t $ is absolutely continuous in $ (\PL,\WL) $ if and only if $ (\mu_t^x)_t $ is absolutely continuous in $ (\Pw,W_2) $ for almost every $ x $.
\bl{AbsContFiberLem}
Let $ T \in (0, \infty) $ or $ T= \infty $,  $ (\mu_t)_{t \in (0,T)} \subset \PL  $ and $ v   \in 	\L^2((0,T)\times \TdR \, ; \, \mu_t dt) 	 $ (or $ t \mapsto \|  v_t\|_{\L^2(\mu_t)}    \in 	\L_\mathrm{loc}^2((0,T)) 	 $ if $ T = \infty $).
Then $ (\mu_t)_t  \in \ACLT$ and $ v $ is the tangent velocity for $ (\mu_t)_t  $ if and only if for almost every $ x \in \Td $, 	$ (\mu_t^x)_t  \in \ACwT$ and $ v^x $ is the tangent velocity for $ (\mu_t^x)_t$ in the Wasserstein sense, i.e.\
there exists a $ \Leb{(0,T)} $-nullset $  \cN_x $ such that 
\begin{enumerate}[(i)] 
	\item 
	$ v^x   \in 	\L^2((0,T)\times \R \, ; \, \mu^x_t dt)	 $ (or $ t \mapsto \|  v_t^x\|_{\L^2(\mu_t^x)}    \in 	\L_\mathrm{loc}^2((0,T)) 	 $ if $ T = \infty $),
	\item 	$ \partial_t \mu_t^x + \partial_\t (\mu_t^x \,  v^x) = 0  $ in $ (0,T)  \times \R $ in the sense of distributions,
	\item $ v_t^x \in \Tanw{\mu_t^x} $ for all $ t \in (0,T) \setminus \cN_x $.
\end{enumerate} 
In particular, $ |\mu'|^2(t) = \| v_t\|_{\L^2(\mu_t)}^2   = \int_{  \Td}  \| v_t^x\|_{\L^2(\mu_t^x)}^2 dx =\int_{  \Td	} |(\mu^x)'|^2(t) \,	dx$  
for almost every $ t $.
\el
\bpr
Assume $ (\mu_t)_t  \in \ACLT$ with tangent velocity $ v  $. 
\emph{(i)} follows from the corresponding integrability condition on $ v $ being the tangent velocity of $ \muc $.
\emph{(ii)} was shown in the proof of Proposition \ref{AbsContProp} \emph{(B)}.
\emph{(iii)} follows from Proposition \ref{AbsContProp} \emph{(A)}.  
By \cite[2.5]{erbar}, these  facts  imply that 
\be{arg1}
 (\mu_t^x)_t  \in \ACwT \ \ \  \text{ and } \ \ \   \| v_t^x\|_{\L^2(\mu_t^x)}  =  |(\mu^x)'|(t)		\  \text{ for almost every } t .
\ee
%

Conversely,
it is an easy observation that (ii) implies that $ \partial_t \mu_t + \partial_\t (\mu_t \,  v) = 0  $ in the sense of distributions. 
Hence, Proposition \ref{AbsContProp} \emph{(B)} yields that $ (\mu_t)_t  \in \ACLT$ and $ \| v_t\|_{\L^2(\mu_t)}  \geq  |\mu'|(t)			 $ for almost every $ t $.
It remains to show that $ v $ is the tangent velocity for $ (\mu_t)_t $.
An easy application of Fubini's theorem  shows that \emph{(iii)} can be reformulated as follows:
For almost every $ t $, there exists a $ \LebTd $-nullset $  \cN_t $ such that 
 $ v_t^x \in \Tanw{\mu_t^x} $ for all $ x \in \Td \setminus \cN_t $.
Using this formulation, we can argue in the same way as in Corollary \ref{LOptTangLem} to conclude that $ v_t \in \TanL{\mu_t} $ for a.e.\ $ t $,
which shows that 
$ v $ is the tangent velocity for $ (\mu_t)_t $.
\epr 

\paragraph{Infinitesimal behaviour.}
The goal of this paragraph is to show differentiability of $ \WL $ along absolutely continuous curves.
We start with the following observation, which, again, is also true in the analogue setting of the Wasserstein distance. 

\bl{DirOptMapsLem} 
Let $ T \in (0, \infty] $ and $ \muc \in \ACLT $ with tangent velocity $ v $. 
Suppose that $ \muc \subset \PLa $.
Then 
\be{DirOptMapsLemEq}
\lim_{h\ra 0} \left\|	\frac{1}{h}(\T_{\mu_t}^{\mu_{t+h}} - \pb^2) - v_t \right\|_{\L^2(\mu_t)} = 0   \ \ \  \text{ 	for almost every 		} t \in (0,T).
\ee
\el
\bpr
Let $ s_h^t := 	\frac{1}{h}(\T_{\mu_t}^{\mu_{t+h}}  - \pb^2)$.
By Lemma \ref{AbsContFiberLem}, we have that   $ (\mu_t^x)_t  \in \ACwT$ with Wasserstein tangent velocity $ v^x $ for almost every $ x $.
Therefore, we can apply \cite[8.4.6]{ambgigsav} to see that for almost every $ x $ there exists a nullset $ \cN_x  $ such that 
\be{DirOptMapsLemEq1}
\lim_{h\ra 0} \left\|	s_h^t(x,\cdot) - v_t^x \right\|_{\L^2(\mu_t^x)} = 0   \ \ \  \text{ 	for all } t \in (0,T)\setminus \cN_x.
\ee
As above, using Fubini's theorem, we can reformulate \eqref{DirOptMapsLemEq1} in such a way that for almost every $ t $, there exists a nullset $ \cN_t  $ such that 
\be{DirOptMapsLemEq2}
\lim_{h\ra 0} \left\|	s_h^t(x,\cdot) - v_t^x \right\|_{\L^2(\mu_t^x)} = 0   \ \ \  \text{ 	for all } x \in \Td \setminus \cN_t.
\ee
In particular, this shows that for almost every $ t $, $ x \mapsto \|	s_h^t(x,\cdot) \|_{\L^2(\mu_t^x)}^2  $ converges to $ x \mapsto \| v_t^x \|_{\L^2(\mu_t^x)}^2  $ point-wise almost everywhere. 
However, since for almost every $ t $
\be{DirOptMapsLemEq3}
\int_{  \Td}   \left\|	s_h^t(x,\cdot)  \right\|_{\L^2(\mu_t^x)}^2 dx =  \frac{1}{h^2} \WL(\mu_t,\mu_{t+h}) \Ra |\mu'|^2(t) = \int_{  \Td} \| v_t^x \|_{\L^2(\mu_t^x)}^2 dx,
\ee
we even have that $ x \mapsto \|	s_h^t(x,\cdot) \|_{\L^2(\mu_t^x)}^2  $ converges to $ x \mapsto \| v_t^x \|_{\L^2(\mu_t^x)}^2  $ in $ \L^1 ( \Td) $  for almost every $ t $.
Hence, for each $ h $, the function $ x \mapsto \|	s_h^t(x,\cdot) -v_t^x \|_{\L^2(\mu_t^x)}^2  $ is majorized by the function $ x \mapsto 4 (\|	s_h^t(x,\cdot) \|_{\L^2(\mu_t^x)}^2  + \| v_t^x \|_{\L^2(\mu_t^x)}^2 ) $, which is a converging sequence in $ \L^1 ( \Td) $.
Therefore, we can apply the (generalized) dominated convergence theorem to obtain that for almost every $ t $
\be{DirOptMapsLemEq4}
\lim_{h\ra 0} \left\|	s_h^t - v_t \right\|_{\L^2(\mu_t)}^2 = 
\lim_{h\ra 0}  \int_{  \Td} \left\|	s_h^t(x,\cdot) - v_t^x \right\|_{\L^2(\mu_t^x)}^2 dx 
=
\int_{  \Td} \lim_{h\ra 0}   \left\|	s_h^t(x,\cdot) - v_t^x \right\|_{\L^2(\mu_t^x)}^2 dx 
=
0,
\ee
which concludes the proof of this lemma. 
\epr

\bp{DerivativeAbsContProp}
Let $ T \in (0, \infty] $ and $ \muc \in \ACLT $ with tangent velocity $ v $. 
Let $ \muc \subset \PLa $ and $ \s \in \PL $.
Then 
\be{DerivativeAbsContPropEq}
\frac{d}{dt}\WL(\mu_t,\s)^2 = 2  \int_{  \TdR	} (\pb^2- \T_{\mu_t}^{\sigma})\, v_t \ d\mu_t    \ \ \  \text{ 	for almost every 		} t \in (0,T).
\ee
\ep
\bpr
As above, the proof relies on the analogous result for $ (\mu_t^x)_t  $ and the dominated convergence theorem.
Let for all $ t\in(0,T) $ and $ h>0 $
\be{DerivativeAbsContPropEq1}
f_h^t: x \to \frac{1}{h^2} W_2(\mu_t^x , \mu_{t+h}^x )^2 + 4 (  W_2(\mu_t^x , \s^x )^2 +	W_2( \mu_{t+h}^x,\s^x )^2 		) .
\ee
It will turn out that $ f_h^t $ is the majorizing sequence that we need. 
Thus, we need to show that $ (f_h^t)_h $ converges in $ \L^1(\Td) $ for a.e.\ $ t $.
Indeed,  we observe that as $ h \ra 0 $ (again, after an application of Fubini's theorem) for almost every $ t $
\be{DerivativeAbsContPropEq2}
f_h^t(x) \Ra |(\mu^x)'|^2(t) + 8  W_2(\mu_t^x , \s^x )^2  		) =: f^t(x)   \ \ \  \text{ 	for almost every 		} x \in \Td.
\ee
Moreover, for almost every $ t $
\begin{align}\label{DerivativeAbsContPropEq3}
	\begin{split}
		\|f_h^t \|_{\L^1(\Td)}&= \frac{1}{h^2} \WL (\mu_t , \mu_{t+h} )^2 + 4 (  \WL(\mu_t , \s )^2 +	\WL( \mu_{t+h},\s )^2 		) \\
		&\Ra |\mu'|^2(t) + 8  \WL(\mu_t , \s )^2 = \|f^t \|_{\L^1(\Td)}.
	\end{split}
\end{align}
\eqref{DerivativeAbsContPropEq2} and \eqref{DerivativeAbsContPropEq3} show that $ \lim_{h\ra 0}f_h^t = f^t $ in $ \L^1(\Td) $ for a.e.\ $ t $.

Note that from   \cite[8.4.7]{ambgigsav} we get that  for almost every $ t $ and $ x $
\be{DerivativeAbsContPropEq4}
\frac{d}{dt} W_2(\mu_t^x , \s^x )^2 = 2  \int_{  \R	} (\Id_{\R}- \T_{\mu_t^x}^{\sigma^x})\, v_t^x \ d\mu_t^x.
\ee
Further, as a consequence of the triangle inequality and Young's inequality, we observe that for all $ t,h $ and $ x $
\be{DerivativeAbsContPropEq5}
\frac{1}{h^2} (	W_2( \mu_{t+h}^x,\s^x )^2 	- W_2(\mu_t^x , \s^x )^2	)  \leq \frac{1}{2} f_h^t(x).
\ee
Therefore, the dominated convergence theorem yields \eqref{DerivativeAbsContPropEq}.
\epr

\subsection{\texorpdfstring{Gradient flows in $ (\PL,\WL) $ for $\lambda$-convex functionals}{Gradient flows in  (P\_2\textasciicircum L ,W\textasciicircum L)}}\label{GradFlowSec}

In this section we introduce the notion of a \emph{subdifferential} for a certain class of functionals in $ \PL $. 
Then we  define \emph{gradient flows} in $ \PL $ for such functionals and  prove in Theorem \ref{ExistenceThm} their existence, uniqueness and some  properties.
 
In this paper, we only consider functionals that satisfy the following convexity property (cf.\  \cite[4.0.1]{ambgigsav}). 
\bd{StrongConvexDef}
Let $ (X, \mathrm{d}) $ be a Polish space. Then $ \phi : X \ra (-\infty,\infty] $ is called \emph{strongly $ \lambda$-convex} if   $ \lambda \in \R $ and for all $ \s, \mu_0,\mu_1 \in D(\phi) $ there exists a curve $ (\gamma_t)_{t\in[0,1]} $ with $ \gamma_0= \mu_0 $, $ \gamma_1= \mu_1 $ such that for all $ 0 <  \tau  <  	\tfrac{1}{\lambda^-  		}	 $ (with the convention that $ 1/0=\infty $), the functional 
\be{StrongConvexDefEq}
\Phi(\tau, \sigma ; \cdot ) := \frac{1}{2\tau } \mathrm{d}(\cdot, \sigma)^2 + \phi (\cdot)
\ee
 is $ (\frac{1}{\tau } + \lambda)$-convex  along $ (\gamma_t)_t $,
 i.e.\ for all $ t\in[0,1] $
\be{StrongConvexDefEq1}
\Phi(\tau, \sigma ; \gamma_t ) \leq (1-t ) \Phi(\tau, \sigma ; \gamma_0 )  + t \, \Phi(\tau, \sigma ; \gamma_1 ) -   
 \frac{1}{2}\left(\frac{1}{\tau } + \lambda\right) t ( 1-t) \, \mathrm{d}(\mu_0, \mu_1)^2 .
\ee
\ed

However, in most of the cases, the following weaker form of convexity will be enough.
\bd{ConvexDef}
$ \phi : X \ra (-\infty,\infty] $ is called \emph{$ \lambda$-convex} if  $ \lambda \in \R $ and for all $  \mu_0,\mu_1 \in D(\phi) $ there exists 
a geodesic $ (\mu_t)_{t\in[0,1]} $ such that $ \phi  $ is $  \lambda $-convex  along $ (\mu_t)_{t}$.
\ed
In our case, i.e.\ if $ X = \PL $, $ (\mu_t)_{t\in[0,1]} $ will always be  the geodesic induced by some 
$ \pi \in \Opt{\mu_0,\mu_1} $ as in \eqref{GeodLemEq}.

\paragraph{Subdifferential calculus.}
Instead of working with  gradients in $ \PL $, we prefer to work with (strong) subdifferentials. 
On the one hand, their properties are easier to verify (as it only needs lower bounds), while on 
 the other hand  they are enough to build the bridge to \eqref{pdeEq}. 
\bd{SubdiffDef}
Let $ \phi : \PL \ra (-\infty,\infty] $ be proper\footnote{This means that $ \phi (\mu) > -\infty  $ for all $ \mu \in \PL $ and there exists $ \mu \in \PL $ such that $ \phi (\mu) < \infty  $.}, $ \lambda $-convex and $ \WL $-l.s.c.
Let $ \mu \in D(\phi) \cap \PLa $ and $ \xi \in \L^2(\mu) $.
Then we say that $ \xi $ \emph{belongs to the subdifferential of $ \phi $ at $ \mu $} and we write $ \xi \in \partial\phi(\mu) $ if 
\be{SubdiffDefEq} 
\phi(\nu ) - \phi(\mu ) \geq \int_{  \TdR } \xi \, ( \T_{\mu}^\nu  -\pb^2) \ d\mu + \frac{\lambda}{2} \WL(\mu,\nu)^2 \ \ \ \forall \, \nu \in D(\phi). 
\ee
Further, we say that $ \xi \in \partial\phi(\mu) $ is a \emph{strong subdifferential of $ \phi $ at $ \mu $} if 
\be{SubdiffDefEq1} 
\phi((\pb^1,T)_\#\mu ) - \phi(\mu ) \geq \int_{  \TdR } \xi \, ( \T  -\pb^2) \ d\mu + o(\| 	\T  -\pb^2 \|_{\L^2(\mu)}) \ \ \text{ as } \| 	\T  -\pb^2 \|_{\L^2(\mu)} \ra 0.
\ee
\ed

\bl{TangentStronSubdiff} 
Let $ \phi : \PL \ra (-\infty,\infty] $ be proper, $ \lambda $-convex and $ \WL $-l.s.c.
Let $ \mu \in D(\phi) \cap \PLa $ and  $ \xi \in \partial\phi(\mu) \cap \TanL{\mu}$.
Then  $ \xi  $ is a strong subdifferential of $ \phi $ at $ \mu $.
\el
\bpr
The proof follows  the same lines as in the Wasserstein case (see \cite[3.2]{erbar}).
Therefore, we omit the details. 
\epr

\bd{MetricSlopeDef} 
Let $ (X,\mathrm{d}) $ be a Polish space. 
Let $ \phi : X \ra (-\infty,\infty] $ be proper and $ \mathrm{d} $-l.s.c.
Then the \emph{metric slope}  $ |\partial\phi|: D(\phi) \ra [0,\infty] $ is defined by 
\be{MetricSlopeDefEq}
|\partial\phi|(\mu ) = \limsup_{\nu \ra \mu } \frac{(\phi ( \mu) - \phi ( \nu))^+}{\mathrm{d}(\mu,\nu)}.
\ee
\ed

Next we show that $ \lambda $-convex functionals are differentiable almost everywhere along curves in $ \ACLT $ and compute the derivative.

\bl{ChainRuleLem} 
Let $ \phi : \PL \ra (-\infty,\infty] $ be proper,  $ \lambda $-convex and $ \WL $-l.s.c.
Let $ T \in (0, \infty] $ and $ \muc \in \ACLT $ with tangent velocity $ v $.
Suppose that 
\be{ChainRuleLemEq}
\int_s^t |\partial\phi|(\mu_r ) |\mu'|(r) \ dr < \infty \ \ \ \ \forall  0 < s< t <T.
\ee
Then  
\begin{enumerate}[(i)]
	\item $t \to \phi(\mu_t)$ is absolutely continuous,
    \item there exists a $\Leb{(0,T)}$-nullset $\cN$ such that for all $\xi \in \partial \phi(\mu_t)$
	\begin{align}\label{ChainRuleLemEq1}
	\frac d{dt} \phi (\mu_t) = \int_{\TdR} \xi  \, v_t \ d\mu_t \ \text{ for all } t\in (0,T) \setminus \cN.
	\end{align}
\end{enumerate}
\el

\bpr 
\emph{(i)} is the  content of   \cite[2.4.10]{ambgigsav}.
To show \emph{(ii)},
let $ \cN  $ be such that 
\eqref{DirOptMapsLemEq} holds, $ \frac{d}{dt}  \phi(\mu_t) $ exists and $ |\partial\phi|(\mu_t ) < \infty $ for all $ t\in (0,T) \setminus \cN $. 
From here we proceed as in   \cite[10.3.18]{ambgigsav}. 
\epr

\paragraph{Gradient flows in  $ (\PL,\WL) $.}
We are now able to define the notion of gradient flows in $ (\PL,\WL) $.
\bd{GradFlowDef}
Let $ \phi : \PL \ra (-\infty,\infty] $ be proper,  $ \lambda $-convex and $ \WL $-l.s.c.
Let $ T \in (0, \infty] $ and $ \muc \in \ACLT $ with tangent velocity $ v $.
Then $ \muc  $ is called \emph{gradient flow for $ \phi $},
if 
\be{GradFlowDefEq}
-v_t \in \partial \phi (\mu_t) \ \ \ \text{ for almost every }t\in(0,T).
\ee
Further, $ \mu_0 \in \PL $ is called \emph{initial value} of $ \muc $ if $ \lim_{t\ra0}\WL(\mu_t,\mu_0) = 0  $.
\ed
Let us first note that, as in  the Wasserstein case, gradient flows in $ (\PL,\WL) $ are equivalent to the solutions of a system of \emph{evolution variational inequalities (E.V.I)}.
\bl{EVI}
Let $ \phi : \PL \ra (-\infty,\infty] $ be proper,  $ \lambda $-convex and $ \WL $-l.s.c.  Let $ T \in (0, \infty] $ and $ \muc \in \ACLT $.
Then $ \muc $ is a gradient flow for $ \phi $ if and only if for all $ \nu \in D(\phi) $ there exists a $\Leb{(0,T)}$-nullset $\cN_\nu$ such that for all $t\in (0,T) \setminus \cN_\nu$
\begin{align}\label{EVIEq}
 \frac{1}{2}\frac{d}{dt}\WL(\mu_t,\nu)^2 \leq \phi(\nu) - \phi(\mu_t) - \frac{\lambda}{2}\WL(\mu_t,\nu)^2.
\end{align}
\el
\bpr
Again, the proof consists in adapting the analogous proof in the Wasserstein case (see \cite[11.1.4]{ambgigsav}), which is based on Proposition \ref{DerivativeAbsContProp}. We omit the details. 
\epr

In the following theorem we obtain  the existence of gradient flows and further properties such as uniqueness, an energy identity and a regularisation estimate.   
The result is limited to the case, when the functional $ \phi $ is  proper, strongly $ \lambda $-convex, $ \WL $-l.s.c and  \emph{coercive}, 
where we say that a functional 
$ \phi : X \ra (-\infty,\infty] $ on a Polish space $ (X, \mathrm{d}) $ is coercive if  
there exists $ \mu^*\in X $ and $ r^*>0  $ such that 
\be{Coercivity}
\inf \{	\phi(\nu) \, | \, 	\nu \in X \, , \, \mathrm{d}(\nu,\mu^*) \leq r^*			\} > -\infty \ \ \  \text{ (cf.\  \cite[(2.4.10)]{ambgigsav})}.
\ee

\bt{ExistenceThm}
Let $ T \in (0, \infty] $ and $ \phi : \PL \ra (-\infty,\infty] $ be proper, strongly $ \lambda $-convex,  $ \WL $-l.s.c and coercive.
Then:
\begin{enumerate}[(i)]
	\item \emph{(Existence)} For each $ \mu_0 \in \overline{D(\phi)} $, there exists a gradient flow for $ \phi $ with initial value $ \mu_0 $.
	\item \emph{($ \lambda $-contractivity and uniqueness)} 
	Let  $ (\mu_t)_t $ and $ (\nu_t)_t $ be gradient flows  for $ \phi $ with initial value $ \mu_0 \in \overline{D(\phi)} $  and $\nu_0 \in \overline{D(\phi)} $, respectively. Then, for all $ t \in (0,T) $
	\be{Contractivity}
	\WL(\mu_t, \nu_t) \leq e^{-\lambda t }	\WL(\mu_0, \nu_0) .
	\ee
	In particular, for each  $ \mu_0 \in \overline{D(\phi)} $, the gradient flow for $ \phi $ with initial value $ \mu_0 $ is unique.
	\item \emph{(Energy identity)}
	Let  $ (\mu_t)_t $ be the gradient flow  for $ \phi $ with initial value $ \mu_0 \in D(\phi) $, then  for all $ t \in (0,T)$ 
\be{EnergyIdEq}
\phi(\mu_t) - 		\phi(\mu_0) + \frac 12 \int_0^t\big(|\partial \phi|^2(\mu_s)  + |\mu'|^2(s)\big)\,ds =0.
\ee
	\item \emph{(Monotonicity along gradient flows)}
		Let  $ (\mu_t)_t $ be the gradient flow  for $ \phi $ with initial value $ \mu_0 \in \overline{D(\phi)} $, then  for almost every $ t \in (0,T)$
			\begin{align}\label{MonotonicityLEm}
		\frac d{dt} \phi (\mu_t) = -\|v_t\|_{\L^2(\mu_t)}^2 .
		\end{align}
		\item \emph{(Regularization estimate)}
		Let  $ (\mu_t)_t $ be the gradient flow  for $ \phi $ with initial value $ \mu_0 \in \overline{D(\phi)} $, then  for all $ t \in (0,T)$ and all $ \nu \in D(\phi) $
	\begin{align}\label{ReguEstEq}
	\phi(\mu_t) \leq 
	\begin{cases}
	\phi(\nu) + \frac \lambda {2 ( e^{\lambda t }-1)} \WL(\mu_0,\nu)^2   &: \lambda \neq 0,\\
	\phi(\nu) + \frac 1 {2  t}  \WL(\mu_0,\nu)^2   &: \lambda = 0.
	\end{cases}
	\end{align}
\end{enumerate}
\et
\bpr
Again, we benefit from the work that was done in \cite{ambgigsav}.

For $ \mu_0 \in \overline{D(\phi)}  $, we introduce the following implicit Euler scheme. 
Let $ \tau >0 $.
Define recursively:
\begin{align}\label{SchemeEq}
\begin{cases}
\mu_0^\tau := \mu_0, \\
\mu_n^\tau 
\in \argmin\limits_{\nu \in \PL}\left( \phi( \nu ) + \frac { 1}{2\tau}\WL(\mu_{n-1}^\tau,\nu)^2  \right) \text{ for } n \in \N.
\end{cases}
\end{align}
\cite[2.2.2]{ambgigsav} shows that this scheme is well-defined. 
Define the piecewise constant interpolating trajectory $ (\bar{\mu}_t^\tau)_{t \in [0,T]} $  by
\begin{align}\label{SchemeEq1}
\begin{cases}
\bar{\mu}_0^\tau := \mu_0,  & \\
\bar{\mu}_t^\tau := \mu_n^\tau & \text{ for } t \in ((n-1)\tau, n \tau] \text{ for all } n \in \N \text{ such that }n \tau \leq T.
\end{cases}
\end{align}
Then  \cite[4.0.4]{ambgigsav} yields the convergence of this scheme with respect to $ \WL $ towards a curve $ \muc \in \ACLT $ with initial value $ \mu_0 $, which solves \eqref{EVIEq} and satisfies \emph{(ii)}.
In addition,  Lemma \ref{EVI} yields \emph{(i)}.

\cite[4.0.4]{ambgigsav} shows that the gradient flow $ \muc $ is a so-called minimizing movement (see  \cite[2.0.6]{ambgigsav} for the definition). 
Hence,   \cite[2.3.3]{ambgigsav} implies \emph{(iii)}. (Note that in our case the object $ |\partial^- \phi| $ from this theorem is just $ |\partial \phi| $ and that the assumption that $ |\partial \phi| $ is a \emph{strong upper gradient} is also fulfilled by \cite[2.4.10]{ambgigsav}.)

\emph{(iv)}
follows from the chain rule given in  Lemma \ref{ChainRuleLem}.

\emph{(v)}
follows from \cite[4.3.2]{ambgigsav} and \cite[(3.1.1)]{ambgigsav}. 
\epr

\subsection{Local McKean-Vlasov equation}\label{McKean}
In this section we apply Theorem \ref{ExistenceThm} to a functional 
$ \cF $ that will be of the form
\begin{align}\label{McKeanVlasovFunctionalF}
\cF(\mu) := 
\cS(\mu)  + \cW(\mu) + \cV(\mu ),
\end{align}
where $ \cS, \cW $ and $ \cV $ are called \emph{entropy}, \emph{interaction energy} and \emph{potential energy}, respectively. 
In order to apply Theorem \ref{ExistenceThm} for $ \cF $, we show separately that each of its  summands $ \cS, \cW $ and $ \cV $  are well-defined, proper, strongly $ \lambda $-convex,  $ \WL $-l.s.c and coercive in the Lemmas \ref{EntropyProp}, \ref{InteractionProp} and \ref{PotentialProp}, respectively.
(It will turn out that $ \cF $ is trivially proper.)
Moreover, we compute  a directional derivative of $ \cF $ (Proposition \ref{DirDerProp}),  analyse the subdifferential of $ \cF $ (Proposition \ref{SubdiffCharacF}) and derive a variational characterisation for gradient flows for $ \cF $ (Theorem \ref{Lyapunov}), which will be a key fact in the forthcoming chapters in this paper.
Finally, we show in Theorem \ref{pde} the equivalence of the gradient flow for $ \cF $ and the weak solution to the partial differential equation \eqref{pdeEq}.


\paragraph{Entropy.}
Define the \emph{entropy} $\cS:\PL \ra (-\infty, \infty]$  by
\begin{align}\label{EntropyFunctional}
\cS(\mu) := 
\begin{cases}
\int_{\TdR}  \log(\rho) d\mu &: \mu \ll \Leb{\TdR}, \ \mu= \rho \, \Leb{\TdR}, \\
\infty   &: \text{else}.
\end{cases} 
\end{align}
A very useful observation is that  for each $ \mu \in \PL $
  \begin{align}\label{EntropyFunctional2}
  \cS(\mu) = 
  \int_{\TdR}  S_1(\mu^x) dx,
  \end{align}
where $ S_1 : \Pw \ra (-\infty, \infty] $ is the entropy functional on $ \Pw $, i.e.
\begin{align}\label{EntropyFunctionalOnR}
S_1(\mu^x) := 
\begin{cases}
\int_{R}  \log(\rho^x) d\mu^x &: \mu^x \ll \Leb{\R}, \ \mu^x= \rho^x \,  \Leb{\R}, \\
\infty   &: \text{else}.
\end{cases} 
\end{align}
This fact will simplify our analysis, since we   benefit from the already known results for $ S_1 $; see e.g.\ in \cite{ambgigsav}.
In the following lemma we show that Theorem \ref{ExistenceThm} is applicable for $ \cS $. 
\bl{EntropyProp}
\begin{enumerate}[(i)]
	\item \emph{(Well-defined)}
	Let $ \mu \in \PL $ and $ \e>0 $. 
	Then there exists $ C_\e > 0  $ such that 
	$ \cS(\mu)  \geq \, - C_\e - \e \int |\t|^2 d\mu \ (>-\infty)$.
	\item \emph{(Coercivity)}
	For all $ r>0 $ we have 
	\be{CoercivityEntropy}
	\inf \left\{	\cS(\nu) \, \Big| \, 	\nu \in \PL \, , \, \int |\t|^2 d\nu \leq r			\right\} > -\infty .
	\ee
	 In particular, $ \cS  $ is coercive.
	 \item \emph{($ \WL$-l.s.c)}
	Let $ (\mu_n)_{n\in\N}  $ be such that $ \sup_n  \int |\t|^2 d\mu_n<\infty$ and $ \mu_n \hra \mu\in\PL $. Then
	 \be{LSCEntropy}
	 \liminf_{n\ra \infty } 	\cS(\mu_n) \geq \cS(\mu).
	 \ee
	 In particular, $ \cS  $ is $ \WL $-l.s.c.
	 \item \emph{(Strong $ 0$-convexity)}
	 $ \cS $ is strongly $ 0$-convex.
\end{enumerate}
\el
\bpr
The corresponding statement for $ S_1 $ (see \cite[(29)]{jko}) and \eqref{EntropyFunctional2} imply \emph{(i)}.

\emph{(ii)} is an immediate consequence of \emph{(i)}.

To show \emph{(iii)},
set $ \nu := e^{-|\t| - \b } d\t dx \in \PL$, where $ \b >0  $ is a normalization constant. 
Recall the definition of the relative entropy given in  \eqref{RelEntropyFunctional}.
Then for $ \mu \in \PL $
\begin{align}\label{EntropyRelEntropy}
\cS(\mu) = \cH(\mu \, | \,  \nu ) - \tilde{\cV} (\mu) ,
\end{align}
where $ \tilde{\cV} (\mu) := \int (|\t| + \b) \, d\mu $.
Since $ \sup_n  \int |\t|^2 d\mu_n<\infty$, \cite[5.1.7]{ambgigsav} implies that
\be{LSC1}
\lim_{n\ra \infty } 	\tilde{\cV} (\mu_n) = \tilde{\cV} (\mu).
\ee
And by the dual representation of $  \cH $ (see  \cite[9.4.4]{ambgigsav}), we have that $  \cH(\cdot \, | \,  \nu ) $ is the supremum of functionals that are continuous with respect to weak convergence. 
Hence, $  \cH(\cdot \, | \,  \nu ) $ is lower semi-continuous  with respect to weak convergence. 
This fact together with \eqref{LSC1}
yields \emph{(iii)}.

It remains to prove \emph{(iv)}.
Let $ \s,\mu_0,\mu_1 \in D(\cS) $ and $ \Phi $ be as in \eqref{StrongConvexDefEq} for the functional $ \cS $.
Analogously, define  $ \Phi_1(\tau , \s^x ;\cdot ) = \tfrac{1}{2\tau} W_2(\s^x, \cdot ) + S_1(\cdot) $.
Then we observe that for all $ \mu \in \PL $
\be{StrongConvexityEntropyEq}
\Phi(\tau , \s ;\mu ) = \int_{\Td} \Phi_1(\tau , \s^x ;\mu^x ) \, dx.
\ee
Moreover, we show at the end of this proof that there exits a measure $ \omega \in \cM_1^\L(\Td\times\R^3) $ such that 
for almost every $ x \in \Td $
\be{StrongConvexityEntropyEq1}
 (\pb_{\R^3}^{1,2}) _\# \omega^x \in \OptW{\s^x,\mu_0^x} \quad \text{and} \quad
 (\pb_{\R^3}^{1,3}) _\# \omega^x \in \OptW{\s^x,\mu_1^x}.
\ee
Set for all $ t \in [0,1] $
\be{StrongConvexityEntropyEq2}
\gamma_t = \left((1-t )\, \pb_{\R^3}^{2}  + t\, \pb_{\R^3}^{3}  \right) _\# \omega^x \, dx \in \cM_1(\Td\times\R\times \R) .
\ee 
Then,  \cite[9.3.9]{ambgigsav} and \cite[9.2.7]{ambgigsav} show that 
for almost every $ x \in \Td $ and for all $ t \in [0,1] $
\be{StrongConvexityEntropyEq3}
\Phi_1(\tau , \s^x ;\gamma_t^x )  \leq (1-t) \Phi_1(\tau , \s^x ;\gamma_0^x ) + t \, \Phi_1(\tau, \sigma^x ; \gamma_1^x ) -   
\frac{1}{2\tau }  t ( 1-t) W_2(\mu_0^x, \mu_1^x)^2 .
\ee
Using \eqref{StrongConvexityEntropyEq}, this implies that $ \cS $ is strongly $ 0$-convex.
It remains to show the existence of the measure $ \omega $.
Let $ \pi_0 \in \Opt{\s,\mu_0} $ and $ \pi_1 \in \Opt{\s,\mu_1} $.
Using the disintegration theorem, we obtain the existence of Borel measurable families $ (\pi_0^{x,m})_{x\in\Td,m\in\R},(\pi_1^{x,m})_{x\in\Td,m\in\R} \subset \cM_1(\R) $ such that 
\be{StrongConvexityEntropyEq4}
\pi_0 = \pi_0^{x,m} \, d\s^x(m)\, dx
 \quad \text{and} \quad
\pi_1 = \pi_1^{x,m} \, d\s^x(m)\, dx.
\ee
Using the measurable selection lemma (\cite[5.22]{vil}), we know that there exists $ (\omega^{x,m})_{x\in\Td,m\in\R} \subset \cM_1(\R^2) $ such that 
\be{StrongConvexityEntropyEq5}
\omega^{x,m} \in \OptW{\pi_0^{x,m},\pi_1^{x,m}} \quad \text{and} \quad
(x,m) \to \omega^{x,m} \text{ is measurable.}
\ee
Define $ \omega := \omega^{x,m} \, d\s^x(m)\, dx \in \cM_1(\Td\times\R^3) $.
It is easy to see that $ \omega  $ fulfils \eqref{StrongConvexityEntropyEq1}. Indeed, for all Borel-measurable $ M,A \subset \R  $
\begin{align}\begin{split}
 \omega^x(M \times A \times \R) = \int_M   \omega^{x,m}(A \times \R) \, d\s^x(m) =  \int_M  \pi_0^{x,m}(A) \, d\s^x(m) = \pi_0^x(M\times A ). 
\end{split}
\end{align}
Therefore, $ (\pb_{\R^3}^{1,2}) _\# \omega^x \in \OptW{\s^x,\mu_0^x} $, since we chose $ \pi_0 \in \Opt{\s,\mu_0} $.
Analogously, one can show that 
$ (\pb_{\R^3}^{1,3}) _\# \omega^x \in \OptW{\s^x,\mu_1^x} $.
 \epr

  \paragraph{Interaction energy.}
  Define the \emph{interaction energy} $\cW:\PL \ra (-\infty, \infty]$  by
  \begin{align}\label{InteractionFunctional}
  \cW(\mu) := 
\frac{1}{2} \int_{\TdR} \int_{\TdR}  W(x,\bar{x},\t,\bar{\t}) \,  d\mu(x,\t) d\mu(\bar{x},\bar{\t}),
  \end{align}
  where 
  $ W \in C^{0,0,1,1} (\Td \times \Td \times \R \times \R)  $ satisfies the following assumptions. 
  \bass{InteractionAssumptions}
  \begin{enumerate}[(1)]
  	\item $ W(x,\bar{x},\t,\bar{\t}) \geq - \alpha ( |(\t,\bar{\t})|^2 + 1) $ for some $ \alpha > 0 $. 
  	\item There exists $ \bar{\lambda} \in \R $ such that for all $ (x,\bar{x}) \in  \Td \times \Td  $, $ (\t,\bar{\t}) \to W(\bar{x},x,\t,\bar{\t}) $ is $ \bar{\lambda}$-convex, i.e.
  	for all $ (\t_1,\bar{\t}_1), (\t_2,\bar{\t}_2) \in \R^2 $
  	\begin{align}\label{WConvex}
  	W(x,\bar{x},(1-t) \t_1 + t \t_2,(1-t) \bar{\t}_1 + t \bar{\t}_2)  
  	\leq &(1-t) W(x,\bar{x},\t_1 ,\bar{\t}_1) + t \, W(x,\bar{x},\t_2 ,\bar{\t}_2)  \nonumber \\
  	&   - \tfrac{\bar{\lambda}}{2}t ( 1-t) \, \big| (\t_1,\bar{\t}_1)- (\t_2,\bar{\t}_2)\big|^2.
  	\end{align}  	
  \end{enumerate}
\eass

  \bl{InteractionProp}
  Suppose that Assumption \ref{InteractionAssumptions} is satisfied. Then $ \cW $ is well-defined, coercive, strongly $ \bar{\lambda}$-convex 
and 
  	$ \WL $-l.s.c.
  \el
  \bpr Assumption \ref{InteractionAssumptions} (1) implies that $ \cW(\mu)  \geq \, - \alpha \int |\t|^2 d\mu -\alpha$ for all $ \mu \in \PL $. This shows that $ \cW $ is well-defined and coercive.

  Let $ (\mu_n)_n $ and $  \mu\in\PL $ be such that $ \lim_{n\ra\infty}  \WL(\mu_n,\mu)=0$.
   From  \cite[Theorem 2.8]{billi} and Lemma \ref{ConvCharacLem}, we obtain that $ \mu_n \otimes  \mu_n \hra \mu  \otimes \mu $ and $ \lim_{n\ra\infty}  \int |\t|^2 d(\mu_n \otimes  \mu_n) = \int |\t|^2 d(\mu \otimes  \mu)  $.
   Therefore, by Assumption \ref{InteractionAssumptions} (1), it is straightforward to see that 
 $ W^- $ is uniformly integrable with respect to $ (\mu_n)_n $. 
 Hence, \cite[5.1.7]{ambgigsav} implies that	$ 
 \liminf_{n\ra \infty } 	\cW(\mu_n) \geq \cW(\mu).
 $

  It remains to show the strong $ \bar{\lambda}$-convexity.
  Let $ \s,\mu_0,\mu_1 \in D(\cW) $ and $ \Phi $ be as in \eqref{StrongConvexDefEq} for the functional $ \cW $.
 Let $ \omega \in \cM_1^\L(\Td\times\R^3) $ and $ (\gamma_t)_{t\in[0,1]} $ be as  in the proof of Proposition \ref{EntropyProp} (iv).  
 Since $ \left(\pb_{\R^3}^{1}, (1-t )\, \pb_{\R^3}^{2}  + t\, \pb_{\R^3}^{3}  \right) _\#\omega^x $ is a coupling of $ \s^x $ and $ \gamma_t^x $ for almost every $ x \in \Td $ and using \eqref{StrongConvexityEntropyEq1}, we obtain that for all $ t \in [0,1] $
 \begin{align}\label{StrongConvexityInterMetric}
\begin{split}
 \int_{  \Td} W_2(\s^x, &\gamma_t^x)^2 \, dx 
 \leq\int_{  \Td} \int_{  \R^3	} |(1-t)\t_2+t\,\t_3 -\t_1|^2 d\omega^x(\t_1,\t_2,\t_3) \, dx   \\
 &= (1-t)\WL(\s,\mu_0)^2 +  t\, \WL(\s,\mu_1)^2   - t(1-t ) \int_{  \Td} \int_{  \R^3	} |\t_2-\t_3|^2 d\omega^x \, dx .  
\end{split}
 \end{align}
  Moreover, Assumption \ref{InteractionAssumptions} (2) implies that 
   \begin{align}\label{StrongConvexityInterFunctional}
   \cW(\gamma_t)&= \frac{1}{2} \int_{(\Td\times \R^3)^2}   W(x,\bar{x},(1-t) \t_2 + t \t_3,(1-t) \bar{\t}_2 + t \bar{\t}_3)  d\omega(x,\t_1,\t_2,\t_3) d\omega(\bar{x},\bar{\t}_1,\bar{\t}_2,\bar{\t}_3) \nonumber
   \\
  & \leq  (1-t)\cW(\mu_0)  +  t\, \cW(\mu_1) -  \frac{\bar{\lambda}}{2} t(1-t ) \int_{  \Td \times \R^3	} |\t_2-\t_3|^2 \,  d\omega (x,\t_1,\t_2,\t_3) .
  \end{align}
  \eqref{StrongConvexityInterMetric} and \eqref{StrongConvexityInterFunctional} yield that
  for all $ \tau \in (0, \tfrac{1}{\bar{\lambda}^-}) $ and for all $ t \in [0,1] $
   \begin{align}\label{StrongConvexityInterEq1}
  \Phi(\tau,\s;\gamma_t) 
  &\leq   (1-t)\Phi(\tau,\s;\gamma_0)  +  t\, \Phi(\tau,\s;\gamma_1) -  \left(\frac{1}{2\tau } + \frac{\bar{\lambda}}{2}\right) t(1-t ) \int_{  \Td \times \R^3	} |\t_2-\t_3|^2 \,  d\omega \nonumber		\\
  &\leq   (1-t)\Phi(\tau,\s;\gamma_0)  +  t\, \Phi(\tau,\s;\gamma_1) -  \left(\frac{1}{2\tau } + \frac{\bar{\lambda}}{2}\right) t(1-t ) \WL(\mu_0,\mu_1),
  \end{align}
  which is also a consequence of \eqref{StrongConvexityEntropyEq1}. This concludes the proof.
\epr
  
    \paragraph{Potential energy.}
  Define the \emph{potential energy} $\cV:\PL \ra (-\infty, \infty]$  by
  \begin{align}\label{PotentialFunctional}
  \cV(\mu) := 
\int_{\TdR}  V d\mu , 
  \end{align}
  where 
    $ V \in C^{0,1} (\TdR)  $ satisfies the following assumptions. 
  \bass{PotentialAssumptions}
  \begin{enumerate}[(1)]
  	\item $ V(x,\t) \geq - \alpha ( |\t|^2 + 1) $ for some $ \alpha > 0 $. 
  	\item There exists $ \hat{\lambda} \in \R $ such that for all $ x \in  \Td  $, $ \t \to V(x,\t) $ is $ \hat{\lambda}$-convex.
  \end{enumerate}
  \eass
  It turns out that under these assumptions, the potential energy is just the special case of the interaction energy, when $ W(x,\bar{x},\t,\bar{\t}) = V(x,\t) +V(\bar{x},\bar{\t}) $.
  Therefore, all the results for the interaction energy carry over to the potential energy and we have nothing to prove here. 
  
  \bl{PotentialProp}
  Suppose that Assumption \ref{PotentialAssumptions} is satisfied. Then, $ \cV $ is well-defined, coercive, strongly $ \hat{\lambda}$-convex and 
  $ \WL $-l.s.c. 
  \el
  
\paragraph{The McKean-Vlasov-functional $ \boldsymbol{\cF} $.}
  From now on, we specify the functionals $ \cW $ and $ \cV $ as follows. 
  \bass{Fass}
   \begin{enumerate}[(1)]
  	\item $W(x,\bar{x},\t,\bar{\t}) = -J(x-\bar{x}) \, \t\,\bar{\t}, $ where $ J:\Td \ra \R $ is continuous and symmetric.
  	  It is easy to see that Assumption \ref{InteractionAssumptions} is satisfied.
  	Indeed, as   an immediate consequence  of Young's inequality, Assumption \ref{InteractionAssumptions} (2) is satisfied for $ \bar{\lambda}:= -\|J\|_\infty $.  	
  	\item  $ V(x,\t) = \Psi(\t),  $ where $ \Psi :\R\ra\R$ is assumed to be a polynomial of even degree  such that 
  	 Assumption \ref{PotentialAssumptions} (2) is satisfied for some $ \hat{\lambda}\in \R $, and 
  	\be{LBPsi}
  	\Psi (\t ) \geq C_\Psi \t^{2\ell} +  C_\Psi' \t^2 -  C_\Psi''  \quad \text{for all $ \t \in \R $},
  	\ee
  	for some $ \ell \in \N $,   $ C_\Psi, C_\Psi''  \geq 0  $ and $ C_\Psi' > \|J\|_\infty $. 
%
  \end{enumerate}
\eass

For example, if $ \Psi $ is  a polynomial of degree $ 2\ell $, then $ \Psi $  satisfies Assumption \ref{Fass}, where, if $ \ell=1 $, we assume
	that the coefficient of degree $ 2 $ is strictly greater than $ \|J\|_\infty $. 
%
%
%
%

Assumption \ref{Fass} implies that $ \cF $ has the form
  \begin{align}\label{McKeanVlasovFunctional}
  \cF(\mu) = 
  \int_{\TdR}  \log(\rho) d\mu +  \int_{\TdR}  \Psi \, d\mu  - \frac{1}{2} \int_{\TdR} \int_{\TdR}  J(x-\bar{x}) \t\bar{\t} d\mu(x,\t) d\mu(\bar{x},\bar{\t})
  \end{align}
  if $ \mu $ has a density $ \rho $ with respect to $ \Leb{\TdR} $ and $ \cF(\mu) = \infty $ otherwise.  
  Note that $ \cF $ is proper (e.g.\ $ \cF( \exp(-\t^2/2) (2\pi)^{-1/2} d\t dx) < \infty $).    
Furthermore, the definition of $ \cF $ can be naturally extended to $ \cM_1(\TdR) $. 
We observe the following lower bound on $ \cF $.
  \bl{LBF}
We have,  for some constant $ C'' >0 $,
 \be{LBFEq}
 \cF(\mu)  \geq  \int_{  \TdR	} \Big(C_\Psi\, |\t|^{2\ell} + ( C_\Psi' - \|J\|_\infty )\, |\t|^2\Big) d\mu - C''  \qquad \text{for all $ \mu \in \cM_1(\TdR) $}. 
\ee
In particular, there exists $ \mu \in \PL $ such that $ \inf_{\sigma \in \cM_1(\TdR)} \cF(\s) = \cF(\mu) $ and $ D(\cF) \subset \{\mu \in \cM_1(\TdR) \, | \,  \int |\t|^2\, d\mu < \infty \} $. 
\el
 \bpr
 Let $ \mu \in \cM_1(\TdR) $ and assume that $ \mu $ has a density $ \rho $, since otherwise the claim is trivial. 
Notice that we can rewrite $ \cF  $ as  
\begin{align}\label{LBFEqq}
\cF(\mu) = \cH\left ( \mu\, \Big| \, e^{-\frac{1}{2} \Psi(\t)} d\t dx \right) + \frac{1}{2} \int_{  (\TdR)^2	} \left(\frac{1}{2} \big(\Psi(\t) +  \Psi(\bar{\t})\big) - J(x-\bar{x}) \t \bar{\t}			\right) d\mu d\mu.
\end{align}
Then, since 
\begin{align}\label{LBFEqq1} 
 \cH\Big( \mu\, &\Big| \, e^{-\frac{1}{2} \Psi(\t)} d\t dx \Big)  
\geq - \log  \int_{\ToR} e^{-\frac{1}{2} \Psi(\t)} d\t dx ,
\end{align}
and by  Young's inequality and \eqref{LBPsi},
\begin{align}\label{Flsc1}
\frac{1}{2} \big(\Psi(\t) +  \Psi(\bar{\t})\big) - J(x-\bar{x}) \t \bar{\t}	&\geq  \frac{1}{2} C_\Psi   (\t^{2\ell}+  \bar{\t}^{2\ell}) +  \frac{1}{2} ( C_\Psi' - \|J\|_\infty ) (\t^2+  \bar{\t}^2) - C_\Psi'' ,
\end{align}
we infer \eqref{LBFEq}.
 
For the second claim, note that \eqref{LBFEq} implies the weak compactness of the level sets of $ \cF $  and that Theorem \ref{FProp} below shows the weak lower semi-continuity of  $ \cF $.
Therefore, the direct method of the calculus of variation is applicable and we infer the existence of a minimizer.
 \epr
 

As a consequence of the observations on $ \cS, \cV $ and $ \cW $, we obtain the following result for $ \cF $.
\bt{FProp}
$ \cF $ is well-defined, proper, coercive, $  (\bar{\lambda}+\hat{\lambda} )$-convex, strongly $ \lambda$-convex for some $ \lambda  \in \R $ and
lower semi-continuous with respect to weak convergence.
In particular, $ \cF  $ is $ \WL$-l.s.c.
Therefore, Theorem \ref{ExistenceThm} is applicable for $ \cF$.
\et
\bpr
It remains to show that $ \cF $ is weakly lower semi-continuous. 
Let 
$ (\mu_n)_{n\in\N} \subset  \cM_1(\TdR) $  and $ \mu\in  \cM_1(\TdR) $ be such that  $ \mu_n \hra \mu $. Without restriction suppose that $ \cF(\mu^n) < \infty $ for all $ n \in \N $ and $ \liminf_{n\ra\infty}\cF(\mu^n) < \infty $.
We show the lower semi-continuity for both summands on the right-hand side of \eqref{LBFEqq} separately.
In the proof of Lemma \ref{EntropyProp} we have already seen that the functional $ \cH\big ( \cdot \, \big| \, \frac{1}{\alpha} e^{-\frac{1}{2} \Psi(\t)} d\t dx \big) $ is weakly lower semi-continuous, where $ \alpha = \int  e^{-\frac{1}{2} \Psi(\t)} d\t dx $. 
Therefore, 
\begin{align}\label{Flsc2} 
\begin{split}
 \liminf_{n\ra \infty }  \cH\Big( \mu^n\, &\Big| \, e^{-\frac{1}{2} \Psi(\t)} d\t dx \Big)  
 =  \liminf_{n\ra \infty }  \cH\Big ( \mu^n \, \Big| \, \tfrac{1}{\alpha} e^{-\frac{1}{2} \Psi(\t)} d\t dx \Big)  - \log(\alpha) \\
 &\geq \cH\Big ( \mu \, \Big| \, \tfrac{1}{\alpha} e^{-\frac{1}{2} \Psi(\t)} d\t dx \Big)  - \log(\alpha) = \cH\Big ( \mu \, \Big| \, e^{-\frac{1}{2} \Psi(\t)} d\t dx \Big)  .
\end{split}
 \end{align}
Moreover,   the integrand of the second summand  in \eqref{LBFEqq} is lower semi-continuous and  bounded from below due to \eqref{Flsc1}. 
Therefore,   \cite[5.1.7]{ambgigsav} yields
\begin{align}\label{Flsc3} 
\begin{split}
\liminf_{n\ra \infty }   \int_{  (\TdR)^2	} &\Big(\frac{1}{2} \big(\Psi(\t) +  \Psi(\bar{\t})\big) - J(x-\bar{x}) \t \bar{\t}			\Big) d(\mu^n\otimes \mu^n) \\
&\geq \int_{  (\TdR)^2	} \left(\frac{1}{2} \big(\Psi(\t) +  \Psi(\bar{\t})\big) - J(x-\bar{x}) \t \bar{\t}			\right) d(\mu\otimes \mu), 
\end{split}
\end{align}
which concludes the proof.
\epr

  \paragraph{Directional derivative.}\label{Dirrr}
In order to find a characterisation of the (strong) subdifferential of $ \cF $, it will be useful to study the infinitesimal behaviour of $ \cF $ along 
curves that are pushed along smooth functions. 
This is the content of the following proposition. 
  \bp{DirDerProp}
Let $ \mu \in D(\cF) $ and $ \beta \in C_c^{2}(\TdR\,;\,\R) $.
For all $ t \in \R $ define 
\be{DirDerDef}
\mu_{t,\beta} = (\pb^1, \pb^2 + t \, \beta)_\# \mu.
\ee
Then 
\be{DirDer}
\restr{ \frac{d}{dt}}{t=0} \cF(\mu_{t,\beta}) = \int_{  \TdR} \left(\beta(x,\t) \left[ \Psi'(\t)	 - 	 \int_{  \TdR}J(x-\bar{x})\bar{\t}      \, d\mu(\bar{x},\bar{\t})\right]	-\partial_\t \beta(x,\t)    \right) d\mu(x,\t).
\ee
\ep
  \bpr
We compute the derivative for each summand separately. We begin with $ \cS $.
    Again, the proof is similar to the Wasserstein case. Indeed, if we consider the function $ \hat{\beta} = (0,\beta) \in C_c^2(\TdR\,;\,\TdR) $, then 
$ 	\mu_{t,\beta} = (\Id_{\TdR} + t \, \hat{\beta})_\# \mu $ for all $ t \in \R $.
Then,  \cite[(38)]{jko} implies that 
\be{DirDerEntropy}
\restr{ \frac{d}{dt}}{t=0} \cS(\mu_{t,\beta}) = - \int_{  \TdR} \mathrm{div} \hat{\beta} \, d\mu = - \int_{  \TdR} \partial_\t \beta \, d\mu .
\ee
To compute the directional derivative of $ \cW $, observe that
\begin{align}\label{DirDerInteraction}
\restr{ \frac{d}{dt}}{t=0} \cW(\mu_{t,\beta}) 
&= 
- \frac{1}{2} \restr{ \frac{d}{dt}}{t=0} \int_{\TdR} \int_{\TdR}    J(x-\bar{x})   (\t + t \beta(x,\t))(\bar{\t} + t \beta(\bar{x},\bar{\t}) )d\mu(x,\t) d\mu(\bar{x},\bar{\t}) \nonumber \\ 
&= 
- \frac{1}{2}   \int_{\TdR} \int_{\TdR} J(x-\bar{x}) \restr{ \frac{d}{dt}}{t=0} (\t + t \beta(x,\t))(\bar{\t} + t \beta(\bar{x},\bar{\t}) )  d\mu(x,\t) d\mu(\bar{x},\bar{\t}) \nonumber \\ 
&= 
-  \int_{\TdR} \int_{\TdR} J(x-\bar{x}) \,  \bar{\t} \, \beta(x,\t) \, d\mu(x,\t) \  d\mu(\bar{x},\bar{\t}),
\end{align}
where we have used the symmetry of the integrand. To exchange differentiation and integration, we have used the Leibniz-integral-rule, which is applicable, since all functions are continuous and $ \beta  $ has compact support. 
In the same way, one computes that 
    \be{DirDerPotential}
    \restr{ \frac{d}{dt}}{t=0} \cV(\mu_{t,\beta}) = \int_{  \TdR}  \beta \, \Psi' \ d\mu.
    \ee
\epr

  \paragraph{Subdifferential of $ \boldsymbol{\cF} $.}
  Note that for all $ \mu \in \PL $, $	(x,\t) \to \int_{  \TdR}J(x-\bar{x})\bar{\t}      \, d\mu(\bar{x},\bar{\t}) \in \L^2(\mu) $. This observation is important in order to compute an element of the subdifferential of $ \cF $ in the following proposition.

  \bl{SubdiffF}
  Let $ \mu \in D(\cF) $. Therefore, $ \mu $ has a density $ \rho $ with respect to $ \Leb{\TdR} $.
  Suppose that $ \partial_\t \rho $ exists weakly in  $  \mathrm{L}^{1}_{\mathrm{loc}}(\TdR) $ and $ \frac{\partial_\t \rho}{\rho} + \Psi' \in \L^2(\mu) $.
  Then 
  \be{SubDiffFEq}
\left(   (x,\t) \to \frac{\partial_\t \rho(x,\t)}{\rho(x,\t)} + \Psi'(\t)	 - 	 \int_{  \TdR}J(x-\bar{x})\bar{\t}      \, d\mu(\bar{x},\bar{\t})   \right) \  \in  \ \partial  \cF (\mu)
  \ee
  \el
  \bpr
We first show that 
 $ 
  \left(   (x,\t) \to -	 \int_{  \TdR}J(x-\bar{x})\bar{\t}      \, d\mu(\bar{x},\bar{\t})   \right)  \in  \partial  \cW (\mu).
 $ 
Note that  for all $ (\t_1,\bar{\t}_1),(\t_2,\bar{\t}_2) \in \R^2 $
    \begin{align}\label{SubDiffWEq1}
-	J(x-\bar{x})& \big(\t_1 \bar{\t_1} - \t_2 \bar{\t}_2 \big)
= - J(x-\bar{x})\Big( \bar{\t}_2 (\t_1-\t_2)  + \t_2 (\bar{\t_1}-\bar{\t}_2) + (\t_1- \t_2)(\bar{\t}_1- \bar{\t}_2)\Big) \nonumber \\
&\geq - J(x-\bar{x})\Big( \bar{\t}_2 (\t_1-\t_2)  + \t_2 (\bar{\t_1}-\bar{\t}_2) \Big) 
+\frac{\bar{\lambda}}{2} \ \big| (\t_1,\bar{\t}_1)- (\t_2,\bar{\t}_2)\big|^2 
\end{align}
This yields for all $ \nu \in D(\cF) \subset D(\cW) $
\begin{align}\label{SubDiffInteraction}
\cW(\nu) -\cW(\mu)
= 
&\frac{1}{2} \int_{(\TdR)^2}     -J(x-\bar{x})   \Big( \T_{\mu}^{\nu}(x,\t) \T_{\mu}^{\nu}(\bar{x},\bar{\t})  -\t \bar{\t}\Big)   d(\mu\otimes\mu) \nonumber \\ 
\geq 
&\frac{1}{2} \int_{(\TdR)^2}     -J(x-\bar{x})    \bar{\t} ( \T_{\mu}^{\nu}(x,\t)  -\t )   d(\mu\otimes\mu)\nonumber  \\ 
&+ 
\frac{1}{2} \int_{(\TdR)^2}     -J(x-\bar{x})     \t ( \T_{\mu}^{\nu}(\bar{x},\bar{\t})  -\bar{\t} )   d(\mu\otimes\mu)  \\ 
&+\frac{\bar{\lambda}}{4} \int_{(\TdR)^2}      \Big| (\T_{\mu}^{\nu}(x,\t), \T_{\mu}^{\nu}(\bar{x},\bar{\t}))  -(\t ,\bar{\t} )\Big|^2
d(\mu\otimes\mu) \nonumber \\ 
= 
& \int_{\TdR}    \left( - \int_{\TdR} J(x-\bar{x})    \bar{\t} d \mu(\bar{x},\bar{\t}) \right)  ( \T_{\mu}^{\nu}(x,\t)  -\t ) \,   d\mu(x,\t) +\frac{\bar{\lambda}}{2} \WL(\mu,\nu)^2. \nonumber
\end{align}

It remains to show that $ \partial_\t \rho/\rho + \Psi' \in  \partial  (\cS+\cV) $. 
Notice that $ \t \mapsto \Psi(\t)- \frac{1}{2} \hat{\lambda}|\t|^2$ is convex. 
Set $ \tilde{V}(\t) := \Psi(\t)- \frac{1}{2} \hat{\lambda}|\t|^2 + \b $, where  $ \b\in\R $ is such that $ \exp(-\tilde V(\t)) d\t $ is a probability measure.
Define $ \tilde{\cV}(\mu) := \int \tilde V d\mu$.
Then, similarly as in \eqref{EntropyRelEntropy} and \eqref{EntropyFunctional2}, we have that 
\begin{align}\label{EntropyRelEntropy2}
\cS(\mu)  +\tilde{\cV} (\mu) = \cH\Big(\mu \, \Big| \,  e^{-\tilde V(\t)} d\t  dx\Big)  =
\int_{\TdR}  \cH\Big(\mu^x\, \Big| \,  e^{-\tilde V(\t)} d\t \Big) dx.
\end{align}
This fact allows us to use the results from the Wasserstein case. 
By taking a compact exhaustion of $ \R $, one can see immediately that there exists a nullset $ \cN \subset \Td $ such that 
for all $ x \in \Td\setminus \cN $ 
\be{SubDiffEVEq}
\rho(x,\cdot) \in \mathrm{W}^{1,1}_{\mathrm{loc}}(\R) , \qquad  \frac{\partial_\t \rho(x,\cdot)}{\rho(x,\cdot)} + \Psi' \in \L^2(\mu^x) \qquad \text{and} \qquad \int_{  \R	}|\t|^2 d\mu^x < \infty.
\ee
Moreover, if we set $ \s^x(\t) =\rho(x,\t) \exp(\tilde{V}(\t)) $, then \eqref{SubDiffEVEq} implies that 
\be{SubDiffEVEq2}
\s^x \in \mathrm{W}^{1,1}_{\mathrm{loc}}(\R) \qquad \text{and} \qquad  \frac{\partial_\t \s^x}{\s^x}  \in \L^2(\mu^x) \qquad \text{ for almost every } x.
\ee
Therefore, \cite[10.4.9]{ambgigsav} is applicable and we obtain that for all $ \nu \in D(\cF) $
\begin{align}\label{SubDiffEVEq3}
  \cH\Big(\nu^x\, \Big| \,  e^{-\tilde V(\t)} d\t \Big) 
  -\cH\Big(\mu^x\, \Big| \,  e^{-\tilde V(\t)} d\t \Big) 
  \geq \int_{\R}     \frac{\partial_\t \s^x}{\s^x}  ( \T_{\mu^x}^{\nu^x}  -\Id_{\R} ) \,   d\mu^x
  \quad \text{ for a.e.\ } x.
\end{align}
  Using \eqref{EntropyRelEntropy2}, this implies that 
 \begin{align}\label{SubDiffEVEq4}
 (\cS +\tilde{\cV})(\nu) - (\cS +\tilde{\cV})(\mu) 
 \geq \int_{\TdR}     \Big(\frac{\partial_\t \rho}{\rho} + \Psi'  -  \hat{\lambda} \pb^2 \Big) ( \T_{\mu}^{\nu}  -\pb^2 ) \,   d\mu.
 \end{align} 
  Since $ \WL(\mu,\nu) = \|\T_{\mu}^{\nu}  -\pb^2 \|_{\L^2(\mu)}  $, we infer
   \begin{align}\label{SubDiffEVEq5}
  (\cS +\cV)(\nu) - (\cS +\cV)(\mu) 
  \geq \int_{\TdR}     \Big(\frac{\partial_\t \rho}{\rho} + \Psi'   \Big) ( \T_{\mu}^{\nu}  -\pb^2 ) \,   d\mu + \frac{\hat{\lambda}}{2}  \WL(\mu,\nu),
  \end{align} 
  which concludes the proof.
  \epr
  
  \
  
   \bp{SubdiffCharacF}
  Let $ \mu = \rho \, \Leb{\TdR} \in D(\cF) $.
  Then the following statements are equivalent. 
  \begin{enumerate}[(i)]
  	\item $ |\partial \cF |(\mu) < \infty $,
  	\item $ \partial_\t \rho $ exists weakly in  $  \mathrm{L}^{1}_{\mathrm{loc}}(\TdR) $ and 
  	there exists 
  	$ w \in \L^2(\mu) $ such that 
  	$ \partial_\t \rho (x,\t) = \rho(x,\t) (w(x,\t) -  \Psi'(\t) + 	 \int_{  \TdR} J(x-\bar{x})\bar{\t}      \, d\mu(\bar{x},\bar{\t}))   $.
  \end{enumerate}
  Moreover, in this case, $ w \in \TanL{\mu}\,  \cap \, \partial \cF (\mu)  $, $  |\partial \cF |(\mu) = \|w \|_{\L^2(\mu)}  $ and $ w $ is the $ \mu$-a.e.\ unique strong subdifferential at $ \mu $.
  \ep
  \bpr
  Again, the proof is very similar to the Wasserstein case (cf. \cite[4.3]{erbar}).
  However, here we include the details, since the statement and the proof will become very crucial for the remainder of this paper. 
  
$ (ii) \Rightarrow (i) $.
  Lemma \ref{SubdiffF} shows that under the conditions of \emph{(ii)}, $ w \in \partial \cF (\mu)  $. Hence, $ |\partial \cF |(\mu) \leq  \|w \|_{\L^2(\mu)} < \infty $ , which is an immediate consequence of the definition of the metric slope (cf.\ \cite[10.3.10]{ambgigsav}).
  
 $ (i) \Rightarrow (ii) $.
Define a linear operator $ L : C_c^{\infty}(\TdR) \ra \R  $ by 
\be{SubdiffCharacFEq1}
L(\beta) := \int_{  \TdR} \left(\beta(x,\t)  \left[ \Psi'(\t)  - 	 \int_{  \TdR}J(x-\bar{x})\bar{\t}      \, d\mu(\bar{x},\bar{\t}) \right]	-\partial_\t \beta(x,\t)     \right) d\mu(x,\t) .
\ee
Let $ \beta \in  C_c^{\infty}(\TdR) $ and $ (\mu_{t,\beta})_t $ be as in Proposition \ref{DirDerProp}.  
Using the   representation  \eqref{MetricCouplRel} and that $ (\pb^1,\pb^2+t\beta,\pb^2)_\# \mu \in \Coupl{\mu_{t,\beta},\mu} $, it is easy to see that $  \WL(\mu_{t,\beta},\mu)  \leq |t| \cdot \|\beta \|_{\L^2(\mu)}  $.
Then, as in \cite[p.\ 13, l.\ 12]{erbar}, via
Proposition \ref{DirDerProp}, we observe that if $ L(\beta) >  0 $,
\begin{align}\label{{SubdiffCharacFEq2}}
L(\beta) &=  \lim_{t\downarrow  0} \frac{ \left( \cF(\mu) - \cF(\mu_{-t,\beta})\right) ^+  }{t} 
 \leq   \limsup_{t\downarrow 0} \frac{\left( \cF(\mu) - \cF(\mu_{t,-\beta})\right) ^+}{  \WL(\mu_{t,-\beta},\mu)}  \|\beta \|_{\L^2(\mu)} 
 \leq |\partial \cF |(\mu) \ \|\beta \|_{\L^2(\mu)}, 
\end{align}
and if $ L(\beta) <  0 $,
\begin{align}\label{{SubdiffCharacFEq22}}
L(\beta)  &=  \lim_{t\downarrow  0} \frac{ \left( \cF(\mu) - \cF(\mu_{t,\beta}) \right) ^+ }{-t}  \geq  - \liminf_{t\downarrow 0} \frac{\left( \cF(\mu) - \cF(\mu_{t,\beta})\right) ^+}{ \WL(\mu_{t,\beta},\mu)}  \|\beta \|_{\L^2(\mu)} 
\geq -|\partial \cF |(\mu) \ \|\beta \|_{\L^2(\mu)}. 
\end{align}
Thus, $ |L(\beta)|  \leq |\partial \cF |(\mu) \ \|\beta \|_{\L^2(\mu)} $. 
Extending $ L $ to the $ \L^2( \mu) $-closure of $ 
C_c^\infty ( \TdR) $, 
the Riesz representation theorem yields the existence of a unique $ w \in \L^2(\mu) $ such that  
\begin{itemize}
	\item $ \int w \beta \rho \, d\t dx =  \int \left(\beta \left[ \Psi' - 	 \int J(\cdot-\bar{x})\bar{\t}      \, d\mu\right]	-\partial_\t \beta    \right) \rho\, d\t dx$
	for all $ \beta \in  C_c^{\infty}(\TdR)$, and 
	\item $  |\partial \cF |(\mu) \geq  \|w \|_{\L^2(\mu)}  $.
\end{itemize}
This shows that the weak derivative $ \partial_\t \rho  $ exists and equals $ \rho (w -  \Psi' + 	 \int_{  \TdR} J(x-\bar{x})\bar{\t}      \, d\mu) $, which clearly belongs to $  \mathrm{L}^{1}_{\mathrm{loc}}(\TdR) $. We infer \emph{(ii)}.

It remains to show the other claims. Let $ \pb^\mathrm{Tan} $ denote the orthogonal projection onto $ \TanL{\mu} $.
Then $ \pb^\mathrm{Tan}(w) \in \partial \cF (\mu)  $, since $ w \in \partial \cF (\mu)  $. Indeed, this follows immediately from the definition of the subdifferential and Corollary \ref{LOptTangLem}.
Hence, by Lemma \ref{TangPropLem}
\be{SubdiffCharacFEq3}
|\partial \cF |(\mu) \leq  \|\pb^\mathrm{Tan}(w)  \|_{\L^2(\mu)}  \leq   \|\pb^\mathrm{Tan}(w) +w - \pb^\mathrm{Tan}(w)   \|_{\L^2(\mu)} = \|w \|_{\L^2(\mu)}  \leq |\partial \cF |(\mu),
 \ee
which, again by Lemma \ref{TangPropLem}, shows   that $ w \in \TanL{\mu}$ and $  |\partial \cF |(\mu) =  \|w \|_{\L^2(\mu)}  $. 

Finally, let $ z  $ be another strong subdifferential of $ \cF $ at $ \mu $.
Then, for all $ \beta \in  C_c^{\infty}(\TdR)$
\be{SubdiffCharacFEq4}
\int w  \beta  \, d\mu = L(\beta ) =   \restr{ \frac{d}{dt}}{t=0} \cF(\mu_{t,\beta}) 
= \lim_{t\downarrow 0} \frac{\cF(\mu_{t,\beta}) - \cF(\mu)}{t} 
\geq \int z \beta  \, d\mu,
\ee
since 
$ z  $ is a strong subdifferential.
Considering  $ \lim_{t\uparrow 0} $, we obtain the other inequality.
Therefore, $ \int w  \beta  \, d\mu =  \int z  \beta  \, d\mu, $ for all $ \beta \in  C_c^{\infty}(\TdR)$, which implies that 
$ z=w $ $ \mu$-a.e.
  \epr 
  \begin{samepage}
  	\vspace{-0.5cm}
  \bc{FSubDiffCor} 
  Let $ \mu  \in D(\cF) $.
  Then  
  \be{FSubCor}
  |\partial \cF |(\mu) = \sup_{\beta \in C_c^\infty ( \TdR),\, \| \beta \|_{\L^2( \mu)} >0} \frac{\left|\int_{\TdR} \left(\beta \left[ \Psi' - 	 \int_{\TdR} J(\cdot-\bar{x})\bar{\t}      \, d\mu(\bar{x},\bar{\t}) \right]	-\partial_\t \beta    \right) d\mu\right|}{\| \beta \|_{\L^2( \mu)}}.
  \ee
  Moreover, if   $ (\mu_n)_{n\in\N}  $  is such that $ \sup_n  \int |\t|^2 d\mu_n<\infty$ and $ \mu_n \hra \mu\in\PL $, then
  	 \be{LSCFisher}
  	 \liminf_{n\ra \infty } 	 |\partial \cF |(\mu_n) \geq  |\partial \cF |(\mu).
  	 \ee 
  \ec
  \bpr
If  $ |\partial \cF |(\mu) < \infty$, then \eqref{FSubCor} follows from  the proof of Proposition \ref{SubdiffCharacF}, since  the right-hand side of \eqref{FSubCor} equals $ \|w \|_{\L^2(\mu)} $. Here we have  used that the extension of the operator $ L $ from \eqref{SubdiffCharacFEq1} has the same operator norm as $ L $.
And if the right-hand side of \eqref{FSubCor} is finite, then  $ L $ is bounded. 
Therefore, repeating the above arguments,  we infer part $ (ii) $ of Proposition \ref{SubdiffCharacF}, which leads to $ |\partial \cF |(\mu) < \infty$ and finally to \eqref{FSubCor}.

The proof of \eqref{LSCFisher} is a straightforward consequence of \eqref{FSubCor}, the fact that $ \beta \in C_c^\infty ( \TdR) $, and \cite[5.1.7]{ambgigsav}.
  \epr   
\end{samepage}

  \paragraph{Variational characterisation of the gradient flow for $ \boldsymbol{\cF} $. } 
  The following characterisation of gradient flows for $ \cF $ is a key fact in order to establish the results in the forthcoming chapters.
  \bt{Lyapunov}
 Let $ T \in (0, \infty) $. Define $ \cJ: C([0,T] \, ; \, \cM_1(\ToR) ) \ra [0,\infty] $ by
      \be{LyapunovEq}
\cJ[ (\nu_t)_t ] :=
	\cF(\nu_T) - 		\cF(\nu_0) + \frac 12 \int_0^T\big(|\partial \cF|^2(\nu_t)  + |\nu'|^2(t)\big)\,dt ,	
\ee
if  $ (\nu_t)_t \in \ACLT $  and 
$ \cJ[ (\nu_t)_t ] = \infty $ 	else.
 Let $ \mu_0 \in D(\cF) $. For any curve $ \muc \in \ACLT $ such that   
  $ \lim_{t \ra 0}\WL(\mu_t,\mu_0) = 0 $ we have that $ \cJ[ \muc ]\geq 0 $.
Equality holds if and only if $ \muc  $ is the gradient flow for $ \cF $ with initial value $ \mu_0 $.
  \et
  \bpr
  Since $ \cF $ is $ (\bar{\lambda} + \hat{\lambda})$-convex,  we can apply \cite[2.4.10]{ambgigsav} to see  that
      \be{LyapunovEq1}
 \cF(\mu_\e) - 		\cF(\mu_T) \leq  \int_\e^T |\partial \cF|(\mu_t) \,  |\mu'|(t) \,dt   \qquad \text{ for all $ \e \in (0,T) $.}
  \ee
  Thus, Young's inequality and the $\WL$-l.s.c.\ of $ \cF $ yield the first claim.
  The ``if''-part of the second claim is the content of Theorem \ref{ExistenceThm} \emph{(iii)}.
  To show  the ``only if''-part,
  assume that $ \cJ[ (\mu_t)_t ] = 0 $. 
  Hence, $ |\partial \cF|(\mu_t) < \infty   $ for almost every $ t $ and  Proposition \ref{SubdiffCharacF} is applicable.
  Let $ (v_t)_t $ be the tangent velocity of $ \muc $ and $ (\rho_t)_t $ be the curve of the probability densities of $ \muc $. Recall that $ \|v_t \|^2_{\L^2(\mu_t)} = |\mu'|(t) $ for a.e.\  $ t $.
  Then, using the chain rule from Lemma \ref{ChainRuleLem} and the characterisation of the metric slope from Proposition \ref{SubdiffCharacF}, we obtain that 
        \be{LyapunovEq2}
\frac{1}{2} \int_0^T \left\|v_t + 	 \frac{\partial_\t \rho_t}{\rho_t} + \Psi'	 - 	 \int_{  \TdR}J(\cdot-\bar{x})\bar{\t}      \, d\mu_t		\right\|^2_{\L^2(\mu_t)}	\,dt  =\cJ[ (\mu_t)_t ] = 0,
  \ee
  which, again by Proposition \ref{SubdiffCharacF}, implies that $ -v_t \in  \partial \cF (\mu_t)$ for a.e.\ $ t $. Therefore, $ \muc  $ is the gradient flow for $ \cF $. 
  \epr

  \paragraph{Local McKean-Vlasov equation.}\label{BridgeToPde}
  Now we are able to build the bridge to  \eqref{pdeEq} in the following theorem. 
  
  \bt{pde}
  Let $ \mu_0 \in D(\cF) $.
  Let $ T \in (0, \infty) $ and $ (\mu_t)_{t\in[0,T]}  \subset \PL $ be  such that $ \lim_{t\ra 0} \WL(\mu_t,\mu_0)=0 $.
  Then  $ \muc $ is the gradient flow for $ \cF $ if and only if 
  \begin{enumerate}[(i)]
  	\item 
  	$ \mu_t=\rho_t \, \Leb{\TdR} $ for all $ t \in [0,T] $,
  	\item   the curve of densities $ (\rho_t)_t$ is a weak solution to 
  	\be{pdeEq1}  
  	\partial_t \rho_t(x,\t) = \partial_{\t\t}^2 \rho_t(x,\t) + \partial_\t \left(	\rho_t(x,\t) \left( \Psi'(\t) - 	\int J(x-\bar{x})\bar{\t}   \rho_t(\bar{x},\bar{\t})   \, d\bar{\t} d	\bar{x}	\right)					\right),
  	\ee 
  	\item $ \int_0^T |\partial \cF|^2(\mu_t) \,dt< \infty $. 
%
  \end{enumerate}  
  \et
  \bpr
  If $ \muc $ is the gradient flow for $ \cF $, then Theorem \ref{ExistenceThm} \emph{(v)} and \emph{(iii)} imply the claims \emph{(i)} and \emph{(iii)}, respectively.
  Claim  \emph{(ii)} follows immediately from Proposition \ref{SubdiffCharacF} and the fact that $ \muc $ satisfies the continuity equation.
  
  Conversely,  assume \emph{(i)}--\emph{(iii)}.
   \emph{(iii)} implies that $ |\partial \cF|(\mu_t)  < \infty  $ for almost every $ t $.
   Therefore, Proposition \ref{SubdiffCharacF} is applicable and we obtain that for almost every $ t $, $ \partial_\t \rho_t $ exists weakly 
    and 
\begin{align}
    w_t :=  \frac{\partial_\t \rho_t}{\rho_t}+  \Psi' - 	 \int J(\cdot-\bar{x})\bar{\t}      \, d\mu_t(\bar{x},\bar{\t})  \, \in \, \TanL{\mu_t}\,  \cap \, \partial \cF (\mu_t)   .
\end{align}
Moreover, 
\emph{(ii)} shows that $ \muc $ solves the continuity equation with respect to $ (-w_t)_t $.
And since \emph{(iii)} also shows that $ w \in \L^2((0,T) \times \TdR ; \mu_t dt ) $, we infer via Proposition \ref{AbsContProp} \emph{(B)} that $ \muc \in \ACLT $ with tangent velocity $ -w $.
And since $ w_t\in \partial \cF (\mu_t)  $  for almost every $ t $, we conclude that $ \muc  $ must be the gradient flow for $ \cF $ with initial value $ \mu_0 $.  
  \epr

%
%

\section{Large deviation principle}
\label{LDPChap}

In this chapter we  derive the large deviation principle for the system from Section \ref{S1.1.1}.
 First, in Section \ref{3.0.1}, we rigorously introduce the model and state some properties. 
 Then we define  in Section \ref{3.0.2} the empirical measure map and in   Section \ref{3.0.3}  we state the main result and its proof. 
 The proof of the lower bound and the recovery sequence are moved to Section \ref{LowerBoundSec} and Section \ref{UpperBoundSec}, respectively.
 For convenience purposes, from now on we restrict to the case $ d=1 $.
 Throughout the remaining part of this paper 
 suppose Assumption \ref{Fass} and let   $ T \in (0,\infty) $. 
 
\subsection{The microscopic system} \label{3.0.1}

Let $ N \in \N $. Recall the definition of $ H^N $ in \eqref{Hamiltonian}. Define 	$ \cH^N : \cP_2(\R^N) \ra (-\infty , \infty] $ by 
\be{RelEntN} 
\cH^N(\cdot) := \cH (\cdot \, | \, \exp(H^N) \Leb{\R^N} )	 .
\ee
Analogously to \eqref{LyapunovEq}, define $ \cJ^N : C([0,T];\cM_1(\R^N))  \ra [0,\infty] $ by 
\be{LyapunovNEq}
\cJ^N[ (\nu_t^N)_{t} ] :=
\cH^N(\nu^N_T) - 		\cH^N(\nu^N_0) + \frac 12 \int_0^T\big(|\partial \cH^N|^2(\nu^N_t)  + |(\nu^N)'|^2(t)\big)\,dt
\ee 
if $ (\nu_t^N)_{t} \in \AC{(0,T)\,;\, \cP_2(\R^N)} $ and $ \cJ^N[ (\nu_t^N)_{t\in[0,T]} ] = \infty $ else. 

\bl{WassGrad}
Recall the parameters from Assumption \ref{Fass}, Lemma \ref{LBF} and Theorem \ref{FProp}. Then, 
\begin{enumerate}[(i)]
		\item $ \cH^N $ is  proper, $ (\bar{\lambda} + \hat{\lambda})$-convex, strongly $ \lambda $-convex,  $W_2$-l.s.c.\ and coercive,
	\item  for all $ \nu^N \in \cM_1(\R^N)  $,  for some constant $ C'' >0 $,
 \be{LBHNEq}
  \frac{1}{N} \cH^N(\nu^N)  \geq  \frac{1}{N}  \int_{  \R^N	} \Big(C_\Psi\, \sum_{i=0}^{N-1}|\t^i|^{2\ell}  +( C_\Psi' - \|J\|_\infty )\, |\varTheta|^2\Big) d\nu^N(\varTheta) - C'', 
\ee
	\item for all $ \mu_0^N \in D(\cH^N) $, there exists a unique curve $ (\mu_t^N)_{t} \in \AC{(0,T)\,;\, \cP_2(\R^N)} $ such that  $ \lim_{t \ra 0} W_2(\mu_t^N,\mu_0^N)=0 $ and 
	$ \cJ^N[ (\mu_t^N)_{t} ] = 0 $.
	We call $ (\mu_t^N)_{t} $   \emph{the Wasserstein gradient flow for 
	$ \cH^N $} with initial value $ \mu_0^N $,
	\item there exists  $ Q^N \in   \cM_1(C([0,T]; \R^N)) $ such that $  (e_t)_\#Q^N = \mu^N_t  $ for all $ t\in [0,T]   $ and 
	$ Q^N $ is the law of the  trajectories on $ [0, T] $ of a solution to 
\begin{align}\label{EqSpinSystemN}	
d\Theta_t^{N} = - \nabla H^N (\Theta_t^N) \,  dt  + \sqrt{2} \, dB_t^{N} \quad \text{and} \quad 
\Theta^N_0 \sim \mu_0^N.
\end{align} 
\end{enumerate}
\el
\bpr
\emph{(i)} follows from  \cite[9.3.9]{ambgigsav}, \cite[9.3.2]{ambgigsav} and  \cite[9.2.7]{ambgigsav}.

To show \emph{(ii)}, let, without restriction,  $ \nu^N \in \cM_1(\R^N)  $ be such that $ \cH^N(\nu^N) < \infty$.  
Then,  
\begin{align}\label{WassGradEqq1}
\begin{split}
\frac{1}{N} \cH^N(\nu^N) &= \frac{1}{N}\cH\left ( \mu^N\, \Bigg| \,  \exp\left(-\frac{1}{2} \sum_{k=0}^{N-1} \Psi(\t^k)\right) d\varTheta \right) \\
&\ \ \ + 
\frac{1}{2N^2} \sum_{k,j=0}^{N-1}  \int_{  \R^N	}\left(\frac{1}{2} \Psi(\t^k) +  \frac{1}{2}\Psi(\t^j) - J\left(\frac{k-j}{N} \right)  \t^{k} \t^{j}\right) d\nu^N(\varTheta) .
\end{split}
\end{align}
Proceeding  as in the proof of Lemma \ref{LBF} 
yields part \textit{(ii)}.

\emph{(iii)}   is a consequence of \cite[11.2.1]{ambgigsav} and part \emph{(i)}. 

\emph{(iv)}  follows from \cite[3.3]{PatLDP}.
\epr

For technical reasons we have to restrict the choice of the sequence $ (\mu_0^N )_N $ of initial values in the following way.  
\bass{AssI}
For all  $ N \in \N $,  $ \mu_0^N \in \cM_1(\R^N) $ is given by $ d\mu_0^N(\varTheta) = \rho_0^N(\varTheta)\,  d \varTheta $, where 
\begin{align}
\rho_0^N(\varTheta) := \prod_{k=0}^{N-1} \kappa\left(\tfrac{k}{N}, \t^k\right) \,  e^{-\Psi(\t^k)} \,,
\end{align}
where $ \kappa:\ToR\ra [0,\infty) $ is upper semi-continuous and  such that 
\begin{itemize}
	  	\setlength\itemsep{-0.1em}
	\item $ \int_{ \R	}\kappa(x,\t) \,  e^{-\Psi(\t)  } \, d\t =1  \text{ for each } x \in \To $,
	\item the restriction $ \kappa: \{\kappa>0\} \ra (0,\infty) $ is a continuous map, where $ \{\kappa>0\} := \{(x,\t) \in \ToR \, | \, 	\kappa(x,\t)>0	\}  $,
	\item $ \kappa(x,\t) \leq C_\kappa \, \exp(\frac{1}{8}C_\Psi \, \t^{2\ell} + \frac{1}{8}(C_\Psi'-\|J\|_\infty) \t^{2} ) $ for some $ C_\kappa >0  $, and 
	\item  either $ \kappa(x,\t) \geq c_\kappa'\,   \exp(-c_\kappa \Psi(\t)  ) $    on  $ \{\kappa>0\}   $ for some $ c_\kappa,c_\kappa' >0  $, or  $ x \to \kappa(x,\t)  $ is constant for all $ \t \in \R $.
\end{itemize}
\eass 
%
\subsection{The empirical measure map.} \label{3.0.2}
For all $ N \in \N $, define the \emph{empirical measure map} $ K^N : \R^N \ra \cM_1(\Tone \times \R) $ by 
  \begin{align}
K^N(\varTheta) = \frac{1}{N} \sum_{k=0}^{N-1} \delta_{\left(\frac{k}{N}, \theta^k\right)}.
\end{align}
Moreover, let 
$ K_T^N: C([0,T] \, ; \, \R^N )  \ra C([0,T] \, ; \, \cM_1(\ToR) )  $ be defined by  
  \begin{align}
K^N_T((\varTheta_t)_{t\in[0,T]}) = (K^N(\varTheta_t))_{t\in[0,T]}.
\end{align}
For technical reasons it will be useful to consider  a  modification of $ K^N $ defined by 
  \begin{align}
\begin{split}
L^N : \R^N &\ra \cM_1^\L(\Tone \times \R)	\\
\varTheta &\to  \sum_{k=0}^{N-1} \Leb{A_{k,N}} \otimes  \delta_{ \theta^k},
\end{split}
\end{align}
where $ (A_{k,N})_{k=0}^{N-1} $ is a
  partition of $ \To $ given by
\be{AkN}
A_{k,N} = [kN, (k+1)N), \qquad k=0,\dots, N-1.
\ee
In the following lemma we show that $ L^N $ is indeed just a small modification of $ K^N $. 
\bl{KvsL}
\begin{enumerate}[(i)]
	\item 
	Let $ N \in \N  $ and $ \varTheta \in \R^N$. Then 
	$ W_2 (K^N(\varTheta), L^N(\varTheta)) \leq \frac{1}{N} $.
	\item Let $ \widetilde{\bW} $ denote the Wasserstein distance on $ \cM_1(\cM_1(\ToR)) $ induced by the distance $ \widetilde{W}   $ on  $ \cM_1(\ToR)  $. 
	Then
	\begin{align}
		\widetilde{\bW} \left( (L^N)_\# \mu^N , (K^N)_\# \mu^N \right ) \leq \frac{1}{N}  \qquad \forall \,  \mu^N \in \cM_1(\R^N) .
	\end{align}
\end{enumerate}
\el
\bpr
Define $ G:\ToR \ra \ToR $ by 
\be{KvsLEq} G(x,\t)  = \sum_{k=0}^{N-1} \Ind_{A_{k,n}}(x) \left(\frac{k}{N}, \t\right)  .
\ee
Then $ (G,\Id_{\ToR})_\# L^N(\varTheta) \in \CouplW{K^N(\varTheta), L^N(\varTheta)} $.
Estimating $ W_2 (K^N(\varTheta), L^N(\varTheta))  $ by the cost with respect to this coupling yields \emph{(i)}.
Finally, 
\emph{(ii)}  follows immediately from part \emph{(i)}.
\epr

\subsection{The large deviation principle.} \label{3.0.3}

\bd{LDPDefi}
Let $(X, \mathrm{d})$ be a Polish space. 
Let $(\varPi_n)_{n \in \N}$ be a family of probability measures on $ X$ and let $I:X \ra [0,\infty]$ be $ \mathrm{d} $-l.s.c. 
Then $(\varPi_n)_n$ is said to satisfy a \emph{large deviation principle (LDP)} on $X$ with \emph{rate function} $I$ if 
\begin{enumerate}[(i)]
	\item for any closed set $ C \subset X$,  $ \limsup_{n \ra \infty} \frac 1{n} \log \varPi_n(C) \leq  -\inf_{x\in C} I(x) $, and	
	\item for any open set $ O \subset X$, $ \liminf_{n \ra \infty} \frac 1{n} \log \varPi_n(O) \geq  -\inf_{x\in O} I(x) $.
\end{enumerate}
\ed

Recall that $ \cM_1(\ToR) $ is equipped with the metric $ \widetilde{W} $.
Let $ C([0,T] \, ; \, \cM_1(\ToR) ) $ be equipped with the supremum norm induced by $ \widetilde{W} $. 
Theorem \ref{LDP} below states the LDP result for the sequence $ \{	(K_T^N)_\# Q^N	\}_N  $ on $  C([0,T] \, ; \, \cM_1(\ToR) )   $, where, for all $ N \in \N $,  $ Q^N $ is the measure from Lemma \ref{WassGrad} \emph{(iv)}.
The rate function will be given by 
\begin{align}\label{RateFct}
I[ (\nu_t)_t ] = 
\begin{cases}
\frac{1}{2} \cJ[ (\nu_t)_{t} ] + \cH(\nu_0\, | \, \mu_0)  &\text{if }(\nu_t)_t \in \AC{[0,T]\, ; \, \PLo},\\
\infty		&\text{else,}
\end{cases}
\end{align}
where $ \mu_0 = \rho_0 \,  \Leb{\ToR}\in \PLo $ with  $ \rho_0(x,\t) = \kappa(x,\t) \,  e^{-\Psi(\t)  } $.
Before we state and prove the LDP result, we need to show the lower semi-continuity of $ I  $. 
\bl{Ilsc}
$ (\nu_t)_{t} \to I[(\nu_t)_{t}] $ is  lower semi-continuous in $ C([0,T] \, ; \, \cM_1(\ToR) ) $.
\el
\bpr
Let $ \lim_{m\ra \infty} (\nu^m_t)_{t}  = (\nu_t)_{t} $ in $ C([0,T] \, ; \, \cM_1(\ToR) ) $.  
In particular,  $\nu^m_t \hra \nu_t  $ for all $ t \geq 0 $. 
Without restriction assume that $ \liminf_{m\ra \infty} I[(\nu^m_t)_{t}]  < \infty $, since otherwise, the claim is trivial.
Moreover, by considering appropriate subsequences, we can  even suppose that $ \sup_{m\in \N} I[(\nu^m_t)_{t}]  < \infty $.
In particular, $ \sup_{m\in \N}\cJ[ (\nu^m_t)_{t} ] $, $ \sup_{m\in \N}\cH(\nu^m_0\, | \, \mu_0) < \infty $, since both terms are non-negative.
The proof is divided into seven steps. 

\noindent
\textbf{Step 1.} [ $ \inf_{m\in \N}\cH(\nu_0^m	\, | \, \mu_0) -  \tfrac{1}{2}\cF(\nu_0^m) > - \infty  $. ]

\noindent
Note that, since $ \sup_{m\in \N} \cH(\nu^m_0\, | \, \mu_0) < \infty  $, $ \kappa $ is strictly positive inside the support of $ \nu_0^m $. 
Then, similarly as in  the proof of Theorem \ref{FProp} 
\begin{align}\label{ILSCEq2}
& \cH(\nu_0^m	\, | \, \mu_0) -  \tfrac{1}{2}\cF(\nu_0^m)  			= \nonumber
\frac{1}{2} \cH\Big( \nu_0^m\, \Big| \,   e^{-\frac{1}{2}  \Psi(\t)} d\t \Big)
\\
&\qquad      +\frac{1}{4}  \int_{  (\ToR)^2	} \Big(	\tfrac{1}{2} \Psi(\t) +\tfrac{1}{2} \Psi(\bar{\t}) 	+ J(x - \bar{x}) \t \bar{\t} - 2 \log	\kappa(x,\t)- 2\log \kappa(\bar{x},\bar{\t})				 	\Big) \, d(\nu_0^m \otimes \nu_0^m) \nonumber\\
&\quad\geq -\frac{1}{2} \log \int   e^{-\frac{1}{2} \Psi(\t)} d\t- \frac{1}{4}C_\Psi'' 
-\log(C_\kappa)   > -\infty.
\end{align}
By using Assumption \ref{Fass} and Assumption \ref{AssI}, we have that
\begin{align}\label{ILSCEq333}
 &\cH\Big( \nu_0^m\, \Big| \,   e^{-\frac{1}{2}  \Psi(\t)} d\t \Big) \geq 
 -  \log \int   e^{-\frac{1}{2} \Psi(\t)} d\t, \quad \text{ and} \\ \label{LBInitBdd1}
& \tfrac{1}{2} \Psi(\t) +\tfrac{1}{2} \Psi(\bar{\t}) 	+ J(x - \bar{x}) \t \bar{\t} - 2 \log	\kappa(x,\t)- 2\log \kappa(\bar{x},\bar{\t})	\geq -  C_\Psi'' 
-4\log(C_\kappa)   .
\end{align}
Combining \eqref{ILSCEq2}, \eqref{ILSCEq333} and \eqref{LBInitBdd1} concludes the claim of Step 1.

\noindent
\textbf{Step 2.} [ $ \sup_{m\in \N}  \int_0^{T}   |(\nu^m)'|^2(r) \, dr  < \infty $ and $ \sup_{m\in \N}  \int_0^{T}     |\partial \cF|^2(\nu_r^m) \, dr  < \infty $. ]

\noindent
Using Step 1, the fact that $ \sup_{m\in \N} I[(\nu^m_t)_{t}]  < \infty $, and   Lemma \ref{LBF}, we infer the claim.
 
 \noindent
 \textbf{Step 3.} [ $ \sup_{m\in \N}  \sup_{t\in [0,T]}  \cF (\nu_t^m)    < \infty $ and  $ \sup_{m\in \N}  \sup_{t\in [0,T]} \int_{\ToR}   |\t|^2  \, d\nu^m  < \infty $. ]
 
 \noindent
Since $ 	|\partial \cF| $ is a so-called strong upper gradient  (\cite[1.2.1 and 2.4.10]{ambgigsav}),  we infer that
\be{ILSCEq111}
\sup_{m\in \N}  \sup_{t\in [0,T]}     \cF(\nu_t^m)   \leq \sup_{m\in \N}  \sup_{t\in [0,T]}  \int_t^T |\partial \cF|(\nu_r^m) \,|(\nu^m)'|(r) \,  dr   +   \cF(\nu_T^m) < \infty ,
\ee 
where we used Step 2 in the last step. The second claim is shown by combining \eqref{ILSCEq111} with Lemma \ref{LBF}.

\noindent
 \textbf{Step 4.} [ $ \liminf_{m\ra \infty} \left(\cF(\nu^m_T) + \frac 12 \int_0^T|\partial \cF|^2(\nu^m_t) \,dt \right) \geq \cF(\nu_T) + \frac 12 \int_0^T|\partial \cF|^2(\nu_t) \,dt $. ]\\
The claim follows from a combination of  Theorem \ref{FProp},  Fatou's lemma, Step 3 and Corollary \ref{FSubDiffCor}.

\noindent
\textbf{Step 5.} [ $ \liminf_{m\ra \infty}\cH(\nu_0^m	\, | \, \mu_0) -  \tfrac{1}{2}\cF(\nu_0^m) \geq \cH(\nu_0^m	\, | \, \mu_0) -  \tfrac{1}{2}\cF(\nu_0^m) $. ]

\noindent
Recall   \eqref{ILSCEq2}, and recall that we have already seen in  the proof of Theorem \ref{FProp} that 
\begin{align}
 \liminf_{m\ra \infty} \ \frac{1}{2}   \cH\Big( \nu_0^m\, \Big| \,   e^{-\frac{1}{2}  \Psi(\t)} d\t \Big)  \  \geq  \frac{1}{2}   \cH\Big( \nu_0\, \Big| \,   e^{-\frac{1}{2}  \Psi(\t)} d\t \Big).
\end{align} 
The integrand in the second term on the right-hand side of \eqref{ILSCEq2} is lower semi-continuous and bounded from below by Assumption \ref{Fass} and Assumption \ref{AssI}.  
Therefore,  analogously to \eqref{Flsc3},   \cite[5.1.7]{ambgigsav} yields the lower semi-continuity of this term.  

\noindent
\textbf{Step 6.} [  $ (\nu_t)_t \in \AC{[0,T]\, ; \, \PLo}$. ]

\noindent 
According to \cite[Lemma 1]{Lisini}, it suffices to show that 
\begin{align}\label{ILSCEq10}
\sup_{0<h<T }  \int_0^{T-h} \frac{1}{h^2} \WL(\nu_t, \nu_{t+h})^2   \, dt < \infty \quad \text{and} \quad 
\int_0^{T} \WL(\nu_t, \delta_0\otimes \Leb{\To})^2   \, dt
< \infty.
\end{align}
Since $ \WL(\nu_t, \delta_0\otimes \Leb{\To})^2 = \int |\t|^2 \, d\nu_t $, Step 3 and \cite[5.1.7]{ambgigsav} imply the second claim in \eqref{ILSCEq10}.
In order to show the first claim in \eqref{ILSCEq10},  note that $ (\nu^m_t)_t \in \AC{[0,T]\, ; \, \PLo}$  for all $ m $.
Then, using Fatou's lemma and  Lemma \ref{LSCLem}, we obtain that 
\begin{align}\label{ILSCEq1}
   \sup_{0<h<T }  \int_0^{T-h} \frac{1}{h^2} \WL(\nu_t, \nu_{t+h})^2   \, dt \nonumber
&\leq \,  \sup_{0<h<T } \liminf_{m\ra \infty }  \int_0^{T-h} \frac{1}{h^2} \WL(\nu^m_t, \nu^m_{t+h})^2   \, dt \\
&\leq  \sup_{0<h<T } \liminf_{m\ra \infty }  \int_0^{T-h}  \frac{1}{h} \int_t^{t+h } |(\nu^m)'|^2(r) \, dr \,  dt		\\
&\leq 				\liminf_{m\ra \infty }    \int_0^{T}   |(\nu^m)'|^2(r) \, dr < \infty\nonumber,
\end{align}
where we have used Fubini's theorem in the last step.

\noindent
\textbf{Step 7.} [  $ \int_0^{T }   |\nu'|^2(t) \, dt
  \leq 	 \liminf_{m\ra \infty }    \int_0^{T}   |(\nu^m)'|^2(r) \, dr $. ]

\noindent 
Let $ \varepsilon \in (0, T/2) $. Then, repeating the arguments from \eqref{ILSCEq1},
\begin{align}\label{ILSCEq11}
\begin{split}
\int_0^{T-\varepsilon}   |\nu'|^2(t) \, dt
\, \leq \,   \liminf_{h\downarrow 0, h<\varepsilon }  \int_0^{T-\varepsilon} \frac{1}{h^2} \WL(\nu_t, \nu_{t+h})^2   \, dt  
\leq 				\liminf_{m\ra \infty }    \int_0^{T}   |(\nu^m)'|^2(r) \, dr.
\end{split}
\end{align}
 Letting $ \varepsilon \downarrow 0 $ concludes the proof. 
\epr

\bt{LDP}
Let $ (\mu_0^N)_N $ satisfy Assumption \ref{AssI}.
For all $ N \in \N $, let $ (\mu_t^N)_{t\in[0,T]} $  be the Wasserstein gradient flow for 
$ \cH^N $ with initial value  $ \mu_0^N $ and  $ (Q^N)_N $ be the corresponding representation measures from Lemma \ref{WassGrad} (iv).
Then the sequence $ \{	(K_T^N)_\# Q^N	\}_N $ satisfies  a large deviation principle on $ C([0,T] \, ; \, \cM_1(\ToR) ) $   with rate function 
$ I $.
\et
\bpr
In \cite[Theorem 3.4 and Theorem 3.5]{Mariani}, it is shown that the above LDP result for  $ \{	(K_T^N)_\# Q^N	\}_N $ is true if and only if  the following three conditions are satisfied.
\begin{itemize}
	\item[\emph{(i)}] The family $ \{	(K_T^N)_\# Q^N	\}_N $ is exponentially tight, i.e.\ for all $ s > 0 $ there exists a compact set $ \cK_s \subset C([0,T] \, ; \, \cM_1(\ToR) ) $ such that 
	\be{ExpTight}
	\limsup_{N\ra \infty } \frac{1}{N} \log \Big(	(K_T^N)_\# Q^N (\cK_s^c)	\Big) \leq -s.
	\ee
	\item[\emph{(ii)}] For all $ (\nu_t)_t \in C([0,T] \, ; \, \cM_1(\ToR) ) $ and for all sequences $ (\varGamma^N)_N \subset \cM_1(\, C([0,T] \, ; \, \cM_1(\ToR) ) \, )  $  that converge  to $ \delta_{(\nu_t)_t} $ weakly in $ \cM_1(\, C([0,T] \, ; \, \cM_1(\ToR) ) \, ) $, it holds
	\be{GammaLB}
	\liminf_{N\ra \infty} \frac{1}{N} 	\cH\left(		\varGamma^N	\, \Big| \, 	(K_T^N)_\# Q^N	\right) \geq I[ (\nu_t)_t ].
	\ee
		\item[\emph{(iii)}] For all $ (\nu_t)_t \in C([0,T] \, ; \, \cM_1(\ToR) ) $ there exists $ (\varGamma^N)_N \subset \cM_1(\, C([0,T] \, ; \, \cM_1(\ToR) ) \, )  $   such that $ (\varGamma^N)_N $ converges to $ \delta_{(\nu_t)_t} $ weakly in $ \cM_1(\, C([0,T] \, ; \, \cM_1(\ToR) ) \, ) $ and 
	\be{GammaUB}
	\limsup_{N\ra \infty} \frac{1}{N} 	\cH\left(		\varGamma^N	\, \Big| \, 	(K_T^N)_\# Q^N	\right) \leq I[ (\nu_t)_t ].
	\ee
\end{itemize}

Fact \emph{(i)} was  proven in \cite[3.29]{PatLDP}.

To prove \emph{(ii)}, note that if the left-hand side of \eqref{GammaLB} is infinite, the claim is trivial. Therefore, we assume without restriction that 
\be{GammaLB1}
\cH\left(		\varGamma^N	\, \Big| \, 	(K_T^N)_\# Q^N	\right) < \infty \quad \text{ for all }N \in \N. 
\ee
This  implies in particular that $ 	\varGamma^N  $ is absolutely continuous with respect to $ (K_T^N)_\# Q^N $ for all $ N $. 
Since the map $ K_T^N $ is injective, we infer that for all $ N $ there is a $ P^N \in  \cM_1(\, C([0,T] \, ; \, \R^N ) \, )  $ such that 
$ \varGamma^N = (K_T^N)_\# P^N  $. Moreover, 
\be{GammaLB2}
\cH\left(	(K_T^N)_\# P^N\, \Big| \, 	(K_T^N)_\# Q^N	\right) = \cH\left(		P^N	\, \Big| \,  Q^N	\right)  \quad \text{ for all }N \in \N,
\ee
which is again a consequence of the injectivity of  $ K_T^N $.
Now we can use \cite[4.1.(i)]{FathiLDP} to observe that
\be{GammaLB3}
 \cH\left(		P^N	\, \Big| \,  Q^N	\right) \geq   \frac{1}{2} \cJ^N[ (\nu^N_t)_{t} ] + \cH\left(		\nu^N_0	\, \Big| \,  \mu^N_0	\right)  \quad \text{ for all }N \in \N,
\ee
where $ \nu^N_t := (e_t)_\# P^N $ for all $ t  $.
In particular, the right-hand side in \eqref{GammaLB3} is finite, which implies that $  (\nu^N_t)_{t} \in \AC{[0,T]; \PwN} $. 
Hence, in order to prove $ (ii) $, it will be enough to show that 
	\be{GammaLB4}
\liminf_{N\ra \infty} \frac{1}{N} \left(  \frac{1}{2} \cJ^N[ (\nu^N_t)_{t} ] + \cH\left(		\nu^N_0	\, \Big| \,  \mu^N_0	\right) \right)  \geq I[ (\nu_t)_t ],
\ee
whenever $ ((K^N)_\# \nu^N_t)_N $ converges to $  \delta_{\nu_t} $ weakly in $ \cM_1(  \cM_1(\ToR)  ) $  for all $ t $, where $  (\nu^N_t)_{t} \in \AC{[0,T]; \PwN} $ for all $ N $. This is the content of Proposition \ref{LowerBound} below.  

It remains to prove \emph{(iii)}.
If $ I[ (\nu_t)_t ] = \infty $, we  take $ \varGamma^N =\delta_{(\nu_t)_t} $ for all $ N $ and \eqref{GammaUB} is trivially satisfied.
So assume that 
$ I[ (\nu_t)_t ] < \infty $. 
In particular, $ (\nu_t)_{t}  \in \AC{[0,T] ; \PLo}$. 
Proposition \ref{UpperBound} below shows  that there exists 
$ (\nu_t^N)_{t}  \in C([0,T] ; \cM_1(\R^N) )$ 
such  that 
$ ((K^N)_\# \nu^N_t)_N $ converges to $  \delta_{\nu_t} $ weakly in $ \cM_1(  \cM_1(\ToR)  ) $ for all $ t $  and
\be{UBPEq}
\limsup_{N\ra \infty }  \frac{1}{N} \left(  \frac{1}{2}\cJ^N[ (\nu_t^N)_{t} ]  + \cH(\nu_0^N	\, | \, \mu_0^N 	)	\right)	\leq   I[ (\nu_t)_{t} ] .
\ee
Further, for all $ N $, \cite[4.1.(ii)]{FathiLDP} yields the existence of $ \tilde{P}^N \subset \cM_1(\, C([0,T] \, ; \, \R^N ) \, )  $   such that 
\be{UBPEq1}
  \frac{1}{2}\cJ^N[ (\nu_t^N)_{t} ]  + \cH(\nu_0^N	\, | \, \mu_0^N 	) = \cH\left(		\tilde{P}^N	\, \Big| \,  Q^N	\right)
\ee
and $ \nu^N_t = (e_t)_\# \tilde{P}^N $ for all $ t  $.
Hence, in order to prove \emph{(iii)}, it only remains to show that $ (\, (K_T^N)_\# \tilde{P}^N\, )_N   $ converges to $ \delta_{(\nu_t)_t} $ weakly in $ \cM_1(\, C([0,T] \, ; \, \cM_1(\ToR) ) \, ) $.

Since $ ((K^N)_\# \nu^N_t)_N $ converges to $  \delta_{\nu_t} $ weakly in $ \cM_1(  \cM_1(\ToR)  ) $ for all $ t $, it suffices to show that $(\,  (K_T^N)_\# \tilde{P}^N\, )_N $ is tight.
Let $ \e > 0 $.
Let $ \cK $ be the compact set from part \emph{(i)} according to the choice $ s = I[ (\nu_t)_t ] / \e $.
Then, via the entropy inequality (see e.g.\ \cite[(3.7)]{Mariani}), \eqref{ExpTight}, \eqref{GammaLB2}, \eqref{UBPEq} and \eqref{UBPEq1}, we obtain
\begin{align}\label{UBPEq2}
\begin{split}
\limsup_{N\ra \infty } \, 	(K_T^N)_\# \tilde{P}^N (\, \cK^c\, )	 &\leq \limsup_{N\ra \infty } \frac{\log 2 + \cH\left(	(K_T^N)_\# \tilde{P}^N\, \Big| \, 	(K_T^N)_\# Q^N	\right)}{\log\left(	1+ 1/ 	(K_T^N)_\# Q^N (\, \cK^c\, )		\right)} \\
&\leq \limsup_{N\ra \infty } \frac{\frac{1}{N} \cH\left(	(K_T^N)_\# \tilde{P}^N\, \Big| \, 	(K_T^N)_\# Q^N	\right)}{- \frac{1}{N} \log\left(		(K_T^N)_\# Q^N (\, \cK^c\, )		\right)} \quad \leq  \ \e,
\end{split}
\end{align}
which implies the tightness of $(\,  (K_T^N)_\# \tilde{P}^N\, )_N $.
\epr


\subsection{Lower Bound}
\label{LowerBoundSec}

\bp{LowerBound}
Let $ (\nu_t)_{t}  \in C([0,T] ; \cM_1(\ToR)) $ and $  (\nu^N_t)_{t} \in \AC{[0,T]; \PwN} $ for all $ N \in \N $. 
Suppose  that $ ((K^N)_\# \nu^N_t)_N $ converges to $  \delta_{\nu_t} $ weakly in $ \cM_1(  \cM_1(\ToR)  ) $ 
for all $ t \in[0,T] $.
Then 
\be{LBEq}
\liminf_{N\ra \infty }  \frac{1}{N} \left(  \frac{1}{2}\cJ^N[ (\nu_t^N)_{t} ]  + \cH(\nu_0^N	\, | \, \mu_0^N 	)	\right)	\geq  \frac{1}{2} \cJ[ (\nu_t)_{t} ] + \cH(\nu_0\, | \, \mu_0).
\ee
\ep
\bpr
Assume that the left-hand side of \eqref{LBEq} is finite. Otherwise, the claim is trivial. 
In particular,  since both summands are non-negative, we have
\be{LBAssEq}
\liminf_{N\ra \infty }  \frac{1}{N}    \cJ^N[ (\nu_t^N)_{t} ]  < \infty \qquad \text{and} \qquad 
\  \liminf_{N\ra \infty }  \frac{1}{N}  \cH(\nu_0^N	\, | \, \mu_0^N 	)< \infty.
\ee
Under this assumption, we show \eqref{LBEq} for each part separately in the forthcoming paragraphs. Hence, the claim
follows from the Lemmas \ref{LBMcK}--\ref{LowerBoundMetrSlope}. 
\epr

\paragraph{Preliminaries.}
We first list some  consequences of \eqref{LBAssEq} in the following lemma.

\bl{LowerBoundAbsCont}
Under the same assumptions as in Proposition \ref{LowerBound} and under \eqref{LBAssEq}, we have  
\be{FinMEq}
\liminf_{N\ra \infty }  \frac{1}{N} \int_0^T |(\nu^N)'|^2(t) \, dt < \infty.
\ee
Moreover, $ (\nu_t)_t $ is an absolutely continuous curve in $ \PLo $. 
\el
\bpr
\textbf{Step 1.} \big[ $  \inf_{N\in \N}  \frac{1}{N} \left(  \cH(\nu_0^N	\, | \, \mu_0^N) - \tfrac{1}{2} \cH^N(\nu_0^N)  		\right) > -\infty $.  \big]

\noindent  
Analogously to \eqref{ILSCEq2} and in view of \eqref{WassGradEqq1} and \eqref{LBInitBdd1} we observe that for all $ N \in \N $

\begin{align}\label{LBInitBd1}
 &\frac{1}{N} \big(  \cH(\nu_0^N	\, | \, \mu_0^N) -  \tfrac{1}{2}\cH^N(\nu_0^N)  		\big)	
= 
\frac{1}{2N} \cH\Big( \nu_0^N\, \Big| \,  \exp\Big(-\frac{1}{2} \sum_{k=0}^{N-1} \Psi(\t^k)\Big) d\varTheta \Big) \nonumber \\
&\quad      +\frac{1}{4} \int_{  \cM_1(\ToR)} \int_{  (\ToR)^2	} \Big(	\tfrac{1}{2} \Psi(\t) +\tfrac{1}{2} \Psi(\bar{\t}) 	+ J(x - \bar{x}) \t \bar{\t}  \nonumber\\
&\qquad\qquad\qquad\qquad\qquad\qquad  - 2\log\kappa(x,\t) - 2\log\kappa(\bar{x},\bar{\t})				 	\Big) \, d\gamma \,  d\gamma \, d(K^N)_\# \nu_0^N(\gamma)\\
&\geq   \frac{1}{2N} \cH\Big( \nu_0^N\, \Big| \,  \exp\Big(-\frac{1}{2} \sum_{k=0}^{N-1} \Psi(\t^k)\Big) d\varTheta \Big) - \frac{1}{4}C_\Psi'' 
 -\log(C_\kappa)  \nonumber \\
&\geq -\frac{1}{2} \log \int  e^{-\frac{1}{2} \Psi(\t)} d\t- \frac{1}{4}C_\Psi'' 
-\log(C_\kappa)   > -\infty.\nonumber
\end{align}

\noindent  
\textbf{Step 2.} \big[  $  \liminf_{N\ra \infty }  \frac{1}{N} \int_0^T |(\nu^N)'|^2(t) \, dt  < \infty $. \big]

\noindent  
Step 1, Lemma \ref{WassGrad} \emph{(ii)} and the finiteness of the  left-hand side of \eqref{LBEq} yield the claim. 


\noindent
\textbf{Step 3.} \big[ $  \liminf_{N\ra \infty }   \frac{1}{N}  \int_{  \R^N	} (|\varTheta|^2 +  \sum_{i=0}^{N-1}|\t^i|^{2\ell}) \, d\nu^N_0  < \infty $. \big]

\noindent 
By a similar computation as in Step 1, \eqref{LBAssEq} yields that 
\begin{align}\label{LB2moment}
\begin{split}
	\infty &> \liminf_{N\ra \infty }  \frac{1}{N} \cH(\nu_0^N	\, | \, \mu_0^N)
	=  \liminf_{N\ra \infty }  \frac{1}{N} 
\cH\Big( \nu_0^N\, \Big| \,  \exp\Big(-\frac{1}{2} \sum_{k=0}^{N-1} \Psi(\t^k)\Big) d\varTheta \Big) 
	\\
	&\quad   + \liminf_{N\ra \infty }  \frac{1}{N}  \sum_{k=0	}^{N-1}  \Bigg(
	- \int_{  \R^N}  \log	\kappa(\tfrac{k}{N}, \t^k) 	d\nu_0^N
	+ \frac{1}{2}
	\int_{  \R^N} \Psi (\t^k)d\nu_0^N
	\Bigg) \\ 
	&\geq C \ + \  \frac{1}{4} \  \liminf_{N\ra \infty }  \frac{1}{N} \int_{  \R^N	} \sum_{k=0	}^{N-1} \Psi (\t^k) \, d\nu^N_0 
\end{split}
\end{align}
for some $ C\in \R  $. 
Finally,  \eqref{LBPsi} implies Step 3.

\noindent
\textbf{Step 4.} \big[  $  \liminf_{N\ra \infty }  \frac{1}{N} \sup_{t\in[0,T]} \int_{  \R^N} |\varTheta|^2 \, d\nu^N_t    < \infty $. \big]

\noindent
Step 2 and 3 imply that 
\begin{align} 
\nonumber
\liminf_{N\ra \infty }  \frac{1}{N}\sup_{t\in[0,T]} &\int_{  \R^N}  |\varTheta|^2 \, d\nu^N_t   
\leq 
\liminf_{N\ra \infty } 4\frac{1}{N} \left(  \sup_{t\in[0,T]} W_2 (\nu_t^N, \nu_0^N)^2 dt +   
\int_{  \R^N} |\varTheta|^2 \, d\nu^N_0  \right) \\
& \leq 
\liminf_{N\ra \infty } 4\frac{1}{N} \left( T   \int_0^T |(\nu^N)'|^2(t) dt +   
\int_{  \R^N} |\varTheta|^2 \, d\nu^N_0  \right) 
< \infty .  
\end{align}    
\textbf{Step 5.} \big[  $ \int_0^T \int_{  \ToR } |\theta|^2 \  d\nu_t(x,\t) \, dt  < \infty $. \big]

\noindent
Using that $ ((K^N)_\# \nu^N_t)_N $ converges to $  \delta_{\nu_t} $ weakly in $ \cM_1(  \cM_1(\ToR)  ) $, and using  \cite[5.1.7]{ambgigsav}, Fatou's lemma and Step 4 we obtain that 
\begin{align}
\int_0^T \int_{  \ToR } |\theta|^2 \ d\nu_t \,   dt 
\leq 
\liminf_{N\ra \infty }  \frac{1}{N} \int_0^T \int_{  \cM_1(\ToR) } \int_{  \ToR } |\theta|^2 \, d\gamma \ d(K^N)_\#\nu^N_t(\gamma ) \, dt  < \infty . 
\end{align}

\noindent
\textbf{Step 6.} \big[ $  \nu_t \in \MLo $ for all $ t \in [0,T] $.  \big]

\noindent
Since $ ((K^N)_\# \nu^N_t)_N $ converges to $  \delta_{\nu_t} $ weakly in $ \cM_1(  \cM_1(\ToR)  ) $, we have that for all $ f \in C_b(\To) $
\begin{align}
\begin{split}
\int_{\ToR} f(x) \, d\nu_t(x,\t) &= 
\lim_{N\ra \infty }  \int_{\cM_1(\ToR)}
\int_{\ToR} f(x) \,  d\gamma(x,\t)  \,d(K^N)_\# \nu_t^N (\gamma) \\
&= 
\lim_{N\ra \infty }  
\int_{\ToR} f\left( \frac{\lfloor  xN \rfloor }{N}
\right) \,  dx = \int_{\ToR} f\left( x
\right) \,  dx.
\end{split}
\end{align}
\textbf{Step 7.} \big[ $ t \mapsto \nu_t $ is absolutely continuous.  \big]

\noindent
Analogously to the proof of Lemma \ref{Ilsc}, it suffices to show that
\begin{align}\label{LBAbsContEq}
\sup_{0<h<T }  \int_0^{T-h} \frac{1}{h^2} \WL(\nu_t, \nu_{t+h})^2   \, dt < \infty \qquad \text{and} \qquad 
\int_0^{T} \int_{\ToR} |\t|^2 \, d\nu_t(x,\t)   \, dt
< \infty.
\end{align}
The second claim was shown in Step 5.
The first claim in \eqref{LBAbsContEq} follows from similar arguments as in \eqref{ILSCEq1}.
Indeed, using Lemma \ref{KvsL} and the Lemmas \ref{LowerSemcWL} and \ref{cWW} below, we observe that  
\begin{align}\label{LBAbsContEq0}
\sup_{0<h<T} \int_0^{T-h} \frac{1}{h^2} \WL(\nu_t, \nu_{t+h})^2  dt 
&\leq   \sup_{0<h<T} \int_0^{T-h} \liminf_{N\ra \infty } \frac{1}{h^2} \bW^\L \Big( (L^N)_\#\nu_{t}^N\, , \, (L^N)_\# \nu_{t+h}^N \Big)^2    dt \nonumber \\
&\leq \sup_{0<h<T} \int_0^{T-h} \liminf_{N\ra \infty } \frac{1}{h^2N} W_2(\nu_{t}^N, \nu_{t+h}^N)^2 \,  dt \\
&\leq 				\liminf_{N\ra \infty }   \frac{1}{N} \int_0^{T}   |(\nu^N)'|^2(r) \, dr < \infty \nonumber,
\end{align}
where we have used Fatou's lemma, Fubini's theorem and Step 2. We conclude the proof.
\epr
\vspace{-0.1cm}
\bl{LowerSemcWL}
Let $ \bW^\L $ denote the Wasserstein distance on $ \cM_1(\MLo) $ induced by $ \WL $.
Let $ (A^N)_N, (B^N)_N \subset  \cM_1(\MLo) $ and $ A, B\in  \cM_1(\MLo) $ be such that $ A^N $ converges to $ A $ and $ B^N $ converges to $ B $ weakly in $ \cM_1(\MLo) $.
Then 
\be{LowerSemcWLEq}
\liminf_{N\ra \infty } \bW^\L(A^N,B^N) \geq  \bW^\L(A,B).
\ee
\el
\bpr
In view of  Lemma \ref{LSCLem}, the claim is an application of \cite[4.3]{vil}.
\epr

\bl{cWW}
Let   $ \mu^N, \nu^N \in  \PwN $. Then 
\be{cWWEq}
\bW^\L \Big( (L^N)_\#\mu^N\, , \, (L^N)_\# \nu^N \Big)   \leq \frac{1}{\sqrt{N}} W_2(\mu^N, \nu^N).
\ee
\el
\bpr
Let $ \pi^N \in \OptW{\mu^N, \nu^N} $. 
Define 
\begin{align}\label{cWWEq1}
\begin{split}
G^N: \R^N \times \R^N &\ra \MLo \times \MLo  \\
(\varTheta, \bar{\varTheta }) &\to   \Big(L^N(\varTheta), L^N(\bar{\varTheta }) \Big).
\end{split}
\end{align}
Set $ \gamma^N = (G^N)_\# \pi^N \in \cM_1(\MLo \times \MLo)$. 
Then $ \gamma^N $ has $ (L^N)_\#\mu^N $ and $ (L^N)_\# \nu^N $ as  marginals. Therefore, 
\begin{align}
\begin{split}
\bW^\L \Big( (L^N)_\#\mu^N\, &, \, (L^N)_\# \nu^N \Big)^2   \leq
\int_{(\MLo)^2} \WL(\s,\bar{\s})^2 \ d \gamma^N(\s,\bar{\s})\\
&= \int_{(\R^N)^2}  \frac{1}{N} \sum_{k=0	}^{N-1} W_2(\delta_{\t^k},\delta_{\bar{\t}^k})^2 \ d \pi^N(\varTheta,\bar{\varTheta})= \, 
\frac{1}{N} W_2(\mu^N, \nu^N)^2,
\end{split}
\end{align}
which concludes the proof.
\epr

\paragraph{McKean-Vlasov-functional.} 
Here we can even show a more        general  statement, which will be useful in the next chapter. 
\bl{LBMcK}
Let $ \mu^N \in \PwN  $ for all $ N \in \N $ and let $ A \in \cM_1(\cM_1(\ToR)) $.
Assume that $ ((K^N)_\# \mu^N)_N   $ converges weakly in $ \cM_1(\cM_1(\ToR)) $ to $ A $.
Then 
\be{LBMcKEq}
\liminf_{N\ra \infty }  \frac{1}{N} \cH^N(\mu^N) \geq  \int_{\cM_1(\ToR)}  \cF(\gamma)\,  d A(\gamma) .
\ee
\el
\bpr
Recall \eqref{WassGradEqq1}. Then we observe that 
\begin{align}\label{HNSplit}
\frac{1}{N} \cH^N(\mu^N) &= \frac{1}{N}\cH\left ( \mu^N\, \Bigg| \,  \exp\left(-\frac{1}{2} \sum_{k=0}^{N-1} \Psi(\t^k)\right) d\varTheta \right) \\
& \ \ \ +\frac{1}{2} \int_{  \cM_1(\TdR)	}\int_{  (\TdR)^2	} \left(\frac{1}{2} \big(\Psi(\t) +  \Psi(\bar{\t})\big) - J(x-\bar{x}) \t \bar{\t}			\right) d\gamma d\gamma\, d(K^N)_\# \mu^N (\gamma) .\nonumber
\end{align}
Similar arguments as in the proof of Theorem \ref{FProp} show that 
\begin{align}\label{HNSplit1}
\begin{split}
\liminf_{N\ra \infty } &\int_{  \cM_1(\TdR)	} \int_{  (\TdR)^2	} \left(\frac{1}{2} \big(\Psi(\t) +  \Psi(\bar{\t})\big) - J(x-\bar{x}) \t \bar{\t}			\right) d\gamma d\gamma\, d(K^N)_\# \mu^N (\gamma) \\
&\geq \int_{  \cM_1(\TdR)	}\int_{  (\TdR)^2	} \left(\frac{1}{2} \big(\Psi(\t) +  \Psi(\bar{\t})\big) - J(x-\bar{x}) \t \bar{\t}			\right) d\gamma d\gamma\, dA (\gamma).
\end{split}
\end{align}
It remains to show that 
\begin{align}\label{HNSplit2}
\begin{split}
\liminf_{N\ra \infty } \frac{1}{N} \cH\left ( \mu^N\, \Bigg| \,  \exp\left(-\frac{1}{2} \sum_{k=0}^{N-1} \Psi(\t^k)\right) d\varTheta \right) 
\geq \int_{  \cM_1(\TdR)	}\cH\left ( \gamma \, \Big| \,  e^{-\frac{1}{2} \Psi(\t)} d\t \right)  \, dA (\gamma).
\end{split}
\end{align}
Let $  \alpha := \int  e^{-\frac{1}{2} \Psi(\t)} d\t$ and for all $ n \in \N $, set 
\begin{align}
 B^N:= (K^N)_\# \left(\frac{1}{\alpha^N} \exp\left(-\frac{1}{2} \sum_{k=0}^{N-1} \Psi(\t^k)\right) d\varTheta  \right)  \ \ \ \text{and} \ \ \   A^N:= (K^N)_\# \mu^N.
\end{align}
Since the map $ K^N $ is injective, we have that 
\begin{align}\label{HNSplit33}
\cH\Big( A^N \, | \, &B^N \Big) =\cH\left ( \mu^N\, \Bigg| \,\frac{1}{\alpha^N}  \exp\left(-\frac{1}{2} \sum_{k=0}^{N-1} \Psi(\t^k)\right) d\varTheta \right)  .
\end{align}
It is an easy adaptation of Sanov's theorem that $ (B^N)_N $ satisfies a large deviation principle with rate function $ \cH\left ( \cdot \, \Big| \, \alpha^{-1} e^{-\frac{1}{2} \Psi(\t)} d\t \right) $; see e.g.\ \cite[Theorem 17]{pelren} for the details.
Therefore, \cite[3.5]{Mariani} implies that 
\begin{align}\label{HNSplit3}
\begin{split}
\liminf_{N\ra \infty } \frac{1}{N} \cH\left ( A^N\, \Big| \,  B^N \right) 
\geq \int_{  \cM_1(\TdR)	}\cH\left ( \gamma \, \Big| \, \alpha^{-1} e^{-\frac{1}{2} \Psi(\t)} d\t \right)  \, dA (\gamma).
\end{split}
\end{align}
\eqref{HNSplit3} and \eqref{HNSplit33} yield \eqref{HNSplit2}. This concludes the proof.
\epr

\paragraph{Initialization.}
\bl{LowerBoundInit}
Under the same assumptions as in Proposition \ref{LowerBound} and under \eqref{LBAssEq}, we have 
\be{LBIEq}
\liminf_{N\ra \infty }  \frac{1}{N} \left(  \cH(\nu_0^N	\, | \, \mu_0^N) - \frac{1}{2} \cH^N(\nu_0^N)  		\right)	\geq \cH(\nu_0\, | \, \mu_0) -\frac{1}{2} \cF(\nu_0)  .
\ee
\el
\bpr
Similarly as in the proof of Lemma \ref{Ilsc} and in Step 1 of the proof of Lemma \ref{LowerBoundAbsCont}, we observe that 
%
%
\begin{align}\label{LBIEqq1}
\begin{split}
&\liminf_{N\ra \infty }  \frac{1}{N} \big(  \cH(\nu_0^N	\, | \, \mu_0^N) -  \tfrac{1}{2}\cH^N(\nu_0^N)  		\big)	\geq  
\liminf_{N\ra \infty }   \frac{1}{N} \cH\Big( \nu_0^N\, \Big| \,  \exp\Big(-\frac{1}{2} \sum_{k=0}^{N-1} \Psi(\t^k)\Big) d\varTheta \Big)  \\
&  \quad   +\frac{1}{2}  \int_{  (\ToR)^2	} \Big(	\tfrac{1}{2} \Psi(\t) +\tfrac{1}{2} \Psi(\bar{\t}) 	+ J(x - \bar{x}) \t \bar{\t}	- 2\log\kappa(x,\t) - 2\log\kappa(\bar{x},\bar{\t})					 	\Big) \, d\nu_0 d\nu_0 ,
\end{split}
\end{align}
where we have used \eqref{LBInitBdd1} and \cite[5.1.7]{ambgigsav}.
Combining \eqref{LBIEqq1} with \eqref{HNSplit2} yields \eqref{LBIEq}. 
\epr

\paragraph{Metric derivative.}
Also here we can show directly a more general statement. 

\begin{lemma}\label{LowerBoundMetrDer}
Let $ (c_t)_{t\in[0,T]} \subset \cM_1(\PLo) $ be absolutely continuous with respect to the metric $ \bW^\L $ from Lemma \ref{LowerSemcWL}.
Let $  (\nu^N_t)_{t} \in \AC{[0,T]; \PwN} $ for all $ N \in \N $. 
Suppose that $ ((K^N)_\# \nu^N_t)_N $ converges to $ c_t $ weakly in $ \cM_1(  \cM_1(\ToR)  ) $ for all $ t \in[0,T] $.
Then, 
	\be{lscHDLEq}
	\liminf_{N\ra \infty} \frac{1}{N} \int_0^{T}   |(\nu^N)'|^2(t) \, dt \geq  \int_0^{T}   |c'|^2(t) \, dt .
	\ee
\end{lemma}
\bpr
Similarly as in \eqref{ILSCEq11} and in \eqref{LBAbsContEq0}, we obtain that for all $ \e \in (0, T/2) $
\begin{align}\label{lscHDLEq1}
	\int_0^{T-\e}   |c'|^2(t) \, dt
	&\leq \,   \liminf_{h\downarrow 0, h<\e }  \int_0^{T-\e} \frac{1}{h^2} \bW^\L(c_t, c_{t+h})^2   \, dt \nonumber\\
	& \leq \liminf_{h\downarrow 0, h<\e }\int_0^{T-\e} \liminf_{N\ra \infty } \frac{1}{h^2} \bW^\L \Big( (L^N)_\#\nu_{t}^N\, , \, (L^N)_\# \nu_{t+h}^N \Big)^2    dt\\
	&\leq \, 	\liminf_{N\ra \infty }   \frac{1}{N} \int_0^{T}   |(\nu^N)'|^2(r) \, dr. \nonumber 
\end{align}
Letting $ \e \downarrow 0 $ concludes the proof.
\epr

\paragraph{Metric slope.} 
Here we postpone the more general statement  to Chapter \ref{HDLChap}.

\bl{LowerBoundMetrSlope} 
Under the same assumptions as in Proposition \ref{LowerBound} and under \eqref{LBAssEq}, we have 
\be{LBMSEq}
\liminf_{N\ra \infty }  \frac{1}{N} \int_0^T |\partial \cH^N|^2(\nu^N_t) \, dt \geq  \int_0^T |\partial \cF|^2(\nu_t) \, dt.
\ee
\el
\bpr 
Similarly as in  Corollary  \ref{FSubDiffCor} one can show that (cf. \cite[4.3]{erbar} or  \cite[10.4.9]{ambgigsav})
\be{LBMSEq2}
|\partial \cH^N |(\nu_t^N) = \sup_{\varphi \in C_c^\infty  ( \R^N \, ; \, \R^N ), \,  \| \varphi \|_{\L^2( \nu_t^N)} >0} \frac{\left|\int_{\R^N} \left(\varphi \nabla H^N	-\mathrm{div} \varphi    \right) d\nu_t^N\right|}{\| \varphi \|_{\L^2( \nu_t^N)}}
\ee
for almost every $ t \in [0,T] $.
Let $ \varphi(\varTheta )=  (\beta (\tfrac{k}{N} \, , \, \t^k))_{k=0	}^{N-1} $ for some arbitrary $ \beta \in C_c^\infty( \ToR )  $ such that $ \| \beta \|_{\L^2( \nu_t)} >0 $.
This is admissible, since 
\be{LBMSEq01} 
 \| \varphi \|_{\L^2( \nu_t^N)}^2 =  N \int_{\cM_1(\ToR)}  \| \beta \|_{\L^2( \gamma)}^2 \,  d(K^N)_\#\nu^N_t (\gamma) 
\ee
and the right-hand side is greater than zero for $ N  $ large enough, since $ ((K^N)_\# \nu^N_t)_N $ converges to $  \delta_{\nu_t} $ weakly in $ \cM_1(  \cM_1(\ToR)  ) $. 
We obtain
\begin{align}\label{LBMSEq1}
	\liminf_{N\ra \infty }  &\frac{1}{N} |\partial \cH^N|^2(\nu^N_t)  \nonumber  \\
	&\geq 
	\liminf_{N\ra \infty }   \frac{ \left(\int_{\cM_1(\ToR)}\int \left(\beta \left[ \Psi' - 	 \int J(\cdot-\bar{x})\bar{\t}      \, d\gamma \right]	-\partial_\t \beta    \right) d\gamma\ d(K^N)_\#\nu^N_t (\gamma) \right)^2}{\int_{\cM_1(\ToR)} \int \beta^2 d\gamma \ d(K^N)_\#\nu^N_t (\gamma)} \\ 
	&= \frac{ 1} {\| \beta \|_{\L^2(\nu_t)}^2 } \left( \int \left(\beta \left[ \Psi' - 	 \int J(\cdot-\bar{x})\bar{\t}      \, d\nu_t \right]	-\partial_\t \beta    \right) d\nu_t \right)^2, \nonumber 
\end{align}
where we used in the last step a combination of \cite[5.1.7]{ambgigsav} and Step 4 of the proof of Lemma \ref{LowerBoundAbsCont}.
Taking the supremum over $ \beta  $ in \eqref{LBMSEq1}, we get via Corollary  \ref{FSubDiffCor}
\begin{align}\label{LBMSEq4}
	\liminf_{N\ra \infty }  \frac{1}{N} |\partial \cH^N|^2(\nu^N_t)  
	\geq  |\partial \cF|^2(\nu_t).
\end{align}
Finally, Fatou's lemma yields \eqref{LBMSEq}. 
\epr

\subsection{Recovery sequence}
\label{UpperBoundSec}

\bp{UpperBound} 
Let $ (\nu_t)_{t\in[0,T]}  \in \AC{[0,T] ; \PLo}$ be such that $ I[ (\nu_t)_t ] < \infty  $.
Then for all $ N \in \N $ there exists
$ (\nu_t^N)_{t\in[0,T]}  \in C([0,T] ; \cM_1(\R^N) )$ 
such  that 
$ ((K^N)_\# \nu^N_t)_N $ converges to $  \delta_{\nu_t} $ weakly in $ \cM_1(  \cM_1(\ToR)  ) $  for all $ t \in[0,T] $ and
\be{UBEq}
\limsup_{N\ra \infty }  \frac{1}{N} \left(  \frac{1}{2}\cJ^N[ (\nu_t^N)_{t} ]  + \cH(\nu_0^N	\, | \, \mu_0^N 	)	\right)	\leq  \frac{1}{2} \cJ[ (\nu_t)_{t} ] + \cH(\nu_0\, | \, \mu_0). 
\ee
\ep
\bpr
First, we observe that, since $ I[ (\nu_t)_t ] < \infty  $, we also have that 
\be{UBAssEq}
\cJ[ (\nu_t)_{t} ]  < \infty \qquad \text{and} \qquad 
\  \cH(\nu_0	\, | \, \mu_0 	)< \infty.
\ee
The recovery sequence will be given as follows. 
Recall the   partition $ (A_{k,N})_{k=0}^{N-1} $  of $ \To $ introduced in 
\eqref{AkN}.
Then define, for all $ N \in \N $ and for all $ t \in[0,T] $, $ \nu_t^N \in \cM_1(\R^N) $ by  
\be{RecSeqDef}
d\nu_t^N(\varTheta) = \prod_{k=0}^{N-1}  N \nu_t(A_{k,N}\times d\t^k ).
\ee
Lemma \ref{RecSeqConv} below shows that $ ((K^N)_\# \nu^N_t)_N $ converges to $  \delta_{\nu_t} $ weakly in $ \cM_1(  \cM_1(\ToR)  ) $  for all $ t $.
We show \eqref{UBEq} for each part separately. Hence, the claim
follows from the Lemmas \ref{RecSeqMc}--\ref{RecSeqM} and Lemma \ref{LBMcK}. 
\epr

\paragraph{Preliminaries.}
%
First we note that   \eqref{UBAssEq} implies that
 $ \nu_0 $ has a density $ f_0 $ with respect to $ \Leb{\ToR} $. 
Moreover, a similar computation as in \eqref{ILSCEq2} shows that 
\be{UBInitBd}
\cH(\nu_0	\, | \, \mu_0) - \frac{1}{2} \cF(\nu_0)  		 > -\infty . 
\ee
Together with Lemma \ref{LBF} and \eqref{UBAssEq}, this yields that
\be{UBMDMSFBd}
\int_0^T \Big(|\nu'|^2(t)	+ 	|\partial \cF|^2 (\nu_t)	\Big) \, dt  < \infty \qquad \text{and} \qquad 
\  \cF(\nu_T)< \infty. 
\ee
Since $ 	|\partial \cF| $ is a strong upper gradient  (\cite[1.2.1 and 2.4.10]{ambgigsav}), from \eqref{UBMDMSFBd} we infer that
\be{UBFBd1}
\int_0^T \cF(\nu_t) dt \leq  \int_0^T\int_t^T |\partial \cF|(\nu_r) \,|\nu'|(r) \,  dr dt + T \cF(\nu_T) < \infty .
\ee 
Therefore, for almost every $ t $,  $ \nu_t $ has a density $ f_t $ with respect to $ \Leb{\ToR} $.
And combining \eqref{UBFBd1} with  the lower bound on $ \cF $ (Lemma \ref{LBF}),  we infer  that
\be{UBt2moment}
\int_0^T \int_{\ToR} (|\t|^2+|\t|^{2\ell}) \, d\nu_t \, dt < \infty.
\ee
Finally, Lemma \ref{LBF} and \eqref{UBInitBd}  yield that 
$ \int (|\t|^2+|\t|^{2\ell}) d\nu_0 < \infty $.

\paragraph{Convergence.}

\bl{RecSeqConv}
Under the same setting as in the proof of Proposition \ref{UpperBound}, we have that for all $ t \in[0,T] $
\be{RecSeqConvEq}
(K^N)_\# \nu_t^N  \text{  converges to }  \delta_{\nu_t}   
\text{  weakly in $ \cM_1(\cM_1(\ToR)) $.}
\ee
\el
\bpr
For all $ N \in \N $ and $ t \in[0,T] $ let $ \Upsilon_t^N = (\vartheta^{k,N}_t)_{k=0, \dots, N-1} $ be a random variable with law $ \nu_t^N $ on a common probability space $ ( \Omega, \cF, \P )$. 

\noindent
\textbf{Step 1.} \big[ Let $ f \in C_b(\ToR) $, then $ \lim_{N \ra \infty } \int f \, d(K^N(\Upsilon_t^N)) = \int f d\nu_t $ a.s.  \big]

\noindent
The proof is a standard application of Kolmogorov's maximum inequality \cite[9.7.4]{dudley}. 
For the sake of completeness, we provide the details. 
Let $ \e > 0 $ and set for all $ N \in \N $
\be{RecSeqConvEq1}
\frac{1}{N}S_N :=  \frac{1}{N} \sum_{k=0}^{N-1} \left( f(\tfrac{k}{N}, \vartheta^{k,N}_t) - \Eb{f(\tfrac{k}{N}, \vartheta^{k,N}_t)}\right). 
\ee
Then
\begin{align}\label{RecSeqConvEq3} 
	\begin{split}
		\sum_{p\in \N} \Pb{\max_{2^{p-1} + 1 \leq n \leq  2^p} |S_n| > n  \e} 
		&\leq 
		\sum_{p\in \N} \Pb{ 		\max_{n\leq 2^p} |S_n| > 2^p \frac \e 2  }\leq 
		\sum_{p\in \N} \frac{4}{\e^2 2^{2p}}
		\mathrm{Var}[S_{2^p}] \\
		&\leq  
		\sum_{p\in \N} \frac{4}{\e^2 2^{p}}
		\| f \|_\infty^2 < \infty .
	\end{split}
\end{align}
Hence, the Borel-Cantelli Lemma yields that $ \lim_{N \ra \infty } S_N/N = 0 $ a.s. 
Finally, 
\begin{align}
	\begin{split}
		\int f \, dK^N(\Upsilon_t^N) 
		&= 
		\frac{1}{N}S_N+  \frac{1}{N} \sum_{k=0}^{N-1}  \Eb{f(\tfrac{k}{N}, \vartheta^{k,N}_t)} = \frac{1}{N}S_N+
		\sum_{k=0}^{N-1} \int_{A_{k,N}} \int_{  \R	}f(\tfrac{k}{N}, \vartheta) d\nu_t \\
		&
		= \frac{1}{N}S_N+
		\int_{  \ToR	}f\left(\tfrac{\lfloor xN \rfloor }{N}, \vartheta\right) d\nu_t (x,\vartheta) \ 
		\Ra  \ \int_{  \ToR	}f\, d\nu_t \quad \text{a.s.}
	\end{split}
\end{align}
\textbf{Step 2.} \big[ $ \Pb{ \lim_{N \ra \infty }\widetilde{W}(K^N(\Upsilon_t^N), \nu_t) = 0 } = 1 $.  \big]

\noindent
The claim follows from Step 1 once we apply the same arguments as in the proof of  \cite[11.4.1]{dudley}.
Recall that we have used those arguments already  to prove Lemma \ref{SeparableLem}.

\noindent
\textbf{Step 3.} \big[ $(K^N)_\# \nu_t^N  $ converges to $  \delta_{\nu_t}$   
weakly in $ \cM_1(\cM_1(\ToR)) $.  \big]

\noindent
Step 2 and \cite[9.2.1]{dudley}
yield that $ \lim_{N \ra \infty } \Pb{ \widetilde{W}(K^N(\Upsilon_t^N), \nu_t) >\e } = 0 $ for all $ \e > 0  $.
Hence, \cite[9.3.5]{dudley} implies the claim. This concludes the proof. 
\epr

\paragraph{McKean-Vlasov-functional.}
\bl{RecSeqMc}
Under the same setting as in the proof of Proposition \ref{UpperBound} and under \eqref{UBAssEq}, we have 
\be{RecSeqMcEq}
\limsup_{N\ra \infty }  \frac{1}{N} \,  \cH^N(\nu^N_t) \,  \leq \,\cF(\nu_t) \quad \text{for almsot every } t \in[0,T] .
\ee
In particular, \eqref{RecSeqMcEq} holds true for $ t=0 $ and $ t=T $.
\el

\bpr
Let $ t $ be such that $ \nu_t  $ has a density $ f_t $ and $ \int |\t|^2 d\nu_t < \infty $.
In particular, $ t=0 $ and $ t=T $ are admissible.
We observe that
\begin{align}\label{RecSeqMc1}
	\frac{1}{N} \,  \cH^N(\nu^N_t)
	&= 	
	\frac{1}{2}  \sum_{k,j=0	}^{N-1}  \int_{A_{k,N}}\int_{A_{j,N}}
	\int_{  \R 	} \int_{  \R 	}   J\left(	\tfrac{k-j}{N}	\right) \t \bar{\t}
	\, d\nu_t(\bar{x}, \bar{\t}) \,  d\nu_t(x, \t) \nonumber
	\\
	\begin{split}
		&\ \ \  + \frac{1}{N}  \sum_{k=0	}^{N-1} 
		\int_{  \R 	}   \log \left( \int_{A_{k,N}} 
		f_t(x, \t) e^{\Psi(\t)}  N dx
		\right) \,
		\int_{A_{k,N}} f_t(x, \t) Ndx \, d\t  \\
		&= 	
		\frac{1}{2}   \int_{(\ToR)^2 	}   J\left(	\tfrac{\lfloor xN \rfloor -\lfloor \bar{x}N \rfloor}{N}	\right) \t \bar{\t}
		\, d\nu_t(\bar{x}, \bar{\t}) \,  d\nu_t(x, \t) 
	\end{split}
	\\
	&\ \ \  +  
	\frac{1}{N}  \sum_{k=0	}^{N-1} 
	\int_{  \R 	}   \log \left( \int_{A_{k,N}} 
	f_t(x, \t) e^{\Psi(\t)}  Ndx
	\right) \,
	\int_{A_{k,N}} f_t(x, \t) e^{\Psi(\t)}  Ndx \, e^{-\Psi(\t)}  d\t. \nonumber
\end{align}
Since $ N \cdot  \Leb{A_{k,N}} $ is a probability measure and $ s \to s \log s  $ is convex on $ (0,\infty) $,  Jensen's inequality yields
\begin{align}\label{RecSeqMc2}
	\begin{split}
		&
		\frac{1}{N}  \sum_{k=0	}^{N-1} 
		\int_{  \R 	}   \log \left( \int_{A_{k,N}} 
		f_t(x, \t) e^{\Psi(\t)}  Ndx
		\right) \,
		\int_{A_{k,N}} f_t(x, \t) e^{\Psi(\t)}  Ndx \, e^{-\Psi(\t)}  d\t  \\
		&\ \ \leq 
		\sum_{k=0	}^{N-1} \int_{  \R 	}  \int_{A_{k,N}} \log \left(  
		f_t(x, \t) e^{\Psi(\t)}  
		\right) \,
		f_t(x, \t) \, dx\,  d\t = \int_{\ToR} \log \left(  
		f_t\,  e^{\Psi}  
		\right) d\nu_t.
	\end{split}
\end{align}
Moreover, the continuity of $ J $, the fact that $ \int |\t|^2 d\nu_t < \infty $ and the dominated convergence theorem yield
\begin{align}\label{RecSeqMc3}
	\lim_{N \ra \infty } \frac{1}{2}   \int_{(\ToR)^2 	}   J\left(	\tfrac{\lfloor xN \rfloor -\lfloor \bar{x}N \rfloor}{N}	\right) \t \bar{\t}
	\, d\nu_t \,  d\nu_t =
	\frac{1}{2}   \int_{(\ToR)^2 	}   J\left(x-\bar{x}\right) \t \bar{\t}
	\, d\nu_t\,  d\nu_t.
\end{align}
Lastly,  \eqref{RecSeqMc1},   \eqref{RecSeqMc2} and \eqref{RecSeqMc3} yield   \eqref{RecSeqMcEq}.
\epr

\paragraph{Initialization.}
\bl{RecSeqInit}
Under the same setting as in the proof of Proposition \ref{UpperBound} and under \eqref{UBAssEq}, we have  
\be{RecSeqInitEq}
\limsup_{N\ra \infty }  \frac{1}{N} \,\cH(\nu_0^N\, | \, \mu_0^N) \leq \, \cH(\nu_0\,| \, \mu_0)  \quad \text{for all } t \in[0,T] .
\ee
\el

\bpr
Set $ \rho_0(x,\t):= e^{-\Psi(\t)}\kappa(x,\t) $. Then, as in the proof of Lemma \ref{RecSeqMc}, we observe that
\begin{align}\label{UBIEq1}
 \frac{1}{N} \cH(\nu_0^N	&|  \mu_0^N)= 
\frac{1}{N}  \sum_{k=0	}^{N-1} 
\int_{  \R 	}   \log \left( \int_{A_{k,N}} 
f_0(x, \t) \rho_0(\tfrac{k}{N},\t)^{-1} N dx
\right) \,
\int_{A_{k,N}} f_0(x, \t)Ndx \, d\t \nonumber \\ 
&\leq \int_{  \ToR 	}   \log \left( 
f_0(x, \t) \rho_0(\tfrac{\lfloor xN\rfloor}{N},\t)^{-1}  
\right) \,
 f_0(x, \t)  \, dx \, d\t.
\end{align}
Under Assumption \ref{AssI}, we either have that $ \rho_0(\tfrac{\lfloor xN\rfloor}{N},\t) \geq c_\kappa'\,  e^{-(c_\kappa+1) \Psi(\t)  }  $ on the set $ \{\rho_0 > 0 \} $ or that $ \rho_0(\tfrac{\lfloor xN\rfloor}{N},\t) = \rho_0(x,\t) $ for all $(x,\t)\in \ToR $.
In the latter case, we trivially obtain \eqref{RecSeqInitEq}.
In the former case, the integrand on the right-hand side of \eqref{UBIEq1} is bounded from above by $ g:=\log 
(f_0\,\,  \exp((c_\kappa+1) \Psi)/c_\kappa')
f_0 $, which is integrable.
Indeed, from \eqref{UBInitBd} and \eqref{UBt2moment} we infer that $ \cH(\nu_0	\, | \, e^{- \Psi(\t)  }dxd\t) $ is finite.
This immediately implies  the integrability of $ g $.
Hence, we can apply the dominated convergence theorem to interchange the integral  and the limit, and the regularity assumptions on $ \rho_0 $ from Assumption \ref{AssI} lead to \eqref{RecSeqInitEq}.
\epr

\paragraph{Metric derivative.}

\bl{RecSeqMD}
Under the same setting as in the proof of Proposition \ref{UpperBound} and under \eqref{UBAssEq}, we have that for all $ N \in \N $ 
\be{RecSeqMDEq}
\frac{1}{\sqrt{N}} \,  |(\nu^N)'|(t) \,  \leq \, |\nu'|(t) \quad \text{for almost every } t \in[0,T] .
\ee
\el
\bpr
Let $ s < t  $. 
Let $ \pi \in \Opt{\nu_s, \nu_t} $ and define $ \gamma \in \cP_2(\R^N \times \R^N) $ by
\be{RecSeqMDEq1}
d\gamma(\varTheta, \bar{\varTheta }) = \prod_{k=0}^{N-1}  N \pi(A_{k,N}\times d\t^k \times d\bar{\t}^k ).
\ee
It is readily checked that $ \gamma \in \CouplW{\nu_s^N,\nu_t^N} $.
Therefore, 
\be{RecSeqMDEq2}
W_2(\nu_s^N,\nu_t^N)^2 
\leq \int_{\R^N \times \R^N} 
|\varTheta- \bar{\varTheta }|^2
d\gamma(\varTheta, \bar{\varTheta }) 
=N \WL(\nu_s, \nu_t)^2,
\ee
which immediately implies \eqref{RecSeqMDEq}.
\epr

\paragraph{Metric slope.}

\bl{RecSeqM}
Under the same setting as in the proof of Proposition \ref{UpperBound} and under \eqref{UBAssEq}, we have 
\be{RecSeqMEq}
\limsup_{N \ra \infty } \frac{1}{N}  \int_0^T   |\partial\cH^N|^2(\nu_t^N) \, dt \,  \leq  \,  \int_0^T   |\partial\cF|^2(\nu_t) \, dt .
\ee
\el

\bpr
Recalling the definition of the weak derivative, one can easily show that for all $ k \leq N-1 $
\be{RecSeqMEq1}
\partial_{\t} \int_{A_{k,N}} f_t(x,\t) \,dx =  \int_{A_{k,N}} \partial_{\t}f_t(x,\t) \,dx  \quad \text{
	for almost every $t $ and $ \t $.}
\ee
Let $ f_t^N  $ be the density of $ \nu_t^N $ with respect to $ \Leb{\R^N} $.
In view of \eqref{RecSeqMEq1}, we observe that
\begin{align}
	&\frac{1}{N}  \int_0^T  \int_{  \R^N	}
	\left| \frac{\nabla f_t^N(\varTheta )}{f_t^N(\varTheta )} + \nabla H^N(\varTheta )
	\right|^2 \, d\nu_t^N (\varTheta )\, dt \nonumber \\
	&= 
	\frac{1}{N} \sum_{k=0	}^{N-1}  \int_0^T  \int_{  \R^N	}
	\left( \frac{N\int_{A_{k,N}} \partial_{\t}f_t(x,\t^k) \,dx   }{N\int_{A_{k,N}} f_t(x,\t^k) \,dx   } + \Psi'(\t^k) - \frac{1}{2N}
	\sum_{j=0	}^{N-1} J\left(\tfrac{k-j	}{N}\right) \t^j 
	\right)^2 \, d\nu_t^N \, dt  \nonumber \\
	&= 
	\frac{1}{N} \sum_{k=0	}^{N-1}  \int_0^T   \int_{  \R^N	}
	\left( \frac{N\int_{A_{k,N}} \partial_{\t}f_t(x,\t^k) \,dx   }{N\int_{A_{k,N}} f_t(x,\t^k) \,dx   } + \Psi'(\t^k)
	\right)^2 \, d\nu_t^N \, dt \label{RecSeqMEq2} \\
	&\quad - 
	\frac{1}{N} \sum_{k=0	}^{N-1}  \int_0^T  \int_{  \R^N	}
	\left( \frac{N\int_{A_{k,N}} \partial_{\t}f_t(x,\t^k) \,dx   }{N\int_{A_{k,N}} f_t(x,\t^k) \,dx   } + \Psi'(\t^k)  \right)   \frac{1}{N}
	\sum_{j=0	}^{N-1} J\left(\tfrac{k-j	}{N}\right) \t^j 
	\, d\nu_t^N \, dt \label{RecSeqMEq3}\\
	&\quad + 
	\frac{1}{N} \sum_{k=0	}^{N-1}  \int_0^T  \int_{  \R^N	}
	\left( \frac{1}{2N}
	\sum_{j=0	}^{N-1} J\left(\tfrac{k-j	}{N}\right) \t^j 
	\right)^2 \, d\nu_t^N \, dt.  \label{RecSeqMEq4}
\end{align}
We treat each term \eqref{RecSeqMEq2}--\eqref{RecSeqMEq4} separately. 
First, we compute
\begin{align}\label{RecSeqMEq21}
	\text{\eqref{RecSeqMEq2}}
	= 
	\frac{1}{N} \sum_{k=0	}^{N-1}  \int_0^T  \int_{  \R	}
	\frac{\left(N\int_{A_{k,N}} \Big(\partial_{\t}f_t(x,\t) + \Psi'(\t) f_t(x,\t)\Big) \,dx  \right)^2 }{N\int_{A_{k,N}} f_t(x,\t) \,dx   } 
	\, d\t \, dt .
\end{align}
In the same way as in the proof of \cite[8.1.10]{ambgigsav}, we are allowed to apply Jensen's inequality for the integrand, since the function $ (x,z) \to x^2/z $ is convex on $ \R\times (0,\infty) $.
Hence,
\begin{align}\label{RecSeqMEq22}
	\begin{split}
		\text{\eqref{RecSeqMEq2}}
		&\leq 
		\sum_{k=0	}^{N-1}  \int_0^T   \int_{  \R	} \int_{A_{k,N}} 
		\frac{\Big(\partial_{\t}f_t(x,\t) + \Psi'(\t) f_t(x,\t)\Big)^2 }{ f_t(x,\t)   } 
		\,dx \, d\t   \, dt \\ 
		&= 
		\int_0^T  \int_{ \ToR	}
		\left(\frac{\partial_{\t}f_t(x,\t)    }{ f_t(x,\t)   }  + \Psi'(\t)\right)^2
		\, d\nu_t \, dt.
	\end{split}
\end{align}
Next, we observe that \eqref{RecSeqMEq3} is equal to
\begin{align}\label{RecSeqMEq31}
	&- 
	\sum_{k,j=0	}^{N-1}  \int_0^T  \int_{  \R^2	} \int_{A_{k,N}}\int_{A_{j,N}}
	\Big(  \partial_{\t}f_t(x,\t)   + \Psi'(\t)   f_t(x,\t) \Big)   
	J\left(\tfrac{k-j	}{N} \right) \bar{\t} f_t(\bar{x},\bar{\t})
	\, d\bar{x} \, dx \, d\bar{\t}\, d\t  dt \nonumber\\ 
	&=- 
	\int_0^T  \int_{  \ToR	} 
	\left(  \frac{\partial_{\t}f_t(x,\t)}{f_t(x,\t)}   + \Psi'(\t)    \right)   
	\int_{  \ToR	} J\left(\tfrac{\lfloor xN \rfloor- \lfloor \bar{x}N \rfloor}{N} \right) \bar{\t}
	\, d\nu_t(\bar{x},\bar{\t}) \, d\nu_t(x,\t) dt \\
	&\Ra 
	\int_0^T  \int_{  \ToR	} 
	\left(  \frac{\partial_{\t}f_t(x,\t)}{f_t(x,\t)}   + \Psi'(\t)    \right)   
	\int_{  \ToR	} J\left(x-\bar{x} \right) \bar{\t}
	\, d\nu_t(\bar{x},\bar{\t}) \, d\nu_t(x,\t) dt,\nonumber
\end{align}
where we have used the continuity of $ J $ and  the dominated convergence theorem, which is applicable, since by Young's inequality, \eqref{UBMDMSFBd} and \eqref{UBt2moment}
\begin{align}\label{RecSeqMEq32}
	\begin{split}
		&\left(  \frac{\partial_{\t}f_t(x,\t)}{f_t(x,\t)}   + \Psi'(\t)    \right)   
		\int_{  \ToR	} J\left(x-\bar{x} \right) \bar{\t}
		\, d\nu_t(\bar{x},\bar{\t}) \\
		&\qquad \leq \frac{1}{2} \left(  \frac{\partial_{\t}f_t(x,\t)}{f_t(x,\t)}   + \Psi'(\t)    \right)^2  + \frac{\|J\|_\infty}{2}\int_{  \ToR	}  \bar{\t}^2
		\, d\nu_t  \ \in \  \L^1([0,T] \times \ToR \, ;\, \nu_t dt).
	\end{split}
\end{align}
For the term \eqref{RecSeqMEq4}, we apply similar arguments to obtain that
\begin{align}\label{RecSeqMEq41}
	\text{\eqref{RecSeqMEq4}}&=  
	\frac{1}{4}\sum_{k,j,l=0	}^{N-1}  \int_0^T  \int_{  \R^3	} \int_{A_{k,N}}\int_{A_{j,N}}\int_{A_{l,N}}
	J\Big(\tfrac{k-j	}{N} \Big) \bar{\t}  \, J\Big(\tfrac{k-l	}{N}\Big)  \hat{\t}
	\, d\nu_t d\nu_t d\nu_t  dt + O\left( \frac{1}{N}\right) \nonumber \\ 
	&=
	\int_0^T  \int_{  \ToR	} 
	\left(   \frac{1}{2} 
	\int_{  \ToR	} J\left(\tfrac{\lfloor xN \rfloor- \lfloor \bar{x}N \rfloor}{N} \right) \bar{\t}
	\, d\nu_t(\bar{x},\bar{\t})  \right)^2 \, d\nu_t(x,\t) dt + O\left( \frac{1}{N}\right)\\
	&\Ra 
	\int_0^T  \int_{  \ToR	} 
	\left(   \frac{1}{2} 
	\int_{  \ToR	} J\left(x - \bar{x} \right) \bar{\t}
	\, d\nu_t(\bar{x},\bar{\t})  \right)^2 \, d\nu_t(x,\t) dt \nonumber.
\end{align}
Hence, \eqref{RecSeqMEq22}, \eqref{RecSeqMEq31}, \eqref{RecSeqMEq41} and Proposition \ref{SubdiffCharacF} show that 
\begin{align} \label{RecSeqMEq5}
	\limsup_{N \ra \infty } \frac{1}{N}  \int_0^T  \int_{  \R^N	}
	\left| \frac{\nabla f_t^N }{f_t^N } + \nabla H^N 
	\right|^2 \, d\nu_t^N  \, dt \leq \int_0^T   |\partial\cF|^2(\nu_t) \, dt .
\end{align}
Finally, it is known that (see for instance, \cite[10.4.9]{ambgigsav}) since the right-hand side (and hence also the left-hand side) of \eqref{RecSeqMEq5} is finite, we have 
\begin{align} \label{RecSeqMEq6}
	\limsup_{N \ra \infty } \frac{1}{N}  \int_0^T  \int_{  \R^N	}
	\left| \frac{\nabla f_t^N }{f_t^N } + \nabla H^N 
	\right|^2 \, d\nu_t^N  \, dt =   \limsup_{N \ra \infty } \frac{1}{N}  \int_0^T   |\partial\cH^N|^2(\nu_t^N) \, dt.
\end{align}
\eqref{RecSeqMEq5} and \eqref{RecSeqMEq6} conclude the proof.
\epr


\section{Hydrodynamic limit}
\label{HDLChap}

In this chapter we derive a law of large numbers for the system introduced in Section \ref{3.0.1}.

\bt{HDL}
For all $ N \in \N $, let $ (\mu_t^N)_{t\in[0,T]} $  be the Wasserstein gradient flow for 
$ \cH^N $ with initial value  $ \mu_0^N $. 
Assume either 
\begin{enumerate}[a)]
	\item Assumption \ref{AssI} on the sequence of initial data $ \{\mu_0^N\}_N $ and let $ \mu_0 = \rho_0 \,  \Leb{\ToR} $ with  $ \rho_0(x,\t) =\kappa(x,\t) \,  e^{-\Psi(\t)  } $, or
	\item $ \mu_0 \in D(\cF) $ and $ \{\mu_0^N\}_N $ is such that $ ((K^N)_\# \mu^N_0)_N $ converges to $  \delta_{\mu_0} $ weakly in $ \cM_1(  \cM_1(\ToR)  ) $   and $ \lim_{N\ra \infty } \frac{1}{N} \cH^N(\mu_0^N) = \cF(\mu_0) $. 
\end{enumerate}
Then, for all $ t \in [0,T] $, $ ((K^N)_\# \mu^N_t)_N $ converges to $  \delta_{\mu_t} $ weakly in $ \cM_1(  \cM_1(\ToR)  ) $, where $ \muc $ is the  gradient flow for $ \cF $ with initial value $ \mu_0 $ and 
\begin{align} \label{PropChaos}
\lim_{N\ra \infty } \frac{1}{N} \cH^N(\mu_t^N ) = \cF(\mu_t) \ \ \  \text{ for all $  t \in [0,T] $}.
\end{align}
Moreover, in the situation of b) and if $C_\Psi >0 $ and $ \ell \geq 2 $ in Assumption \ref{Fass}, then we even have that  	$ ((K^N)_\# \nu^N_t)_N $ converges to $  \delta_{\nu_t} $ weakly in $ \cM_1\big(\,  (\cP_2(\ToR), W_2)  \,\big) $ for all $ t \in [0,T] $. 
\et
\bpr
In the situation of \emph{a)}, the proof follows immediately from Theorem \ref{LDP}, since the corresponding rate function in the LDP  result  has a unique minimum at $ \muc $ by Theorem \ref{Lyapunov}.
So assume  \emph{b)}. 
First notice that
\be{HDLEq1}
\sup_{N\in \N}  \frac{1}{N} \int_0^T |(\mu^N)'|^2(t) \, dt , \ \  \sup_{N\in \N}  \frac{1}{N} \int_0^T |\partial \cH^N|^2(\mu^N_t) \, dt , \  \  \sup_{N\in \N}  \frac{1}{N} \cH^N(\mu_T^N)   \ \ < \infty ,
\ee
 since $ \cJ^N[(\mu_t^N)_t] = 0  $ for all $ N $,   $ \lim_{N\ra \infty } \frac{1}{N} \cH^N(\mu_0^N) = \cF(\mu_0) $ and by  Lemma \ref{WassGrad} \emph{(ii)}.
Moreover,  arguing as in \eqref{UBFBd1}, 
we infer that 
\be{HDLEq122}
\sup_{N\in \N}  \frac{1}{N}\int_0^T \cH^N(\mu^N_t) \, dt  \leq \sup_{N\in \N}  \frac{1}{N} \left( \int_0^T \int_t^T |\partial \cH^N|(\mu^N_r) \cdot |(\mu^N)'|(r) \,  dr dt  + T \cH^N(\mu^N_T) \right) < \infty ,
\ee 
and for all $ t \in [0,T] $,
\be{HDLEq12}
\sup_{N\in \N}  \frac{1}{N} \cH^N(\mu^N_t) \leq \sup_{N\in \N}  \frac{1}{N} \left( \int_t^T |\partial \cH^N|(\mu^N_r) \cdot |(\mu^N)'|(r) \,  dr  +  \cH^N(\mu^N_T) \right) < \infty .
\ee 
By Lemma \ref{WassGrad} \emph{(ii)} this implies that
\begin{align}\label{HDLEq13}
\begin{split}
\sup_{N\in \N}  \frac{1}{N} \int_0^T &\int_{\R^N}  \sum_{i=0}^{N-1}|\t^i|^{2\ell}\, d\mu^N_t(\varTheta) \, dt  < \infty , \text{ and }\\                                                                                                 
\sup_{N\in \N}  \frac{1}{N} &\int_{\R^N}  \sum_{i=0}^{N-1}|\t^i|^{2\ell}\, d\mu^N_t(\varTheta) \,  < \infty \ \ \ \text{ for  all  }\, t\in [0,T].
\end{split}
\end{align}

\noindent
\textbf{Step 1.} \big[ Compactness.   \big]

\noindent
Lemma \ref{CompactnessHDL}  yields the existence of a subsequence $ \{(\mu_t^n)_t\}_n $ and a continuous curve $ (c_t)_t \in C([0,T]\,;\,\cM_1(\cM_1(\ToR))) $ such that for all $ t \in [0,T] $, $ ((K^n)_\# \mu^n_t)_n $ converges to $  c_t $ weakly in $ \cM_1(  \cM_1(\ToR)  ) $, and if 
 $ \ell \geq 2 $, we even  have that this convergence holds  weakly in $ \cM_1\big(\,  (\cP_2(\ToR), W_2)  \,\big)$.

\noindent
\textbf{Step 2.} \big[ Superposition.  \big]

\noindent
Lemma \ref{SuperpositionHDL} below shows that there exists  a measure $ \varUpsilon \in \cM_1(\AC{[0,T]; \PLo}) $ such that 
 $(e_t)_\#\varUpsilon = c_t \text{ for all }t\in [0,T]$. 

\noindent
\textbf{Step 3.} \big[ Lower semi-continuity.  \big]

\noindent
Assumption b) and the Lemmas \ref{LBMcK}, \ref{LowerBoundMetrDer} and \ref{lscMHDL} show that 
\be{LSCHDL}
\int_{\AC{[0,T]; \PLo}} \Ind_{\mu_0}(\eta_0) \cdot  \cJ[\, (\eta_t)_t \,]\ d\varUpsilon((\eta_t)_t  ) \leq \liminf_{n \ra \infty  } \frac{1}{n} \cJ^n [\, (\mu_t^n)_t \, ]. 
\ee

\noindent
\textbf{Step 4.} \big[ Convergence towards $  \delta_{\mu_t} $ for all $ t \in [0,T] $.  \big]

\noindent
Step 3 shows that 
\be{ConclHDL}
\int_{\AC{[0,T]; \PLo}} \Ind_{\mu_0}(\eta_0) \cdot  \cJ[\, (\eta_t)_t \,]\, d\varUpsilon( (\eta_t)_t  ) \leq 0.
\ee
Since the integrand on the left-hand side is non-negative (see Theorem \ref{Lyapunov}), this implies that $ \Ind_{\mu_0}(\eta_0) \cdot  \cJ[\, (\eta_t)_t \,] = 0 $ 
for $ \varUpsilon $-a.e.\  $ (\eta_t)_t $. 
However, the uniqueness claim in Theorem \ref{Lyapunov} yields that  $ \varUpsilon $ must be concentrated on $ \muc$. 
Together with Step 1 and Step 2, this shows that $ (K^n)_\# \mu_t^n $ converges to $ \delta_{\mu_t} $  weakly in $ \cM_1(  \cM_1(\ToR)  ) $, respectively weakly in $ \cM_1\big(\,  (\cP_2(\ToR), W_2)  \,\big)$, for all $ t \in [0,T] $.
Since the limit is unique, we also get the convergence of the full sequence.

\noindent
\textbf{Step 5.} \big[ Proof of \eqref{PropChaos}.  \big]

\noindent
From the previous steps, we infer that $ \lim_{N\ra \infty} \frac{1}{N} \cJ^N [\, (\mu_t^N)_t \, ] = \cJ[\, (\mu_t)_t \,]  $.
Using again the Lemmas \ref{LBMcK}, \ref{LowerBoundMetrDer} and \ref{lscMHDL} and Assumption \emph{b)}, this easily implies that $ \lim_{N\ra \infty } \frac{1}{N} \cH^N(\mu_T^N) = \cF(\mu_T) $. 
We can now replace $ T $ by some arbitrary $ t \in (0,T) $ and repeat the above proof to obtain
  \eqref{PropChaos}. 
\epr

\begin{lemma}[Compactness]\label{CompactnessHDL}
Let $ (\mu_t^N)_t \in \AC{[0.T] \, ; \, \PwN} $ for all $ N \in \N $. Assume that 
\be{HDLCompactnessEq1}
\sup_{N\in \N}  \frac{1}{N} \int_0^T |(\mu^N)'|^2(t) \, dt  < \infty \quad \text{ and } \quad \sup_{N\in \N}  \frac{1}{N} \int_{\R^N}  \sum_{i=0}^{N-1}|\t^i|^{2\ell}\, d\mu^N_t(\varTheta)  < \infty \ \  \forall\, t\in [0,T].
\ee
Then there exists a subsequence $ \{(\mu_t^n)_t\}_n $ and a curve $ (c_t)_t \in C([0,T];\cM_1(\cM_1(\ToR))) $ such that
for all $ t \in [0,T] $, $ ((K^n)_\# \mu^n_t)_n $ converges to $  c_t $ weakly in $ \cM_1(  \cM_1(\ToR)  ) $.
If 
$ \ell \geq 2 $, we even  have that $ ((K^n)_\# \mu^n_t)_n $ converges to $  c_t $ weakly in $ \cM_1\big(\,  (\cP_2(\ToR), W_2)  \,\big)$.
 
\end{lemma}
\bpr
In view of Lemma \ref{KvsL}, it is equivalent to show the claim with $ L^N $ replacing $ K^N $. 
Recall the definitions of $ \widetilde{\bW} $ (Lemma \ref{KvsL}) and $ \bW^\L $ (Lemma \ref{LowerSemcWL}). 
Let $ \widetilde{\bW} _2   $ denote the Wasserstein distance on $ \cM_1\big(\,  (\cP_2(\ToR), W_2)  \,\big)$ induced by the distance $ W_2/(1+W_2) $.
Then $ \widetilde{\bW} _2   $ metrizes the weak convergence in $ \cM_1\big(\,  (\cP_2(\ToR), W_2)  \,\big)$.
In Lemma \ref{cWW} we have seen  that
\be{CompactnessHDLEq}
\bW^\L \Big( (L^N)_\#\mu_s^N\, , \, (L^N)_\#\mu_t^N \Big)   
\leq \frac{1}{\sqrt{N}} W_2(\mu_s^N, \mu_t^N) \quad \text{ for all $ 0\leq  s < t \leq T $}.
\ee
By \eqref{HDLCompactnessEq1} and since both $ \widetilde{\bW}_2   $ and $ \widetilde{\bW} $ are dominated by $ \bW^\L $, this implies the equi-continuity of the sequence $ \{((L^N)_\#\mu_t^N)_t\}_N  $ with respect to $ \widetilde{\bW}_2   $ and $ \widetilde{\bW} $. 
Moreover, again by  \eqref{HDLCompactnessEq1}, we have that for all $ t \in [0,T] $
\be{CompactnessHDLEq2}
\sup_{N\in \N}  \int_{\cM_1(\ToR)}  \int_{\ToR} |\t|^{2\ell}\, d\gamma \ d(L^N)_\#\mu_t^N(\gamma)
= \sup_{N\in \N}  \frac{1}{N} \int_{\R^N}  \sum_{i=0}^{N-1}|\t^i|^{2\ell}\, d\mu^N_t(\varTheta)  < \infty.
\ee
This shows that $ ((L^N)_\#\mu_t^N)_N  $ is relatively compact in $ \cM_1(\cM_1(\ToR)) $ for all $ t \in [0,T] $ with respect to $ \widetilde{\bW} $ and with respect to $ \widetilde{\bW}_2   $ if $ \ell \geq 2 $ (cf. \cite[6.8 (iii)]{vil}). 
Thus, we can apply  the (extended) Arzel\'a-Ascoli theorem  (\cite[Chapter 7, Theorem 17]{Kelley}) to conclude the proof. 
\epr

\vspace{-0.3cm}
\begin{lemma}[Superposition]\label{SuperpositionHDL}
Consider the same setting as in Lemma \ref{CompactnessHDL} and assume in addition that 
\be{SuperpositionHDLEq00}
\sup_{N\in \N}  \frac{1}{N} \int_0^T \int_{\R^N} |\varTheta|^{2}\, d\mu^N_t(\varTheta) \, dt  < \infty .
\ee
Then $ (c_t)_t $ is  absolutely continuous with respect to $  \bW^\L  $, and  there exists  a measure $ \varUpsilon \in \cM_1(\AC{[0,T]; \PLo}) $ such that 
	\be{Super}
	(e_t)_\#\varUpsilon = c_t \text{ for all }t\in [0,T]  \ \text{ and } \ \int |\eta'|^2(t)\,	d\varUpsilon((\eta_t)_t ) = |c'|^2(t)  \text{ for a.e.\ }t\in [0,T].
	\ee
\end{lemma}
\bpr 
Note that $ \MLo $ is a closed subspace of $ \cM_1(\ToR) $.
Therefore, the Portmanteau theorem (\cite[11.1.1]{dudley}) yields that for  almost every  $ t\in [0,T] $
\begin{align}
	c_t ( \MLo) \geq \limsup_{n\ra \infty} \, (L^n)_\#\mu_t^n( \MLo) = 1.
\end{align}
Hence,  $ c_t $ is supported in  $ \MLo $ for  almost every  $ t $. 
Moreover, \eqref{HDLCompactnessEq1} and \cite[5.1.7]{ambgigsav} show that 
$ c_t $ is supported in  $ \PLo $ for all $ t $.
To show the absolute continuity of $ (c_t)_t $  with respect to $ \bW^\L $
we proceed as in  the proofs of the  Lemmas \ref{Ilsc} and  \ref{LowerBoundAbsCont}.
We have that
\begin{align}\label{SuperpositionHDL2}
\sup_{0<h<T} \int_0^{T-h} \frac{1}{h^2} \bW^\L(c_t, c_{t+h})^2  dt                                                                                                                                                  
&\leq   \sup_{0<h<T} \int_0^{T-h} \liminf_{n\ra \infty } \frac{1}{h^2} \bW^\L \Big( (L^n)_\#\mu_{t}^n\, , \, (L^n)_\# \mu_{t+h}^n \Big)^2    dt \nonumber  \\
&\leq 				\liminf_{n\ra \infty }   \frac{1}{n} \int_0^{T}   |(\mu^n)'|^2(r) \, dr < \infty .
\end{align}
Moreover, by \eqref{SuperpositionHDLEq00},
\begin{align}\label{SuperpositionHDL3}
\int_0^{T}  \bW^\L(c_t, \delta_{\delta_0\otimes \Leb{\To}})^2  dt 
&=
\int_0^{T}  \int_{\PLo} \WL(\gamma, \delta_0\otimes \Leb{\To})^2 \,  dc_t(\gamma)  dt 
=
\int_0^{T}  \int \int |\t|^2 d\gamma \,dc_t(\gamma) dt \nonumber\\
&\leq 
\liminf_{n\ra \infty }  \int_0^{T}  \int_{\PLo} \int_{\ToR} |\t|^2 d\gamma \,d(L^n)_\#\mu_t^n(\gamma) \, dt\\
&=
\liminf_{n\ra \infty } \frac{1}{n} \int_0^T \int_{\R^n} |\varTheta|^{2 }\, d\mu^n_t(\varTheta) \, dt  < \infty.\nonumber
\end{align}
By  \cite[Lemma 1]{Lisini}, \eqref{SuperpositionHDL2} and \eqref{SuperpositionHDL3} yield the absolute continuity of $ (c_t)_t $.
Finally, \cite[Theorem 5]{Lisini} shows that this already implies the second claim. 
\epr

%
%

\vspace{-0.3cm} 
\begin{lemma}[Lower semi-continuity, metric slope]\label{lscMHDL}
	Let $ \mu^n \in \Pwn \cap D(\cH^n) $ for all $ n \in \N $.
	Assume that  $ (K^n)_\# \mu^n \hra c $ for some $ c \in \cM_1(\cM_1(\ToR)) $. 
	Then
	\be{lscMHDLEq}
	\liminf_{n\ra \infty} \frac{1}{n}|\partial \cH^n|^2 (\mu^n)  \geq \int_{\cM_1(\ToR)} |\partial \cF|^2(\s) \, dc(\s).
	\ee
\end{lemma}
\bpr
We use the same strategy as in the proof of \cite[3.5]{Mariani}.
As in \cite[3.9]{Mariani}, let $ \{	(E_{\delta,l}^i)_{i=0}^{N_{\delta,l}}  \}_{\delta>0, l \in \N} $ be a sequence of subsets of $ \cM_1(\ToR) $ such that 
\begin{enumerate}[1)]
  	\setlength\itemsep{-0.3em}
	\item $  \lim_{l \ra \infty } c\left(\cup_{i=1}^{N_{\delta,l}} E_{\delta,l}^i \right) = 1 $ and  $ \cup_{i=0}^{N_{\delta,l}} E_{\delta,l}^i = \cM_1(\ToR)			$ , 
		\item $ E_{\delta,l}^i \cap E_{\delta,l}^j = \emptyset $ if $ j \neq i $,
		\item $ \widetilde{W}(\s, \eta ) < \delta $ for all $ \s, \eta \in E_{\delta,l}^i $ and $ i=1, \dots, N_{\delta,l} $,
		\item $ c( \partial E_{\delta,l}^i ) = 0  $ for all  $ i=1, \dots, N_{\delta,l} $, 
		\item each $ E_{\delta,l}^i $ has non-empty interior, 
		\item $ (E_{\delta,l}^i)_{i=0}^{N_{\delta,l}} $ is finer than $ (E_{\delta',l'}^i)_{i=0}^{N_{\delta',l'}} $ if $ \delta \leq \delta' $ and $ l \geq l' $.
\end{enumerate}
For the proof of the existence of such a sequence, we refer to \cite[3.9]{Mariani}.
Assume that the left-hand side of \eqref{lscMHDLEq} is finite, since the claim would be trivial otherwise.
Let $ (\mu^m)_m $ be a subsequence such that 
\be{lscMHDLEq2}
	\lim_{m\ra \infty} \frac{1}{m}|\partial \cH^m|^2 (\mu^m)=	\liminf_{n\ra \infty} \frac{1}{n}|\partial \cH^n|^2 (\mu^n) \quad \text{and} \quad
	\sup_{m\in \N} \frac{1}{m}|\partial \cH^m|^2 (\mu^m) < \infty.
\ee 
In particular, by \cite[10.4.9]{ambgigsav}, this implies that 
\be{lscMHDLEq3}
|\partial \cH^m|^2 (\mu^m) = \int_{\R^m} 	\left| \frac{\nabla \rho^m}{\rho^m} + \nabla H^m \right|^2		\, d\mu^m,
\ee 
where for all $ m $,  $ \rho^m $ denotes the  density of $ \mu^m $ with respect to $ \Leb{\R^m} $.
For each $ m,\delta,l,i $, define the measure $ \mu^{m,\delta,l,i} \in \Pwn $ by 
\be{lscMHDLEq1}
\int_{\R^m} f\, d\mu^{m,\delta,l,i}  = \frac{1}{(K^m)_\# \mu^m(E_{\delta,l}^i)} \int_{(K^m)^{-1}(E_{\delta,l}^i)} f\, d\mu^{m}
\ee
for all measurable and bounded $ f:\R^N \ra \R $.
Then
\begin{align}
	\lim_{m\ra \infty} \frac{1}{m}|\partial \cH^m|^2 (\mu^m)&=
		\lim_{m\ra \infty} \sum_{i=0}^{N_{\delta,l}} \frac{1}{m} \int_{(K^m)^{-1}(E_{\delta,l}^i)} 	\left| \frac{\nabla \rho^{m,\delta,l,i}}{\rho^{m,\delta,l,i}} + \nabla H^m \right|^2		\, d\mu^{m,\delta,l,i}  \cdot (K^m)_\# \mu^m(E_{\delta,l}^i) \nonumber \\
		&=  \sum_{i=0}^{N_{\delta,l}} \lim_{m\ra \infty} \frac{1}{m} |\partial \cH^m|^2 (\mu^{m,\delta,l,i}) \cdot c(  E_{\delta,l}^i ),
\end{align}
where we have used property 4).
If we define a piecewise constant function $ I_{\delta,l} $ by
\begin{align}
I_{\delta,l}(\gamma) =  \lim_{m\ra \infty} \frac{1}{m} |\partial \cH^m|^2 (\mu^{m,\delta,l,i}),  \quad \text{ if } \gamma \in  E_{\delta,l}^i 
\end{align}
and use Fatou's Lemma, we obtain that
\begin{align}
\lim_{m\ra \infty} \frac{1}{m}|\partial \cH^m|^2 (\mu^m)\geq \int_{\cM_1(\ToR)} \liminf_{l\ra\infty} \liminf_{\delta\ra 0} I_{\delta,l}(\gamma)  \, dc( \gamma).
\end{align}
By a straightforward modification of the proof of Lemma \ref{LowerBoundMetrSlope}, we can show that for $ c $-a.e.\ $ \gamma $
\begin{align}
\liminf_{l\ra\infty} \liminf_{\delta\ra 0} I_{\delta,l}(\gamma)  \, \geq |\partial \cF|^2( \gamma).
\end{align}
This concludes the proof.
\epr



\paragraph*{Acknowledgement}
The first author thanks  Matthias Erbar, Max Fathi, Dmitry Ioffe and Andr\'e Schlichting for numerous useful discussions. 
Special thanks  to Lorenzo Dello Schiavo for providing many good ideas and proofreading a lot of parts of this work. 
Moreover, we would like to thank the anonymous referees for reading the paper with
great care and for their valuable comments.

\renewcommand{\baselinestretch}{0.85}\normalsize
\bibliography{TEX-Bib}{} 
\addcontentsline{toc}{section}{References}

\end{document}